\documentclass[12pt]{amsart}
\usepackage{amsmath,amssymb,amsthm}
\usepackage{ifthen,verbatim,here}

\usepackage{mathrsfs}
\usepackage{yhmath}
\usepackage{booktabs}
\usepackage{longtable}
\usepackage{color}
\usepackage{url}
\usepackage[english]{babel}
\usepackage{verbatim,here}
\usepackage[utf8]{inputenc}
\usepackage[T1]{fontenc}
\usepackage{floatflt,graphicx,graphics}
\usepackage{a4wide}

\numberwithin{equation}{section}
\nonstopmode
\theoremstyle{plain}

\newtheorem{theorem}{Theorem}[section]

\theoremstyle{definition}
\newtheorem{remark}{Remark}[section]
\newtheorem{example}{Example}[section]

{\qed\bigskip}

\newcounter{alphabet}

\newcommand{\be}{\begin{eqnarray}}
\newcommand{\ee}{\end{eqnarray}}
\newcommand{\ba}{\begin{array}}
\newcommand{\ea}{\end{array}}
\newcommand{\ben}{\begin{eqnarray*}}
\newcommand{\een}{\end{eqnarray*}}

\newcommand{\C}{{\mathbb C}}
\newcommand{\R}{{\mathbb R}}

\newcommand{\B}{\mathbb{B}}

\newcommand{\capa}{\operatorname{cap}}

\newcommand{\moda}{\operatorname{mod}}

\newcommand{\D}{{\Delta}}
\newcommand{\E}{{\mathcal E}}

\newcommand{\K}{\mathcal{K}}

\renewcommand{\Im}{{\operatorname{Im}}}
\renewcommand{\Re}{{\operatorname{Re}}}

\newcommand {\M} {\mathsf{M}}

\renewcommand{\i}{\mathrm{i}}

\newcommand{\bI}{{\bf I}}
\newcommand{\bM}{{\bf M}}
\newcommand{\bN}{{\bf N}}
\newcommand{\CC}{{\mathbb C}}

\newcommand{\col}{\colon\,}

\newcommand{\g}{\gamma} 

\font\fFt=eusm10 
\font\fFa=eusm7  
\font\fFp=eusm5  
\def\K{\mathchoice
{\hbox{\,\fFt K}}
{\hbox{\,\fFt K}}
{\hbox{\,\fFa K}}
{\hbox{\,\fFp K}}}
\font\fFt=eusm10 
\font\fFa=eusm7  
\font\fFp=eusm5  
\def\E{\mathchoice
{\hbox{\,\fFt E}}
{\hbox{\,\fFt E}}
{\hbox{\,\fFa E}}
{\hbox{\,\fFp E}}}

\newcounter{minutes}
\divide\time by 60
\newcounter{hours}
\multiply\time by 60
\addtocounter{minutes}{-\time}
\begin{document}

\bibliographystyle{amsplain}
\title[On computation of capacities and conformal invariants]
{On computation of capacities and conformal invariants}

\def\thefootnote{}
\footnotetext{
\texttt{\tiny File:~\jobname .tex,
          printed: \number\year-\number\month-\number\day,
          \thehours.\ifnum\theminutes<10{0}\fi\theminutes}
}
\makeatletter\def\thefootnote{\@arabic\c@footnote}\makeatother

\author[M.M.S. Nasser]{Mohamed M. S. Nasser}
\address{Wichita State University \\
Department of Mathematics, Statistics, and Physics \\
Wichita, KS 67260-0033, USA}
\email{mms.nasser@wichita.edu}

\author{Matti Vuorinen}
\address{
University of Turku\\
Department of Mathematics and Statistics\\
FI--20014 University of Turku, Finland\\  ORCID: 0000-0002-1734-8228
}
\email{vuorinen@utu.fi}

\begin{abstract}
We give a survey of computation of the conformal capacity of planar condensers, generalized capacity, and logarithmic capacity with emphasis on our recent work 2020-2025.
We also discuss some applications of our method based on the boundary integral equation with the generalized Neumann kernel 
to the computation of several  other conformal invariants: harmonic measure, modulus of a quadrilateral, reduced modulus, hyperbolic capacity, and elliptic capacity.  
Here the solution of mixed Dirichlet-Neumann boundary value problem for the Laplace equation has a key role. At the end of the paper we give a topicwise structured list  to our  extensive bibliography on
constructive complex analysis and potential theory.

\end{abstract}

\maketitle

\def\cprime{$'$} \def\cprime{$'$} \def\cprime{$'$}
\providecommand{\bysame}{\leavevmode\hbox to3em{\hrulefill}\thinspace}
\providecommand{\MR}{\relax\ifhmode\unskip\space\fi MR }
\providecommand{\MRhref}[2]{%
  \href{http://www.ams.org/mathscinet-getitem?mr=#1}{#2}
}
\providecommand{\href}[2]{#2}

\tableofcontents

\section{Introduction}\label{sec:int}

Conformal capacity of condensers is one of the key notions of potential theory with many 
applications to geometric function theory \cite{gor,hkm,hkv,gmp}, to PDEs, to conformal invariants,
moduli of curve families. Thus the roots of the study of conformal capacity are in classical
function theory, in the works of Koebe, Bergman, Gr\"otzsch, Teichm\"uller, Ahlfors,  Beurling,  Fuglede  and many others \cite{ah, ab, f,o,tsu, ki}.

The book of  P\'{o}lya  and  Szeg\"{o} \cite{ps} was a landmark work with many results on isoperimetric problems of mathematical physics. Their topic of interest was to study extremal problems formulated in terms of domain characteristics, 
what  they called
domain functionals, like perimeter, center of mass, torsion constant, principal frequency, and, in particular, the capacity of condensers. The classical isoperimetric problem is to maximize the area of a planar domain with a given perimeter, here the domain functionals are area and perimeter. The authors of \cite{ps} studied many extremal problems for pairs of domain functionals, the problem was to minimize/maximize one of the domain functionals under the constraint that the other domain functional was constant. The results were then summarized  in numerous tables with the relevant numerical data \cite{ps}.  
It is often the case that the extremal configurations
exhibit symmetry and therefore symmetrization methods compatible with domain functionals
can provide guides for the solution of constrained extremal problems expressed
in terms of domain functionals. In the case of conformally invariant domain
functionals, such as the conformal capacity of a condenser, it is convenient
to simplify the problem by means of auxiliary conformal mapping, if possible.
Both methods, symmetrization and use of auxiliary conformal mappings are standard
tools in the study of conformal capacity and its applications to geometric function theory. For purely
numerical approximation of the conformal capacity
of a condenser, the problem is reduced to the Dirichlet
problem for the Laplace operator. The impact of computers on constructive complex analysis is apparent
in the following collections of papers \cite{beck,to,tr,pasa} published in the second half of the 20'th century, since the publication of \cite{ps}. The aforementioned fundamental work
was continued by the prominent early pioneers of numerical conformal mappings,  Gaier \cite{gai}, \cite{p1} and  Henrici \cite{hen}. For the relevant literature, see Section~\ref{sec:topic}.

The advent of PCs was a quantum leap in the development of computational methods.  
Numerical work with methods such as the Schwarz-Christoffel formula became much easier and faster, see  Trefethen \cite{tr},  Driscoll-Trefethen \cite{dt},  DeLillo, Elcrat,  Kropf, and  Pfaltzgraff~\cite{dekp}.
Classical formulas for conformal mappings are usually expressed in terms
of special functions like elliptic functions and elliptic integrals which are
not practical for manual calculations \cite{as, avv,bb,law}.   Such calculations became practical  when software packages like MATLAB  were introduced  and  suitable software was available.  
Methods such as the fast Fourier transform (FFT) \cite{fft} and the fast multipole method (FMM)~\cite{gro, gg} reduced the computational time significantly and the theory and practice of PDE solvers developed fast.
All this progress marked the beginning of the modern era for constructive methods and, in particular, for 
numerical conformal mapping and capacity computation. 
Some of the surveys are \cite{weg05, gut, bds, dt}, \cite[pp. 13-16]{ps10}, \cite[pp.3-12]{ky}. Learning material includes the textbook Crowdy \cite{c} and the
lecture notes with solved problems Papamichael \cite{p2}.

A condenser is a pair  $(G,E)$  where $G \subset \mathbb{R}^n, n\ge 2,$ is a domain and $E \subset G$ is a nonempty compact set.
Condenser capacity is an important tool in geometric function theory \cite{du,gor, gmp,hkm,hkv}.  
The capacity we study is the \emph{conformal capacity of a condenser}. 
The conformal capacity of a condenser  $(G,E)$  is defined as \cite{hkv}
\begin{align}\label{def_condensercap-1}
{\rm cap}(G,E)=\inf_{u\in A}\int_{G}|\nabla u|^n \, dm,
\end{align}
where $A$ is the class of $C^{\infty}_0(G)$ functions $u: G \to [0,\infty) $
such that $u(x)\geq1$ for all $x \in E$ and $dm$ is the $n$-dimensional Lebesgue measure. 
Below usually $n=2.$
This capacity is also related to the modulus $\M(\Delta)$ of the family $\Delta$
of all curves in $G \setminus E$ joining the set $E$ with the boundary $\partial G$ \cite{vais}, \cite[Thm 9.6]{hkv} as follows
\begin{equation}
\label{capmod}
{\rm cap}(G,E) =\M(\Delta) \,.
\end{equation}
The modulus of a curve family was introduced by Ahlfors and  Beurling in their landmark paper~\cite{ab} in 1950 for the case $n=2$ and it has become a key tool
in geometric function theory. Fuglede~\cite{f} then extended these notions
for $n \ge 3.$  Both the modulus and
the conformal capacity are conformally invariant. Considering the many applications of these notions 
Gol'dshte\v{i}n-Reshetnyak \cite{gor}, Heinonen-Kilpel\"ainen-Martio \cite{hkm}, Hariri-Kl\'en-Vuorinen \cite{hkv}, it is
surprising that the exact values of condenser
capacities are only known in very few cases, and
therefore it is natural to look for methods for
numerical approximation.

As mentioned above, this paper is a survey of some of our recent work during the past five years on numerical approximation of conformal invariants, based on the boundary integral equation method
developed by the first author during the past two decades Nasser ~\cite{Nas-ETNA}. We apply this method  to
investigate conformal invariants that the second author has studied in   Hakula-Rasila-Vuorinen \cite{hrv1,hrv2,hrv3} and also in
\cite{avv,hkv}. The topic of this paper relies on the pioneering work of mathematicians of earlier generations, researchers cited above, our teachers, and colleagues.  At the end of this paper
we give  a topicwise organized list  to literature which we hope will be useful to the reader. The
structure of this paper appears in the list of contents and there is no need to repeat it here. 
What is noteworthy is the remarkable feature of the boundary integral equation method: the same method works for the numerical approximation of a very large class of conformal invariants, usually with
rather small changes in coding. In the cases we tested, the accuracy is  the same as in the case of other methods. The computational
results produced for this paper are given in numerous figures and some numerical
tables. These may deal with the rate of convergence of the
iterations, error estimates of the computations, or just indicate the values computed.

\section{Preliminary results}

The conformal capacity of a condenser is conformally invariant. It is therefore
natural to take this invariance into account when analyzing the values of the capacity.
The hyperbolic geometry is well-suited for the purpose because of its conformal invariant
character.

\subsection{Hyperbolic geometry}\label{hg}
We recall some basic formulas and notation for hyperbolic geometry from~\cite{b}.
The Euclidean balls with center $x \in {\mathbb{R}^n}$ and radius $r>0$ are denoted $B^n(x,r)$ and its boundary sphere is $S^{n-1}(x,r)$. For brevity we write $ {\mathbb{B}^n}= B^n(0,1)$.

For $a,b\in \mathbb{B}^2$, the hyperbolic distance $\rho_{\B^2}(a,b)$ between $a$ and $b$ is defined via the formula \cite{b}
\begin{equation}\label{rhoB}
\sinh \frac{\rho_{\B^2}(a,b)}{2}=
\frac{|a-b|}{\sqrt{(1-|a|^2)(1-|b|^2)}}.
\end{equation}
The hyperbolic disk with center $x \in {\mathbb{B}^2} $ and radius $R>0$
is $B_{\rho}(x,R)=\{ z: \rho_{\B^2}(x,z)<R \}.$ We often use the
connection between the hyperbolic disk and Euclidean disk \cite[p. 56, (4.20)]{hkv}
\begin{equation}\label{hkv4.20}
\begin{cases}
B_\rho(x,R)=B^2(y,r)\;,&\\
\noalign{\vskip5pt}
{\displaystyle y=\frac{x(1-t^2) }{ 1-|x|^2t^2}\;,\;\;
r=\frac{(1-|x|^2)t }{ 1-|x|^2t^2}\;,\;\;t=\tanh(R/2)\;.}&
\end{cases}
\end{equation}
For a simply connected domain $D \subset {\mathbb R}^2$ with nondegenerate boundary and
$x,y \in D$ one can define the hyperbolic metric $\rho_D(x,y)$ as 
$$\label{rhodef}
\rho_{{\mathbb B}^2}(f(x), f(y))
$$ where
$f: D \to {\mathbb B}^2 = f(D)$ is a conformal mapping. This definition yields a well-defined metric, independent of the conformal mapping $f$
 \cite{b,garmar,kela}.

\subsection{Special functions.}
	For $|z|<1$, 
	the Gaussian hypergeometric function is defined by the equality
	$$
	{}_2F_1(a,b;c;z)=\sum_{k=1}^\infty
	\frac{(a)_k(b)_k}{(c)_k}\,\frac{z^k}{k!},
	$$
	{where} $(q)_k$ denotes the Pochhammer symbol, i.e.
	$(q)_k=q(q+1)\ldots(q+(k-1))$ for every natural $k$ and $(q)_0=1$~{\cite{as}}.
	The complete elliptic integral of the first kind \cite{ab,avv,bb}
	\begin{equation}\label{ek}
		{\K}
		(r)=\int_0^1\frac{dt}{\sqrt{(1-t^2)(1-r^2t^2)}}\,,\quad r \in (0,1),
	\end{equation}
	is, in fact, a special case of the Gaussian hypergeometric
	function; we have
	$$
	{\K}
	(r)=\frac{\pi}{2}\,\,{}_2F_1({\textstyle\frac{1}{2}},{\textstyle\frac{1}{2}};1;r^2).
	$$
	We also use the complete elliptic integral of the second kind
	
	$$
	{\E}
	(r)=\frac{\pi}{2}\,\,{}_2F_1({\textstyle-\frac{1}{2}},{\textstyle\frac{1}{2}};1;r^2).
	$$
	\noindent 
	The decreasing homeomorphism $\mu:(0,1)\to (0,\infty)$
	\begin{equation}
	\label{mudef}
	\mu(r)=\frac{\pi}{2} \frac{{\K}(r')}{{\K}(r)},\quad r'=\sqrt{1-r^2},\quad
	0<r<1,
	\end{equation}
	is recurrent in the study of conformal invariants,  its properties
	are studied in \cite{avv}. The functions $\K(r)$ and 
	$\mu(r)$ can be computed
	by means of a simple recursion based on the Landen transformation
	\cite{krv}. 
In this paper, the values of the function $\mu$ and its inverse are computed as described in~\cite{nv4}.

\subsection{Gr\"otzsch condenser} \label{gro} 
The condenser $(\mathbb{B}^n, [0,r]), 0<r<1,$
is called the  Gr\"otzsch condenser. Its capacity
is denoted by $\gamma_n(r).$ For $n=2$ it is one of the few condensers whose capacity is known
\begin{equation} \label{grocap}
\gamma_2(r) \equiv {\rm cap}\, (\mathbb{B}^2, [0,r])= 2 \pi/\mu(r)\,.
\end{equation}
The function $\gamma_2(r)$ has an important role in geometric function theory \cite{hkv} and  short tables with numerical values of the functions $\K(r), \E(r),$ and $\mu(r)$ can be found in \cite[pp. 459-460]{avv}. 
Numerical approximations for the values of $\gamma_3(r)$ were computed by Samuelsson in \cite{sv}.

\subsection{Canonical domains} 
In the study of conformal mappings of multiply connected domains in the extended complex plane $\mathbb{C} \cup \{\infty\}$, a variety of canonical domains onto which a given domain can be mapped have been extensively investigated in the literature (see~\cite{c,hen,Sch} for details). 
Most of these canonical domains are slit domains. Thirty-nine of these canonical slit domains have been catalogued by Koebe in his classical paper~\cite{koe}. 
There are many other canonical slit domains which have not been listed in~\cite{koe} such as the canonical domain obtained by removing rectilinear slits from a strip~\cite[p.~128]{wen}, the parabolic slit domain~\cite{jen,kuh,neh}, the elliptic slit domain~\cite{jen,neh,Sch}, hyperbolic slit domain~\cite{jen,neh}, etc. 
For more details, see~\cite{kuh} and the references cited therein.
Conformal mappings onto slit domains are closely connected to several fundamental concepts in potential theory, including Green's functions, modified Green's functions, and harmonic measures~\cite{c,Sch}, and also play a significant role in addressing various problems in applied mathematics~\cite{c-slit,c}.
 
An important canonical multiply connected domain which is not a slit domain is the circular domain, i.e., a domain all of whose boundary components are circles, see~\cite{hen,koe10,mar}.
Circular domains are ideal for using Fourier series and FFT~\cite{bdhw,dhp}.
Further, analytic formulas exist for several applied problems in circular domains (see the recent monograph~\cite{c} and the references cited therein). See also~\cite{kyyg}.

\subsection{FEM methods} Some of the first studies of numerical approximation of conformal invariants are probably due to Gaier \cite{Ga,gai}
and  Henrici \cite{hen}. For a comprehensive
survey of conformal invariants, see Kuz\'{}mina \cite{kuz}.
In his PhD thesis,  Samuelsson \cite{sam,sv} implemented  his AFEM method in the C++ -language,  a variant of the adaptive FEM method. He applied AFEM to study,  among other things, 
the Gr\"otzsch capacity $\gamma_3(r)$ (some of the results of \cite{sv} 
 are given also in
\cite[pp.244-245]{avv}).  This capacity computation was
reduced to  numerical solution of the  Dirichlet-Neumann problem of a PDE, which for dimensions $n\ge 3$ is
 non-linear and for $n=2$  the linear Laplace equation. Later on, Betsakos-Samuelsson-Vuorinen \cite{bsv} used  AFEM  to approximate
the capacities of several polygonal planar condensers with  simple geometry. The authors of \cite{ps10}
compared these results to the results obtained by other authors \cite[pp.100-103]{ps10}, and reported a $10^{-6}$  agreement of the results.  Further work was carried out by Hakula, Rasila, and Vuorinen in
\cite{hrv1,hrv2, hrv3}, now using Hakula's Mathematica implementation of the $hp$-FEM method. Later on
 Hakula-Nasser-Vuorinen  \cite{hnv0,hnv1,hnv2} carried out numerical studies on conformal invariants with a systematic comparison
between the  $hp$-FEM method and the boundary integral equation method  and again, the numerical
agreement of the results was very good, in some cases of the order $10^{-13}$ \cite{hnv0}.

\subsection{Extremal problems of Gr\"otzsch and Teichm\"uller}  
For a proper subdomain $G$ of $\overline{\mathbb R}^n$ and for $x\in G$, we define $C_x=\gamma_x([0,1))$ where $\gamma_x : [0,1)\to G$ is a curve such that $\gamma_x(0)=x$ and $\gamma_x(t)\to \partial G$ when $t\to 1$. Then for $x,y\in G$ with $x\ne y$, we define\label{page104}
\begin{equation}\label{8.2.}
   \lambda_G(x,y)=\inf_{C_x,C_y}\M\bigl(\D(C_x,C_y;G)\bigr)
\end{equation}
where $\D(C_x,C_y;G)$ is the family of all curves joining $C_x$ and $C_y$ in $G$. Further, for all $x,y\in G$, we define
\begin{equation}\label{8.4.}
  \mu_G(x,y)=\inf_{C_{xy}}\M\bigl(\D(C_{xy},\partial G;G)\bigr)
\end{equation}
where the infimum is taken over all continua $C_{xy}$ such that $C_{xy}=\g([0,1])$ and $\gamma:[0,1]\to G $ is a curve with $\g(0)=x$ and $\g(1)=y$.
Because the modulus is conformally invariant,  it is clear that $\lambda_ G$ and $\mu_G$ are invariant under conformal mappings of $G$.
That is, 
\[
\lambda_{f(G)}(f(x),f(y))=\lambda_G(x,y) \quad
{\rm and }\quad \mu_{f(G)}(f(x),f(y))=\mu_G(x,y) ,
\] 
{if  $f\col G\to f(G)$  is conformal and $x,y\in G$ are distinct.
It is easy  to verify that $\mu_G$ is a metric
if ${\rm cap} ( \partial G)>0$.  If ${\rm cap} (\partial G)>0$, 
$\mu_ G$ is called the {\it modulus metric\/} in $G$. It is also true, but not trivial,  that $\lambda_G^{1/(1-n)}$ is a metric \cite[Thm 10.3]{hkv} in $G.$

In two special cases these functions were studied by
Gr\"otzsch and Teichm\"uller \cite[p.72]{ah}.
In the case $G= \mathbb{B}^2,$  explicit formulas for
$\mu_G$  and  $\lambda_G$ in terms of the hyperbolic metric are given in  \cite[Thm 10.4]{hkv}. In the case
$G= \mathbb{R}^2 \setminus \{0\},$ the formula for
 $\lambda_G$ is given in \cite[p.72]{kuz} and \cite[pp.313-315]{avv}. See also Solynin-Vuorinen  \cite{solvu}.

\section{The boundary integral equation with the generalized Neumann kernel}\label{sec:bie}

The boundary integral equation with the generalized Neumann kernel method has been proposed in~\cite{mn03,wmn,wn,Nas-cr} to solve the Riemann-Hilbert problem in simply and multiply connected domains. 
This method has been then used to solve several problems in complex analysis and potential theory that can be reformulated as a Riemann-Hilbert problem.
One of the novel applications of the integral equation method is the numerical computation of conformal mappings of simply and multiply connected domains. It has been used to compute the conformal mapping from multiply connected domains onto more than forty canonical slit domains~\cite{lmny,Nas-Siam1,Nas-jmaa11,Nas-jmaa13,Nas-cmft15,Nas-PlgCir,Nas-Siam2}. Only the right-hand side of the integral equation is different from one canonical domain to another.
Further, the integral equation has been used effectively to numerically compute several objects in potential theory such as capacities and conformal invariants~\cite{GrNa24,GrNag,hnv0,hnv1,hnv2,knv,LSN17,nnv,nrrvwyz,nrv,nv1,nv2,nv3,nv4}. 

A fast and accurate numerical method for solving the integral equation, based on using the fast multipole method (FMM) and the generalized minimal residual (GMRES) method, has been presented in~\cite{Nas-ETNA}. The MATLAB implementation of this method can be easily modified for a wide range of domains, including, domains with close-to-touching boundaries, domains with non-convex boundaries, domains with piecewise smooth boundaries, and domains of high connectivity. 

In this paper, we review the applications of the boundary integral equation   method with the generalized Neumann kernel to compute capacities and several conformal invariants.

\subsection{The integral equation}\label{sec:mul-ie}

Assume $G$ is a given bounded or unbounded multiply connected domain bordered by $m+1$ smooth Jordan curves $\Gamma_k$, $k=0,1,\ldots,m$.  If $G$ is bounded, then we assume that $\Gamma_0$ is the external boundary component and encloses all the other boundary components $\Gamma_k$, $k=1,\ldots,m$. The total boundary 
\[
\Gamma=\partial G=\bigcup_{k=0}^{m}\Gamma_k,
\]
is oriented such that $G$ is on the left of $\Gamma$.

Each boundary component $\Gamma_k$ is parametrized by a $2\pi$-periodic complex function $\eta_k(t)$ such that  $\eta'_k(t)\ne0$, $t\in J_k=[0,2\pi]$, $k=0,1,\ldots,m$. The total parameter domain $J$ is the disjoint union of the $m+1$ intervals $J_0,J_1,\ldots,J_{m}$, 
\[
J = \bigsqcup_{k=0}^{m} J_k=\bigcup_{k=0}^{m}\{(t,k)\,:\,t\in J_k\}.
\]
That is, the elements of $J$ are ordered pairs $(t,k)$ where $k$ is an auxiliary index indicating which of the intervals contains the point $t$. 
A parametrization of the whole boundary $\Gamma$ is then defined by
\begin{equation}\label{e:eta-1}
\eta(t,k)=\eta_k(t), \quad t\in J_k,\quad k=0,1,\ldots,m.
\end{equation}
For a given $t$, the value of an auxiliary index $k$ such that $t\in J_k$ will be always clear from the context. So we replace the pair $(t,k)$ on the left-hand side of~(\ref{e:eta-1}) by $t$~\cite{Nas-ETNA,wn}. Thus, the function $\eta$ in~(\ref{e:eta-1}) is written as
\begin{equation}\label{eq:mul-eta}
\eta(t)= \left\{ \begin{array}{l@{\hspace{0.5cm}}l}
\eta_0(t),&t\in J_0,\\
\eta_1(t),&t\in J_1,\\
\hspace{0.3cm}\vdots\\
\eta_{m}(t),&t\in J_{m}.
\end{array}
\right.
\end{equation}

Let $A : J \to {\mathbb C}\setminus \{0\}$ be the complex function 
\begin{equation}\label{eq:A}
A(t)=
\left\{ \begin{array}{l@{\hspace{0.5cm}}l}
\eta(t)-\alpha,&\textrm{if $G$ is bounded},\\
1,&\textrm{if $G$ is unbounded},
\end{array}
\right. 
\end{equation}
where $\alpha$ is a given point in the domain $G$.
The generalized Neumann kernel $N(s,t)$ is defined for $(s,t)\in J\times J$ by~\cite{mn03,wmn,wn}
\begin{equation}\label{eq:N}
N(s,t) =
\frac{1}{\pi}\Im\left(\frac{A(s)}{A(t)}\frac{\eta'(t)}{\eta(t)-\eta(s)}\right).
\end{equation}
The kernel $N(s,t)$ is continuous with~\cite{wn}
\begin{equation}\label{eq:Ntt}
N(t,t) =
\frac{1}{\pi}\left(\frac{1}{2}\Im\frac{\eta''(t)}{\eta'(t)}-\Im\frac{A'(t)}{A(t)}\right).
\end{equation}
The integral equation with the generalized Neumann kernel involves also the following kernel 
\begin{equation}\label{eq:M}
M(s,t) =
\frac{1}{\pi}\Re\left(\frac{A(s)}{A(t)}\frac{\dot\eta(t)}{\eta(t)-\eta(s)}\right), \quad (s,t)\in J\times J,
\end{equation}
which has a singularity of cotangent type.
When $s,t\in J_k$ are in the same parameter interval $J_k$, then the kernel $M(s,t)$ has the representation
\begin{equation}\label{eq:Mtt-1}
M(s,t) = -
\frac{1}{2\pi}\cot\left(\frac{s-t}{2}\right)+M_1(s,t), 
\end{equation}
with a continuous kernel $M_1(s,t)$, which takes the values on the diagonal~\cite{wn}
\begin{equation}\label{eq:Mtt-2}
M_1(t,t) =
\frac{1}{\pi}\left(\frac{1}{2}\Re\frac{\eta''(t)}{\eta'(t)}-\Re\frac{A'(t)}{A(t)}\right).
\end{equation}

Let $H$ be the space of all real H\"older continuous functions defined on the boundary $\Gamma$. We define the integral operators $\bN$ and $\bM$ on $H$ by 
\begin{equation}\label{eq:bN}
\bN\rho(s) = \int_J N(s,t) \rho(t) dt, \quad s\in J,
\end{equation}
and
\begin{equation}\label{eq:bM}
\bM\rho(s) = \int_J  M(s,t) \rho(t) dt, \quad s\in J.
\end{equation}
The integral operator $\bN$ is compact and the integral operator $\bM$ is singular. Both operators $\bN$ and $\bM$ are bounded on $H$ and both operators map $H$ into $H$. Further details can be found in~\cite{wmn,wn}. 

\begin{theorem}[{\cite{nv1}}]\label{thm:ie}
For a given function $\gamma\in H$, there exits a unique function $\rho\in H$ and a unique piecewise constant function
\begin{equation} \label{eqn:nu}
\nu(t)= \left\{ \begin{array}{l@{\hspace{0.5cm}}l}
\nu_0,&t\in J_0,\\
\nu_1,&t\in J_1, \\
\hspace{0.3cm}\vdots\\
\nu_m,&t\in J_m, \\
\end{array}
\right.
\end{equation}
with real constants $\nu_0,\nu_1,\ldots,\nu_m$, such that
the formula
\begin{equation}\label{eqn:Af}
f(\eta(t))=\frac{\gamma(t)+\nu(t)+\i\rho(t)}{A(t)}, \quad t\in J,
\end{equation}
defines the boundary values of an analytic function $f$ in $G$ with $f(\infty)=0$ for unbounded $G$. The function $\rho$ is the unique solution of the integral equation
\begin{equation}\label{eqn:ie}
(\bI - \bN) \rho = - \bM \gamma 
\end{equation}
and the piecewise constant function $\nu$ is given by
\begin{equation} \label{eqn:h}
\nu = ( \bM \rho - (\bI - \bN) \gamma )/2.
\end{equation}
\end{theorem}

For simplicity, the function $\nu(t)$ in~\eqref{eqn:nu} will be denoted by
\[
\nu(t)=(\nu_0,\nu_1,\ldots,\nu_m), \quad t\in J.
\]
This notation will be adopted for any piecewise constant function defined on $J$.

\begin{remark}
Doubly connected domains and simply connected domains are particular cases of the above domain $G$ when $m=1$ and $m=0$, respectively. Thus, Theorem~\ref{thm:ie} is valid for doubly connected and simply connected domains. 
In this paper, the integral equation~\eqref{eqn:ie} will be used to solve problems in simply, doubly, and multiply connected domains.
Note that, in the case of a simply connected domain, $\Gamma=\Gamma_0$, $\eta(t)=\eta_0(t)$, the function $\nu(t)=\nu_0$ in~\eqref{eqn:h} is a constant function. 
\end{remark}

\subsection{Numerical solution of the integral equation}

The integral operators $\bN$ and $\bM$ can be best discretized by the Nystr\"om method with the trapezoidal rule since the integrals in~\eqref{eqn:ie} and~\eqref{eqn:h} are over $2\pi$-periodic functions~\cite{atk}.
The smoothness of the integrands in~\eqref{eqn:ie} and~\eqref{eqn:h} depends on the smoothness of the boundary $\Gamma$.		
Further, the function $\bM\gamma$ on the right-hand side of the integral equation~\eqref{eqn:ie} is continuous if the function $\gamma$ is H\"older continuous. 
The stability and convergence of the Nystr\"om method is based on the compactness of the operator $\bN$ in the space of continuous functions equipped with the sup-norm, on the convergence of
the trapezoidal rule for all continuous functions, and on the theory of collectively
compact operator sequences (cf.~\cite{atk}). The numerical solution of the integral equation will converge with a similar rate of convergence as the trapezoidal rule~\cite[p.~322]{atk}.
For smooth boundaries of class $C^{q+2}$ and smooth integrand of class $C^q$, the trapezoidal Nystr\"om method converges with order $O(1/n^q)$ where $n$ is the number of mesh points~\cite{kre90}. 
If the parametrization of the boundary $\Gamma$ is of the class  $C^{\infty}$, then the rate of the convergence of the numerical solution of the integral equation is $\mathcal{O}(e^{-c n})$ where $c$ is a positive constant.

A MATLAB function {\tt fbie} for solving the integral equation~\eqref{eqn:ie} and computing the piecewise constant function $\nu$ in~\eqref{eqn:h} is presented in~\cite{Nas-ETNA}. This MATLAB function is based on discretizing the integrals in~\eqref{eqn:ie} and~\eqref{eqn:h} by the Nystr\"om method with the trapezoidal rule.
For a multiply connected domain $G$ of connectivity $m+1$, discretizing the integral equation~\eqref{eqn:ie} yields an $(m+1)n\times (m+1)n$ linear system where $n$ is the number of mesh points on each boundary component of $G$.
The linear system is then solved by the Generalized Minimal Residual (GMRES) method using the MATLAB function {\tt gmres} where the matrix-vector product in {\tt gmres} is computed by the Fast Multipole Method (FMM) via the MATLAB function {\tt zfmm2dpart} from the MATLAB toolbox FMMLIB2D~\cite{gg} (see also~\cite{gro}). 
The computational cost of the method is $O((m+1)n\log n)$.
See~\cite{Nas-ETNA} for details.

In the numerical computations presented in this paper, the GMRES method is used without restart, with tolerance $10^{-14}$, and with $100$ as the maximal number of allowed iterations. The tolerance of the FMM in {\tt zfmm2dpart} is chosen to be $0.5\times10^{-15}$.

\subsection{Domains with corners}\label{sec:sim-cr}

The integral equation with the generalized Neumann kernel~\eqref{eqn:ie} can be used for domains with piecewise smooth boundaries without cusps~\cite{Nas-cr}.
In this case, the integral operator $\bN$ is not compact, but this operator can be written as a sum of a compact operator and a bounded non-compact operator with norm less than one in suitable function spaces~\cite{Nas-cr}. Hence, we can apply the Fredholm theory to the integral equation with the generalized Neumann kernel~\eqref{eqn:ie} although the operator $\bN$ is not compact~\cite{kre90}. 

For domains with corners, the solution of the integral equation has a singularity in its derivative in the vicinity of the corner points~\cite[p.~390]{atk} and this causes that the equidistant trapezoidal rule yields only poor convergence~\cite{kre90}. 
To achieve a satisfactory accuracy, we discretize the integral equation using a graded mesh quadrature and then applying the Nystr\"om's method~\cite{atk,kre90,kre91}. 
To use such a graded mesh method, we assume that the boundary $\Gamma$ is parametrized by a function $\hat\eta(t)$, $t\in J$. The function $\hat\eta(t)$ is assumed to be smooth with $\hat\eta'(t)\ne0$ for all values of $t\in J$ such that $\hat\eta(t)$ is not a corner point. We assume that $\hat\eta'(t)$ has only the first kind discontinuity at these corner points. If $\hat\eta(\hat t)$ is a corner point, we define $\hat\eta'(\hat t)=\hat\eta'(\hat t-0)$.
Then, as noted in~\cite{kre91}, using the graded mesh method suggested in~\cite{kre90} for discretizing the integral equation is equivalent to parameterizing the boundary $\Gamma$ by 
\[
\eta(t)=\hat\eta(\delta(t)), \quad t\in J,
\]
where the function $\delta(t)$ is defined in~\cite[pp.~696--697]{LSN17} which is chosen to remove the discontinuity in the derivatives of the solution of the integral equation at the corner points. 
To compute discrete parameterizing of the boundaries of polygonal and polycircular domains, we can use the two MATLAB functions \verb|polygonp.m| and \verb|plgcirarcp.m|, respectively, which are based on the method described above.  These two functions are available in \url{https://github.com/mmsnasser/polycircular}. 
With the parametrization $\eta(t)=\hat\eta(\delta(t))$, the integral equation can be solved using the MATLAB function {\tt fbie} as in the case of smooth domains. 
For the numerical computations presented in this paper, the grading parameter $p$ in Kress's method~\cite{kre90} is chosen to be $p=3$. In this case, the numerical results presented in~\cite{hnv0,nrrvwyz} and the numerical results presented in this paper illustrate that the rate of convergence is $\mathcal{O}(n^{-q})$ with $q\le p$.

\section{Computation of the capacity of generalized condensers}\label{sec:cgc}

\subsection{Capacity of generalized condensers}

The conformal capacity of a condenser is one of the key notions of potential theory of elliptic partial differential equations \cite{gor,hkm} and it has numerous applications to geometric function theory, both in the plane and in higher dimensions~\cite{du,  gor, hkv, hkm}. For the basic facts about capacities, the reader is referred to \cite{du, gor, hkv, hkm}.

Consider the \emph{generalized condenser} $C=(\Omega,E,\delta)$ where $\Omega \subset \CC$ is a bounded domain, $E=\cup_{j=1}^{m}E_j$ where $E_1,\ldots,E_m$ are $m$ compact disjoint non-empty subsets of $\Omega$, $\delta=\{\delta_k\}_{k=1}^{m}$ is a collection of $m$ real numbers. 
The \emph{conformal capacity} of this generalized condenser is defined as  
\cite{du,gor, hkm,hkv} 
\begin{align}\label{def_condensercap}
\capa(C)=\inf_{u\in A}\int_{\Omega}|\nabla u|^2 dm,
\end{align}
where $A$ is the class of $C^\infty_0(G)$ functions $u: \Omega\to \R$ 
with  $u(x) \ge \delta_k$ for all $x \in E$ and $dm$ is the $2$-dimensional 
Lebesgue measure. 
This is a special case of the generalized condenser defined in~\cite{du} (see also~\cite{nv1}). Further, when $\delta_1=\cdots=\delta_m=1$, the generalized condenser $C$ is the classical condenser defined in \S\ref{sec:int} above~\cite{du,  gor, hkv, hkm}. For this case, we will denote the condenser by $C=(\Omega,E)$.

Here we assume that $\Gamma_1=\partial E_1, \ldots, \Gamma_m=\partial E_m$ and $\Gamma_0=\partial \Omega$ are piecewise smooth Jordan curves. 
Hence $G=\Omega\backslash E$ is a multiply connected domain of connectivity $m+1$ with $\Gamma=\partial G=\cup_{k=0}^m\Gamma_k$ and the infimum in~\eqref{def_condensercap} is attained by a harmonic function $u$. This extremal function $u$ is the unique solution of the Laplace equation in $G$ with boundary values given by $u=\delta_k$ on $\Gamma_k$, $k=1,\ldots,m$, and $u=0$ on $\Gamma_0$~\cite{du}.
The capacity can be then expressed in terms of the extremal function $u$ as
\begin{equation}\label{eq:cap}
	\capa(C)=\iint \limits_{G}|\nabla u|^2 dxdy,
\end{equation}
which, using Green's formula~\cite[p. 4]{du}, implies that
\begin{equation}\label{eq:cap-n}
	\capa(C)=\int_{\Gamma} u\frac{\partial u}{\partial {\bf n}}ds
\end{equation}
where $\partial u/\partial {\bf n}$ denotes the directional derivative of $u$ along the outward normal. 
Since the Dirichlet integral is conformally invariant, the cases for which $E_1, \ldots, E_m$ are slits can be handled with the help of auxiliary conformal mappings which transform the slits to smooth Jordan curves.


\subsection{The numerical method}

A boundary integral equation with the generalized Neumann kernel has been presented in~\cite{nv1} for the numerical computation of the capacity $\capa(C)$ as well as the values of the potential function $u(z)$ for $z\in G$. 
The harmonic function $u$ is the real part of an analytic function $F$ in $G$ which is not necessarily single-valued. Assume that $\alpha_k$ is an auxiliary point in the interior of $\Gamma_k$ for each $k=1,2,\ldots,m$, then the function $F$ can be written as~\cite{Gak,garmar,Mik64}
\begin{equation}\label{eq:F-u}
F(z)=g(z)-\sum_{k=1}^{m} a_k\log(z-\alpha_k)
\end{equation}
where $g$ is a single-valued analytic function in $G$ and $a_1,\ldots,a_{m}$ are undetermined real constants such that~\cite[\S31]{Mik64}
\begin{equation}\label{eq:ak}
a_k=\frac{1}{2\pi}\int_{\Gamma_k}\frac{\partial u}{\partial{\bf n}}ds, \quad k=1,2,\ldots,m.
\end{equation}
Since $u=0$ on $\Gamma_{0}$ and $u=\delta_k$ on $\Gamma_k$ for $k=1,2,\ldots,m$, then in view of~\eqref{eq:cap-n} and~\eqref{eq:ak}, we have
\begin{equation}\label{eq:cap-k}
\capa(C)=\sum_{k=1}^{m} \int_{\Gamma_k}\delta_k\frac{\partial u}{\partial{\bf n}}ds
=2\pi\sum_{k=1}^{m}\delta_ka_k.
\end{equation}
Equation~\eqref{eq:cap-k} gives us a simple formula for computing the capacity of the generalized condenser $C=(\Omega,E,\delta)$ in terms of the values of the constants $a_1,\ldots,a_m$ and $\delta_1,\ldots,\delta_m$. The constant $2\pi\delta_k a_k$ can be considered as the contribution of the compact set $E_k$ to the capacity $\capa(C)$ for $k=1,2,\ldots,m$.

The function $g(z)$ can be written as
\[
g(z)=(z-\alpha)f(z)+c
\]
where $\alpha$ is a given point in the domain $G$, $f(z)$ is a single-valued analytic function in $G$, and $c$ is an undetermined constant. Without loss of generality we can assume that $c$ is a real constant. Then, the function $F$ can be written as
\begin{equation}\label{eq:F-log}
F(z)=(z-\alpha)f(z)+c-\sum_{k=1}^{m} a_k\log(z-\alpha_k),
\end{equation}
and hence the function $u(z)$ is given for $z\in G$ by
\begin{equation}\label{eq:F-u2}
u(z)=\Re[(z-\alpha)f(z)]+c-\sum_{k=1}^{m} a_k\log|z-\alpha_k|.
\end{equation}

For each $k=1,2,\ldots,m$, let the function $\gamma_k$ be defined by
\begin{equation}\label{eq:gam-k}
	\gamma_k(t)=\log|\eta(t)-\alpha_k|.
\end{equation}
Since $u=0$ on $\Gamma_{0}$ and $u=1$ on $\Gamma_k$ for $k=1,2,\ldots,m$, then the function $f(z)$  is the solution of the Riemann-Hilbert problem
\begin{equation}\label{eq:rhp}
\Re[A(t)f(\eta(t))]=\hat\nu(t)-c+\sum_{k=1}^{m} a_k\gamma_k(t),
\end{equation}
where $\hat\nu(t)=(0,\delta_1, \ldots,\delta_m)$. For the coefficient function $A$ given by~\eqref{eq:A}, the solution of the Riemann-Hilbert problem~\eqref{eq:rhp} is unique~\cite{Nas-ETNA}.

Assume that the boundary $\Gamma$ is parametrized by the function $\eta(t)$ in~\eqref{eq:mul-eta} and assume the function $A(t)$ is defined by~\eqref{eq:A}, i.e., $A(t)=\eta(t)-\alpha$ since $G$ is bounded. Let the kernels $N(s,t)$ and $M(s,t)$ of the integral operators $\bN$ and $\bM$, respectively, be formed with these functions $\eta(t)$ and $A(t)$.
For $k=1,\ldots,m$, it follows from Theorem~\ref{thm:ie} that the integral equation
\begin{equation}\label{eq:ie}
	\mu_k-\bN\mu_k=-\bM\gamma_k,
\end{equation}
has a unique solution $\mu_k$, the function $\nu_k$ given by
\begin{equation}\label{eq:h}
	\nu_k=[\bM\mu_k-(\bI-\bN)\gamma_k]/2\,.
\end{equation}
This is a piecewise constant function, i.e., $\nu_k(t)=(\nu_{0,k}, \nu_{1,k}, \ldots,\nu_{m,k})$ where $\nu_{0,k},\nu_{1,k},\ldots,\nu_{m,k}$ are real constants, and
\begin{equation}\label{eq:fk-eta}
A(t)f_k(\eta(t))=\gamma_k(t)+\nu_k(t)+\i\mu_k(t), \quad t\in J,
\end{equation}
are boundary values of an analytic function $f_k(z)$ in $G$. Then the function
\begin{equation}\label{eq:f-fk}
f(z) = \sum_{k=1}^m a_k f_k(z), \quad z\in\Omega\cup\Gamma,
\end{equation}
is the unique solution of the Riemann-Hilbert problem~\eqref{eq:rhp} in $G$ if and only if
\begin{equation}\label{eq:c-hk}
\hat\nu(t)-c=\sum_{k=1}^m a_k \nu_k(t), \quad t\in J.
\end{equation}
It follows from~\eqref{eq:c-hk} that the values of the $m+1$ unknown real constants $a_1,\ldots,a_{m},c$ can be computed by solving the $(m+1)\times(m+1)$ linear system
\begin{equation}\label{eq:sys-method}
	\left[\begin{array}{ccccc}
		\nu_{0,1}    &\nu_{0,2}    &\cdots &\nu_{0,m}      &1       \\
		\nu_{1,1}    &\nu_{1,2}    &\cdots &\nu_{1,m}      &1       \\
		\vdots       &\vdots       &\ddots &\vdots       &\vdots  \\
		\nu_{m,1}    &\nu_{m,2}    &\cdots &\nu_{m,m}      &1       \\
	\end{array}\right]
	\left[\begin{array}{c}
		a_1    \\a_2    \\ \vdots \\ a_{m} \\  c 
	\end{array}\right]
	= \left[\begin{array}{c}
		0 \\  \delta_1 \\  \vdots \\ \delta_m  
	\end{array}\right].
\end{equation}
The linear system~\eqref{eq:sys-method} has a unique solution~\cite{nv1}. The size of the system~\eqref{eq:sys-method} is usually small and hence can be solved by the MATLAB ``$\backslash$'' function. 

By computing the functions $\mu_k$ and $\nu_k$ for $k=1,\ldots,m$ numerically, we obtain approximations of the boundary values of the analytic function $f_k(z)$ through~\eqref{eq:fk-eta}. We obtain also the entries of the coefficient matrix of the linear system~\eqref{eq:sys-method}. By solving the linear system~\eqref{eq:sys-method} for the $m+1$ real constants $a_1,\ldots,a_{m},c$, the value of the capacity $\capa(C)$ can be computed by~\eqref{eq:cap-k} and the boundary values of the function $f(z)$ can be computed through
\[ 
A(t)f(\eta(t))=\sum_{k=1}^m a_k\left(\gamma_k(t)+\nu_k(t)+\i\mu_k(t)\right), \quad t\in J.
\]
Then, the values of $f(z)$ for $z\in G$ can be computed by the Cauchy integral formula. A MATLAB function {\tt fcau} for fast and accurate computation of the Cauchy integral formula is presented in~\cite{Nas-ETNA}. Then the values of $u(z)$ can be computed for $z\in G$ by~\eqref{eq:F-u2}.

\begin{example}[Circular ring]\label{ex:ring}
For $a \in (0,1)$, let $\Omega=\B^2$ and $E=\{z\in\CC\,:\, |z|\le a\}$.  
The exact formula for the capacity of the condenser $C=(\Omega,E)$ is given by
\[
\capa(C)=\frac{2\pi}{\log(1/a)}.
\]
We use the above integral equation method to compute the approximate values of $\capa(C)$ for several values of $n$ where $n$ is the number of mesh points on each boundary component of $G=\Omega\backslash E$. The relative error in the computed values are presented in Figure~\ref{fig:err-1} (left). It is clear that the method converges exponentially as the boundary of the domains are $C^\infty$-smooth.
\end{example}

\begin{example}[Square in square: \cite{bsv,ps10,p2}]\label{ex:sq}
For $a \in (0,1)$, let $\Omega= (-1,1)\times  (-1,1)$ and $A= [-a,a]\times  [-a,a]$.  
The exact formula for the capacity of the condenser $C=(\Omega,A)$ is given in~\cite{bsv} as
\[
\capa(C)=\frac{4 \pi}{\mu(r)} =\frac{4}{\pi}\mu((u/v)^2)
\]
where the second equality follows from~\cite[Exer. 7.33(3)]{hkv} and \eqref{mudef}, and where
\[
c=\frac{1-a}{1+a}\,,\qquad u=\mu^{-1}(\frac{\pi c}{2})\,,\qquad  
v=\mu^{-1}(\frac{\pi }{2c})\,,\qquad r= \left(\frac{u-v}{u+v}\right)^2\,.
\]
We use the above integral equation method and the relative error in the computed approximate values are presented in Figure~\ref{fig:err-1} (right). The boundary of the domain has corners and the method converges algebraically with $O(n^{-3})$.
\end{example}

\begin{example}[Disk with a polygonal hole]\label{ex:dpg}
For $a \in (0,1)$, let $\Omega= \B^2$ and $E$ is the closure of the set of points interior to the regular polygon with the $\ell$ vertices 
\[
z_k = a\,e^{\i(k-1)\theta},\quad \theta=\frac{2\pi}{\ell}, \quad k=1,2,\ldots,\ell, \quad \ell\ge 3.
\]

The area of $E$ is $\frac{1}{2}\ell a^2\sin\theta$ and the perimeter of $E$ is $2\ell a\sin(\theta/2)$. Let $r_1$ be such that the area of the disk $B^2(0,r_1)$ is equal to the area of $E$ and let $r_2$ be such that the perimeter of the disk $B^2(0,r_2)$ is equal to the perimeter of $E$, i.e.,
\[
r_1 = a\sqrt{\sin\theta}\sqrt{\frac{\ell }{2\pi}}, \quad
r_2 = \frac{\ell a}{\pi}\sin\frac{\theta}{2}.
\]

The exact value of the capacity of the condenser $C_k=(\Omega,\overline{B^2(0,r_k)})$ is $\capa(C_k)=2\pi/\log(1/r_k)$, $k=1,2$. The exact formula for the capacity of the condenser $C=(\Omega,E)$ is unknown. We compute the values of $\capa(C)$ using the integral equation method for $0.05<a<0.95$ and the obtained numerical results are presented in Figure~\ref{fig:err-plg}. It is clear that the condenser with inner disk will always have larger capacity than the condenser with inner polygon with the same area and smaller capacity than the polygon with the same perimeter. Further, the capacity of a condenser with an inner disk can be considered as a good approximation for the capacity of a condenser with inner regular polygon with the same area as the inner disk.

\end{example}

\begin{example}[Lens shaped plate: \cite{hnv0,mc}]\label{ex:lens}

For $a \in (0,1)$, let $\Omega=\B^2$ and $E$ be the closure of the lens domain bordered by the two circular arcs, the first one passes through the points $a$, $-\i s$, $-a$, and the second one passes through the points $-a$, $\i s$, $a$.
In~\cite{hnv0} the capacity of this condenser was computed with the integral equation and $hp$-FEM methods and in \cite[(5.1), (5.2)]{hnv0} upper and lower bounds were given in terms of the hyperbolic perimeter of the lens shaped set.

The exact value of the capacity of the condenser $C=(\Omega,E)$ for $s=0$ is $\capa(C)=2\pi/\mu(2a/(1+a^2))$ and for $s=a$ is $\capa(C)=2\pi/\log(1/a)$. The exact formula for the capacity is unknown for $0<s<a$. We compute the values of $\capa(C)$ for $0<s<a$ and the obtained numerical results are presented in Figure~\ref{fig:err-len} (left) for $a=0.5$. Assume that the circular arc passes through the points $-a$, $\i s$, $a$ is on a circle with radius $r$, then the values of $\capa(C)$ are approximated in~\cite{mc} by
\begin{equation}\label{eq:len-est}
\capa(C) \approx \frac{2\pi}{\log\left(\frac{2(\pi-\theta)}{\pi a}\right)}, \quad \theta = \arctan\left(\frac{a}{r-s}\right).
\end{equation}
The values of these estimations are presented also in Figure~\ref{fig:err-len} (left) and the differences between the estimations and the computed values of the capacity are presented in Figure~\ref{fig:err-len} (right). The estimated values agree with the computed values of the capacity in particular when $s$ is close to $a$. 
\end{example}

\begin{figure}[ht] %
	{\includegraphics[width=0.325\textwidth]{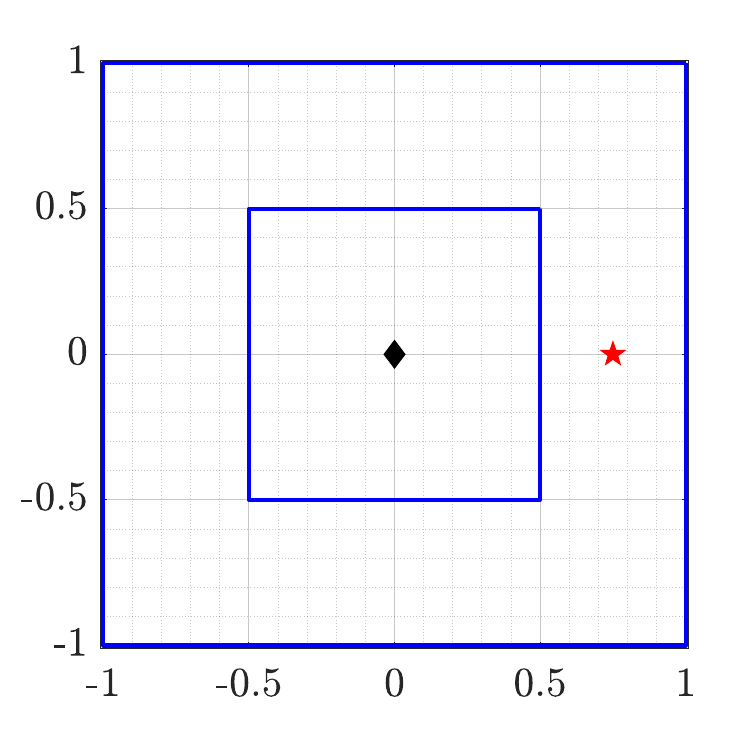}}
	{\includegraphics[width=0.325\textwidth]{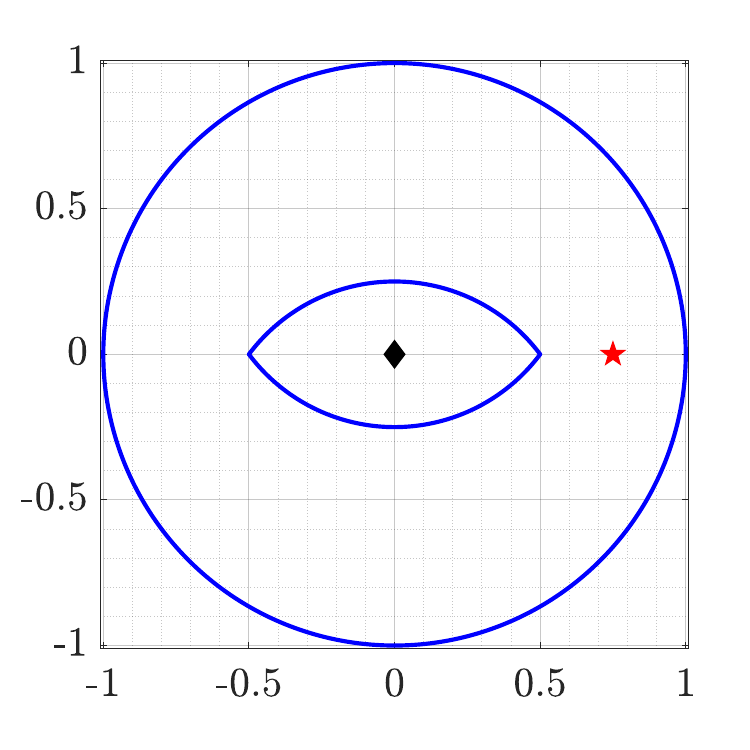}}
	{\includegraphics[width=0.325\textwidth]{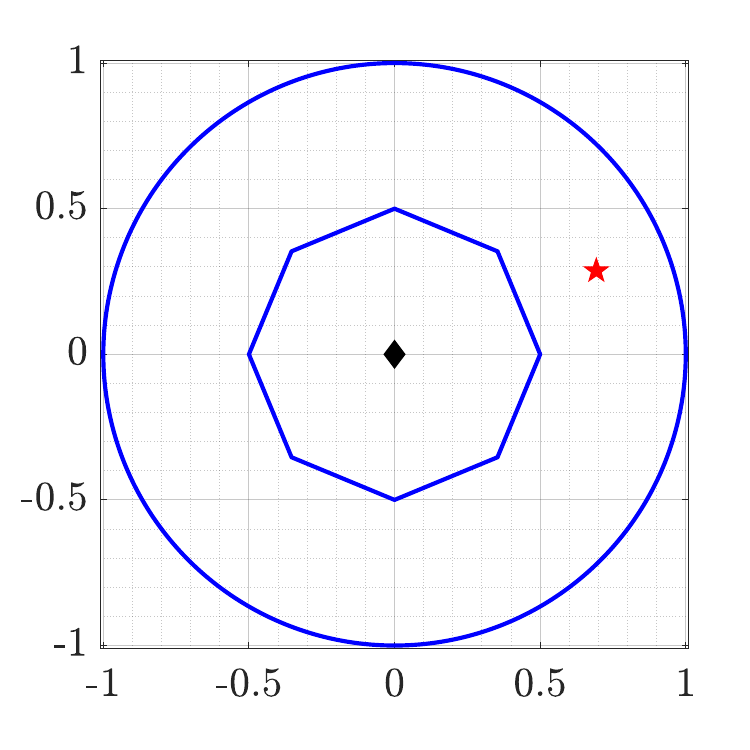}}
\caption{The domains of the condensers considered in Examples~\ref{ex:sq}, \ref{ex:lens}, and~\ref{ex:dpg}.}
	\label{fig:cap-dom}
\end{figure}

\begin{figure}[ht] %
	\centering{
	\hfill{\includegraphics[width=0.45\textwidth]{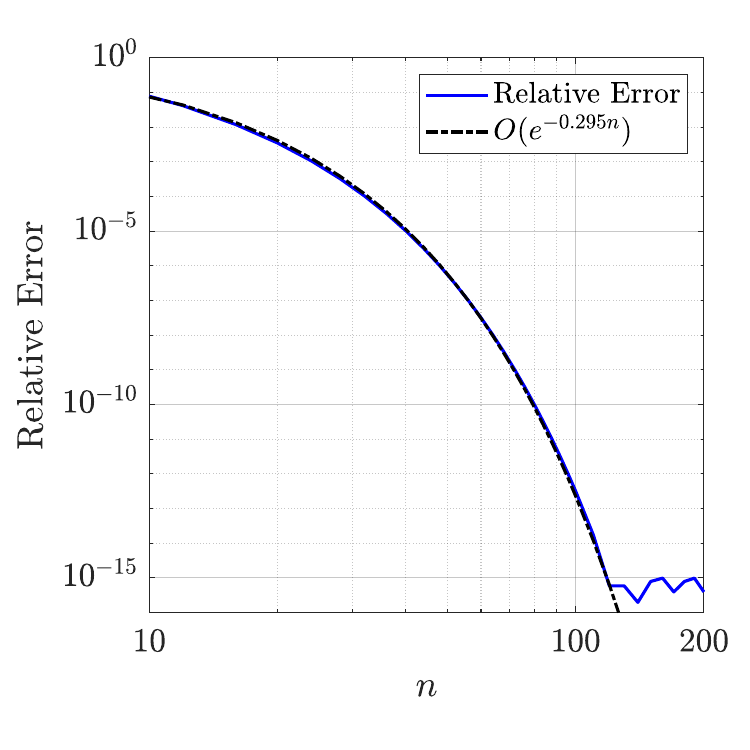}}
	\hfill{\includegraphics[width=0.45\textwidth]{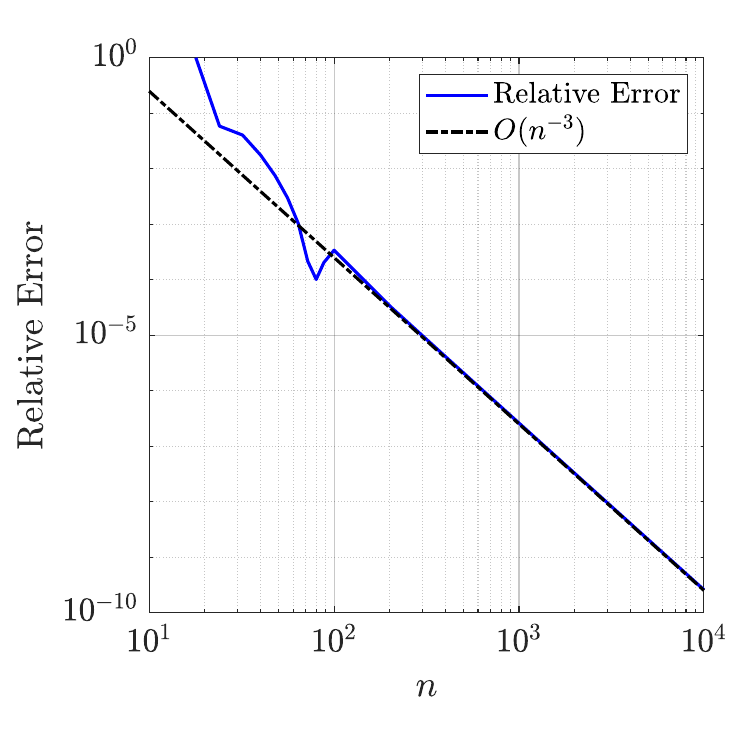}}
	\hfill}
\caption{The relative error in the computed values of the capacities of the condensers in Examples~\ref{ex:ring} (left) and \ref{ex:sq} (right) where $n$ is the number of mesh points on each boundary component of $G$.}
	\label{fig:err-1}
\end{figure}

\begin{figure}[ht] %
	\centering{
	\hfill{\includegraphics[width=0.45\textwidth]{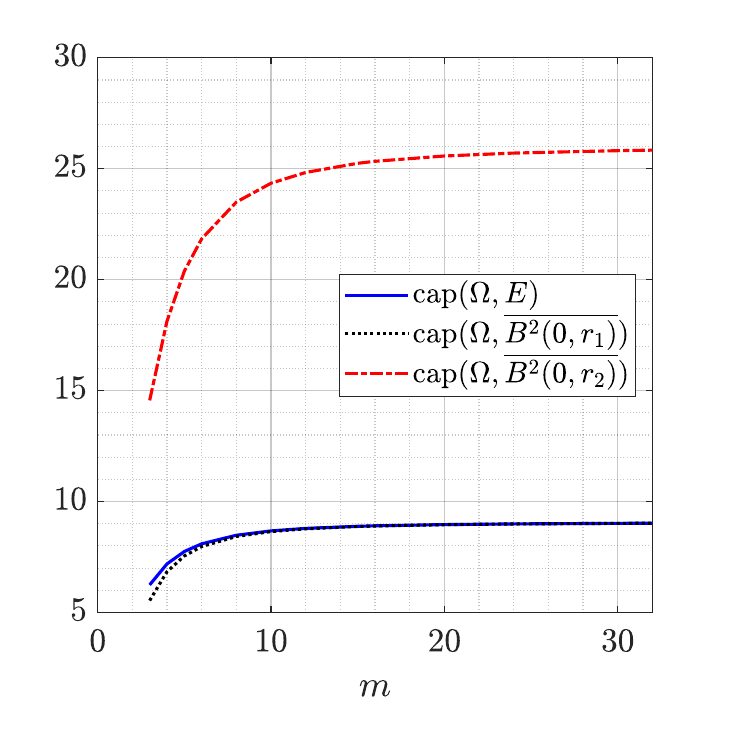}}
	\hfill}
\caption{The computed values of the capacity of the condenser $(\Omega,E)$ in Examples~\ref{ex:dpg} and the capacities of $(\Omega,\overline{B^2(0,r_k)})$, $k=1,2$.}
	\label{fig:err-plg}
\end{figure}

\begin{figure}[ht] %
	\centering{
	\hfill{\includegraphics[width=0.45\textwidth]{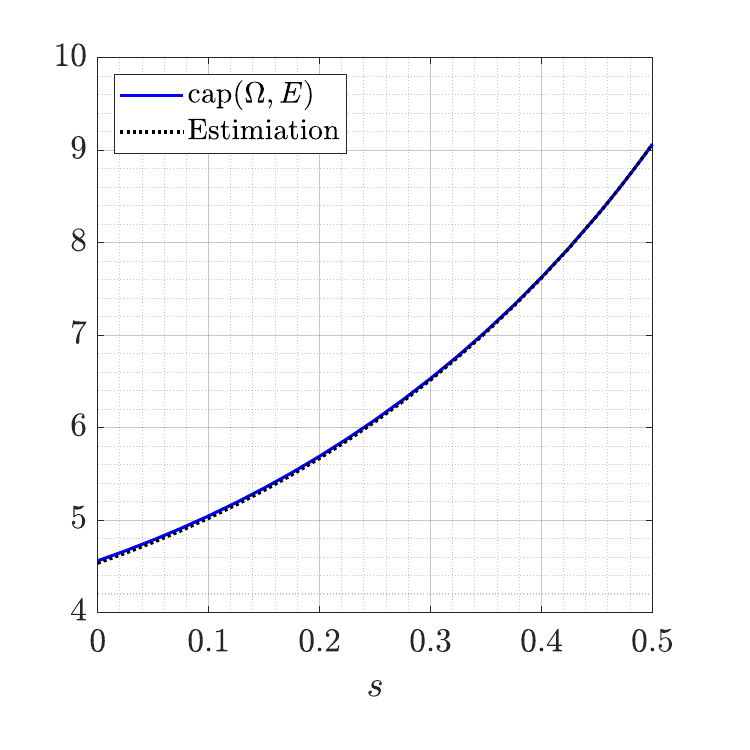}}
	\hfill{\includegraphics[width=0.45\textwidth]{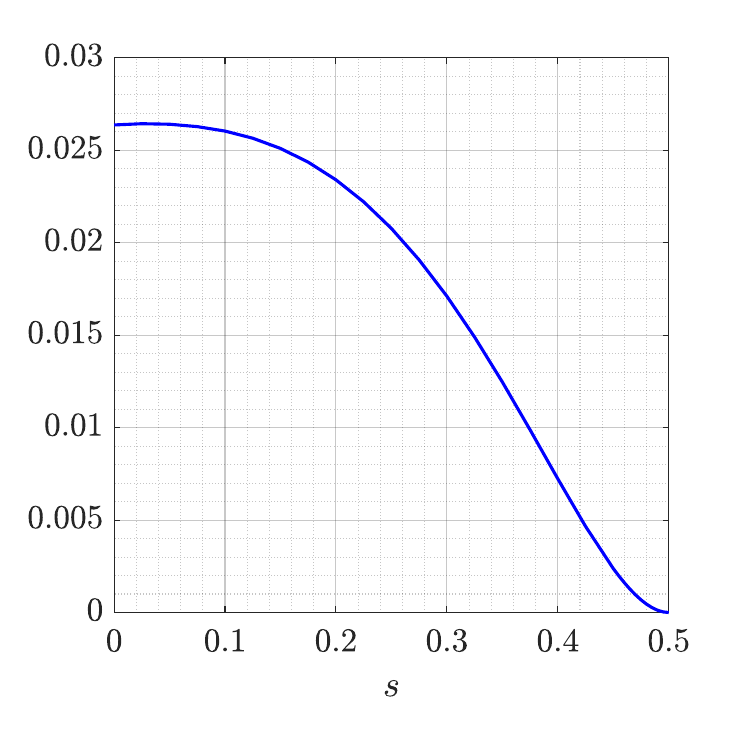}}
	\hfill}
\caption{Left: The computed values of the capacity of the condenser in Examples~\ref{ex:lens}. Right: absolute values of the differences between the computed values of the capacity and the estimate~\eqref{eq:len-est}.}
	\label{fig:err-len}
\end{figure}

\begin{example}[A disk with $7$ circular holes]\label{ex:7d}
Let $\Omega=\B^2$ and $E=\cup_{k=1}^7E_k$ where $E_k=\{z\in\CC\,:\, |z-a_k|\le r\}$,
\[
a_k=(0.1+k/10)e^{\i(k-1)\pi/2}, \quad k=1,2,\ldots,7,
\]
and $0<r<\sqrt{0.0325}\approx0.18$ (see Figure~\ref{fig:7d-1} (left) for $r=0.1$). We choose also $\delta_1=\cdots=\delta_7=1$, and hence, by~\eqref{eq:cap-k}, the capacity of the condenser $C=(\Omega,E)$ is given by $\capa(C)=2\pi\sum_{k=1}^7a_k$ where $a_k$ are defined by~\eqref{eq:ak}.
The exact value of the capacity of the condenser $C$ is unknown. 

We use the above integral equation method to compute the approximate values of $\capa(C)$ for several values of $r$ and the computed values of the constants $a_1,\ldots,a_7$ are presented in Figure~\ref{fig:7d-1} (right). 
The constant $2\pi a_k$ is the contribution of the compact set $E_k$ to the capacity of the generalized condenser $C=(\Omega,E)$.
It is clear that the values of $a_1,\ldots,a_7$ depend on the location of the disk $E_k$ as well as on its radius.
\end{example}

\begin{figure}[ht] %
	\centering{
	\hfill{\includegraphics[width=0.45\textwidth]{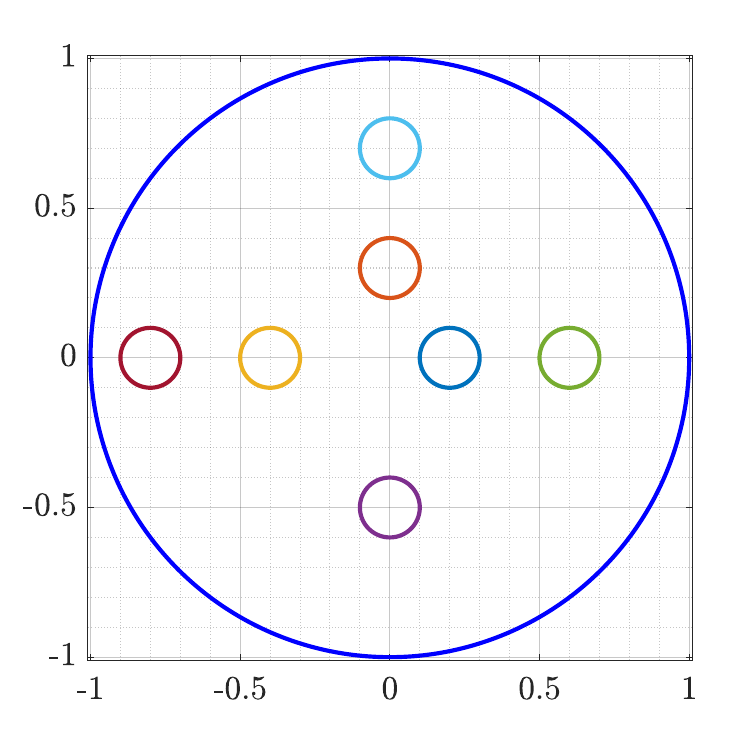}}
	\hfill{\includegraphics[width=0.45\textwidth]{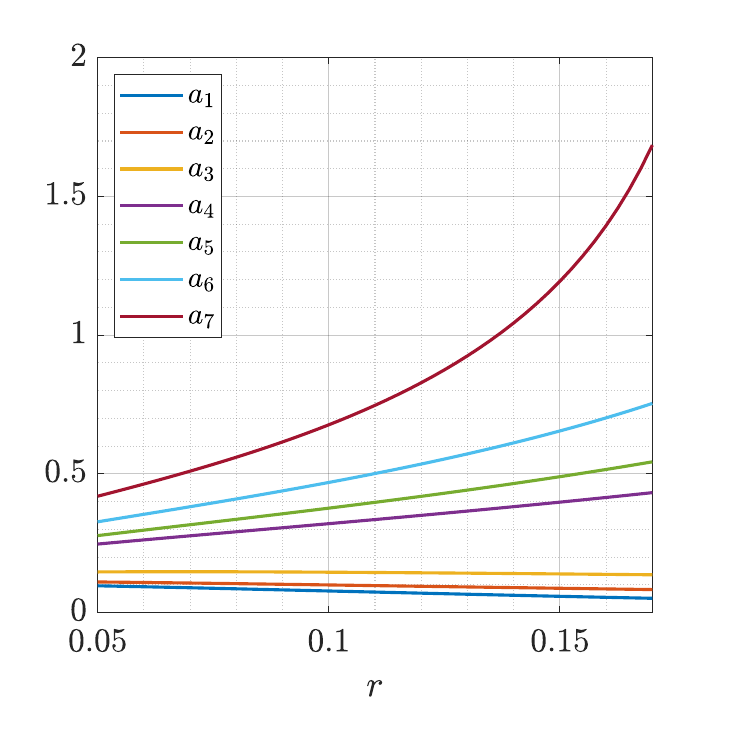}}
	\hfill}
\caption{Left: The domain of the condenser in Examples~\ref{ex:7d}. Right: The computed values of the constants $a_1,\ldots,a_7$ for $0.01\le r\le 0.17$.}
	\label{fig:7d-1}
\end{figure}

\section{Computation of logarithmic capacity}\label{sec:lc}

\subsection{Logarithmic capacity}

Let $E$ be a compact set ($E$ is not a single point) whose complement $G=\overline{\C}\backslash E$ is an unbounded multiply connected domain of connectivity $m+1$ bordered by $m+1$ piecewise smooth Jordan curves $\Gamma_0,\Gamma_1,\ldots,\Gamma_m$ and let $\Gamma=\cup_{k=0}^{m}\Gamma_k$. Let $g_G(z)$ be the Green function of $G$ with pole at infinity. Then the logarithmic capacity of $E$, denoted here by $\capa_l(E)$, is defined by~\cite{LSN17,r}
\begin{equation}\label{eq:cap-log-1}
\capa_l(E)=\lim_{z\to\infty}\exp\left(\log|z|-g_G(z)\right).
\end{equation}

If $G$ is simply connected, then its Green's function is given by $g_G(z)=\log|\Phi(z)|$, where $w=\Phi(z)$ is the uniquely determined
conformal map from $G$ onto the unbounded domain $\{z\in\overline{\C}:|z|>1\}$ exterior to the unit circle with normalization $\Phi(\infty)=\infty$ and $\Phi'(\infty)>0$.

Walsh~\cite{walsh} proved a direct generalization of the classical Riemann mapping theorem in which he replaced the exterior of the unit circle by a lemniscatic domain of the form
\[
\Omega=\{z\in\overline{\C}:|U(z)|>\kappa\},
\]
where
\[
U(z)=\prod_{j=0}^{m}(z-\beta_j)^{\ell_j},
\]
$\beta_0,\beta_1,\ldots,\beta_m\in\C$ are pairwise distinct, $\ell_0,\ell_1,\ldots,\ell_m$ are positive real numbers with $\sum_{j=0}^m\ell_j=1$, and $\kappa>0$. Then~\cite{LSN17}
\[
g_G(z)=\log|U(z)|-\log(\kappa),
\]
which implies that
\begin{equation}\label{eq:cap-log}
\capa_l(E)=\lim_{z\to\infty}\exp\left(\log|z|-\log|U(z)|+\log(\kappa)\right) = \kappa.
\end{equation}

Analytic formulas of logarithmic capacity are known only for a few cases. For example, the logarithmic capacity of a disk of radius $r$ is $r$ and the logarithmic capacity of the an ellipse with semi-axes $a$ and $b$ is $(a+b)/2$ (see~\cite[Table~1]{LSN17}, \cite[pp.172-173]{lan} for more examples).

\subsection{The numerical method}

The integral equation~\eqref{eq:ie} has been used in~\cite{LSN17} to develop a fast and accurate numerical method for computing the logarithmic capacity. 
Let $\Gamma$ be parametrized by the function $\eta(t)$ in~\eqref{eq:mul-eta}, let the function $A(t)$ be defined by~\eqref{eq:A}, i.e., $A(t)=1$ since $G$ is unbounded, and let the kernels $N(s,t)$ and $M(s,t)$ of the integral operators $\bN$ and $\bM$, respectively, be formed with these functions $\eta(t)$ and $A(t)$. 
For each $j=0,1,\ldots,m$, we choose an auxiliary point $\alpha_j$ in the interior of the Jordan curve $\Gamma_j$ and define the function $\gamma_j(t)$ by
\[
\gamma_j(t) = -\log|\eta(t)-\alpha_j|, \quad t\in J.
\]
Let $\rho_j$ be the unique solution of the integral equation~\eqref{eq:ie} and let the piecewise constant function $\nu_j(t)=(\nu_{0,j},\nu_{1,j},\ldots,\nu_{m,j})$ be given by~\eqref{eq:h}, $j=0,1,\ldots,m$. Then $\log(\kappa)$ is computed by solving the following uniquely solvable linear system~\cite{LSN17} 
\begin{equation}\label{eq:sys-log-cap}
	\left[\begin{array}{ccccc}
		\nu_{0,0}    &\nu_{0,1}    &\cdots &\nu_{0,m}      &-1       \\
		\nu_{1,0}    &\nu_{1,1}    &\cdots &\nu_{1,m}      &-1       \\
		\vdots       &\vdots       &\ddots &\vdots       &\vdots  \\
		\nu_{m,0}    &\nu_{m,1}    &\cdots &\nu_{m,m}      &-1       \\
		1            &1            &\cdots &1              & 0       \\
	\end{array}\right]
	\left[\begin{array}{c}
		\ell_1    \\  \ell_2    \\ \vdots \\ \ell_m \\  \log(\kappa) 
	\end{array}\right]
	= \left[\begin{array}{c}
		0 \\  0 \\  \vdots \\ 0  \\ 1  
	\end{array}\right].
\end{equation}
By computing $\log(\kappa)$, we obtain the logarithmic capacity by~\eqref{eq:cap-log}.

\begin{example}[Circular domain]\label{ex:log-cap}
Let $E$ be the union of the five disjoint disks $E_1=\{z\in\C:|z|\le0.8\}$, $E_{2,3}=\{z\in\C:|z\pm r|\le1\}$, and $E_{4,5}=\{z\in\C:|z\pm r\i|\le1\}$, $r>1.8$. The disks are close to each other for small values of $r$ and far for large values of $r$.
Then $G=\overline{\C}\backslash E$ is an unbounded multiply connected domain $G$ of connectivity $5$.
We use the integral equation method with $n=2^{10}$ to compute the values of the logarithmic capacity $\capa_l(E)$ for several values of $r$, $2\le r\le 8$. The obtained values of  $\capa_l(E)$ vs. $r$ are shown in Figure~\ref{fig:log-cap}. Since $\capa_l(E_1)=0.8$ and $\capa_l(E_k)=1$ for $k=2,\ldots,5$, then $\sum_{k=1}^5 \capa_l(E_k)=4.8$. It is clear from Figure~\ref{fig:log-cap} that $\capa_l(\cup_{k=1}^5E_k) < \sum_{k=1}^5 \capa_l(E_k)$ for small values of $r$ and $\capa_l(\cup_{k=1}^5E_k) > \sum_{k=1}^5 \capa_l(E_k)$ for large values of $r$, i.e., the logarithmic capacity is not subadditive (see also~\cite{PP}).
\end{example}

\begin{figure}[ht] %
	\centering{
	\hfill{\includegraphics[width=0.45\textwidth]{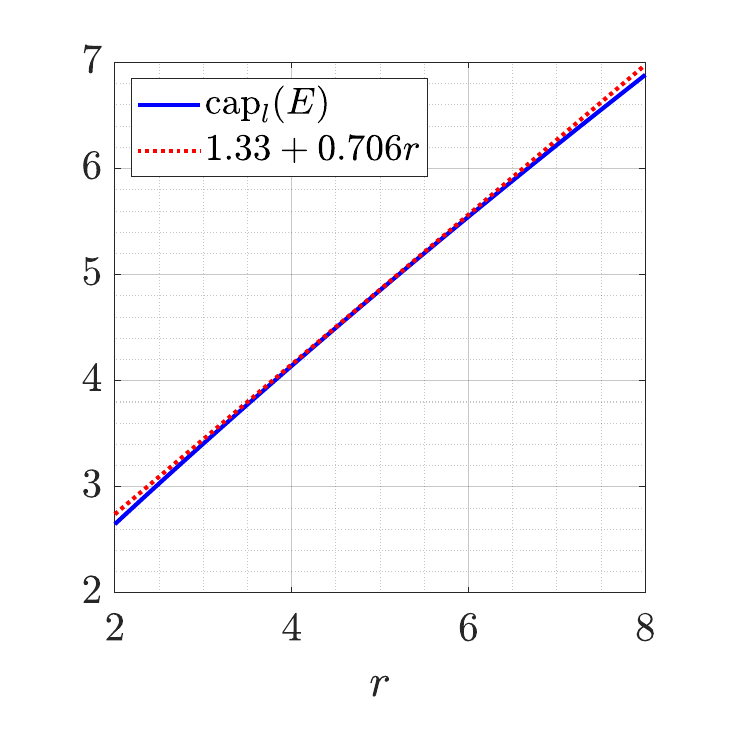}}
	\hfill{\includegraphics[width=0.45\textwidth]{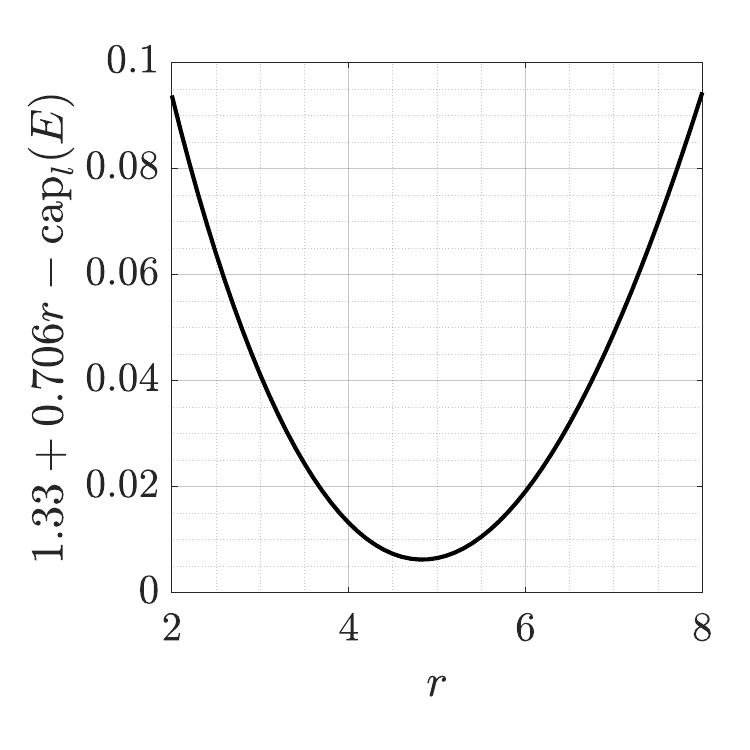}}
	\hfill}
\caption{The values of $\capa_l(E)$ for the compact set $E$ in Example~\ref{ex:log-cap} for several values of $r$.}
	\label{fig:log-cap}
\end{figure}


\section{Hyperbolic and elliptic capacities}\label{sec:hec}

In this section, we assume that $G$ is a doubly connected domain, i.e., $m=1$ in Section~\ref{sec:bie}. We assume that $G$ is bordered by the piecewise smooth Jordan curves $\Gamma_0$ and $\Gamma_1$ and $\Gamma=\Gamma_0\cup\Gamma_1$. We assume also that $\Gamma$ is parametrized by the function $\eta(t)$ in~\eqref{eq:mul-eta}, the function $A(t)$ is defined by~\eqref{eq:A}, and the kernels $N(s,t)$ and $M(s,t)$ of the integral operators $\bN$ and $\bM$, respectively, are formed with these functions $\eta(t)$ and $A(t)$.

\subsection{Conformal mapping onto an annulus}\label{sec:doub-cm}

The conformal mapping $w=\Phi(z)$ from the doubly connected domain $G$ onto the annulus $\{w\in\C\;:\;q<|w|<1\}$ can be computed using the following method from~\cite[\S4.1]{Nas-Siam1}.

For bounded $G$, let $w=\Phi(z)$ be the unique conformal mapping with the normalization
\[
\Phi(\alpha)>0,
\]
where $\alpha$ is a given auxiliary point in $G$. Let $z_1$ be a given point in the simply connected domain interior to $\Gamma_1$ and let the function $\gamma$ be defined by
\begin{equation}
\gamma(t) = -\log\left|\frac{\eta(t)-z_1}{\alpha-z_1}\right|, \quad t\in J.
\end{equation}
Let also $\rho(t)$ be the unique solution of the integral equation~\eqref{eqn:ie} and let the piecewise constant function $\nu(t)=(\nu_0,\nu_1)$ be given by~\eqref{eqn:h}.
Then the function $f$ with the boundary values~\eqref{eqn:Af} is analytic in the domain $G$ and the conformal mapping $\Phi$ is
given by
\begin{equation}\label{eq:Phi-an-b}
\Phi(z)=c\left(\frac{z-z_1}{\alpha-z_1}\right)e^{(z-\alpha)f(z)}, \quad z\in G\cup\Gamma,
\end{equation}
where $c=e^{-\nu_0}=\Phi(\alpha)>0$ and the modulus $q$ is given by
\begin{equation}\label{eq:R-an-b}
q=e^{\nu_1-\nu_0}.
\end{equation}

For unbounded $G$, let $w=\Phi(z)$ be the unique conformal mapping with the normalization
\[
\Phi(\infty)>0.
\]
Let $z_0$ be a given point in the simply connected domain interior to $\Gamma_0$, let $z_1$ be a given point in the simply connected domain interior to $\Gamma_1$, and let the function $\gamma$ be defined by
\begin{equation}
\gamma(t)=-\log\left|\frac{\eta(t)-z_1}{\eta(t)-z_0}\right|, \quad t\in J.
\end{equation}
Let also $\rho(t)$ be the unique solution of the integral equation~\eqref{eqn:ie} and let the piecewise constant function $\nu(t)=(\nu_0,\nu_1)$ be given by~\eqref{eqn:h}.
Then the function $f$ with the boundary values~\eqref{eqn:Af} is analytic in the unbounded domain $G$ with $f(\infty)=0$ and the conformal mapping $\Phi$ is given by
\begin{equation}\label{eq:Phi-an-u}
\Phi(z)=c\left(\frac{z-z_1}{z-z_0}\right)e^{f(z)}, \quad z\in G\cup\Gamma,
\end{equation}
where $c=e^{-\nu_0}=\Phi(\infty)>0$ and the modulus $q$ is given by
\begin{equation}\label{eq:R-an-u}
q=e^{\nu_1-\nu_0}.
\end{equation}

By computing the functions $\rho$ and $\nu$, we obtain approximations of the boundary values of the analytic function $f(z)$ by~\eqref{eqn:Af}. Then the boundary values of the mapping function $\Phi(z)$ can be computed from~\eqref{eq:Phi-an-b} for bounded domains and from~\eqref{eq:Phi-an-u} for unbounded domains. The values of $\Phi(z)$ for $z\in G$ can be computed by the Cauchy integral formula.

\subsection{Hyperbolic capacity}\label{sec:doub-hc}

Let $E$ be a compact and connected set (not a single point) in the unit disk $\B^2$.
The hyperbolic capacity of $E$, $\capa_h(E)$, is defined by~\cite[p.~19]{Vas02}
\begin{equation}\label{eq:hyp-cap}
\capa_h(E)=\lim_{n\to\infty}\left[\max_{z_1,\ldots,z_n\in E}\prod_{1\le k<j\le n}
\left|\frac{z_k-z_j}{1-z_k\overline{z_j}}\right|
\right]^{\frac{2}{n(n-1)}}.
\end{equation}
For the hyperbolic capacity, we assume $G$ is the bounded doubly connected domain defined by $G=\B^2\backslash E$. The domain $G$ can be mapped conformally onto an annulus $q<|w|<1$. Hence the hyperbolic capacity $\capa_h(E)$ is given by~\cite{DuKu}
\begin{equation}\label{eq:hyp-cap3}
\capa_h(E)=q.
\end{equation}
The hyperbolic capacity is invariant under conformal mappings.

\begin{example}[Hyperbolic capacity of an ellipse]\label{ex:hyc-ellipse}
Let $E$ be the closed region bordered by the ellipse
\[
\eta_1(t)=0.75\cos t-\i 0.5\sin t, \quad 0\le t\le 2\pi.
\]
Let $G=\B^2\backslash E$, $\alpha=0.75\i$, and $z_1=0$. Then, we used the above method to compute the conformal mapping from $G$ onto the annulus $q<|w|<1$, and hence $\capa_h(E)=q$. The obtained approximate results (with $n=2^{10}$) of $\capa_h(E)$ is $0.634497711721981$.
\end{example}

\subsection{Elliptic capacity}\label{sec:doub-ec}

Let $E$ be a compact and connected set (not a single point) in the unit disk $\B^2$.
We define the antipodal set $E^\ast=\{-1/\overline{a}\,:\,a\in E\}$. Since we assume $E\subset\B^2$, we have $E\cap E^\ast=\emptyset$ (in this case, the set $E$ is called ``elliptically schlicht''~\cite{DuKu}). The elliptic capacity of $E$, $\capa_e(E)$, is defined by~\cite{DuKu}
\begin{equation}\label{eq:ell-cap}
\capa_e(E)=\lim_{n\to\infty}\left[\max_{z_1,\ldots,z_n\in E}\prod_{1\le k<j\le n}
\left|\frac{z_k-z_j}{1+z_k\overline{z_j}}\right|
\right]^{\frac{2}{n(n-1)}}.
\end{equation}

To compute the elliptic capacity, we assume $G$ is the doubly connected domain between $E$ and $E^\ast$. Such a domain $G$ can be bounded (if $0\in E$) or unbounded (if $0$ is in the exterior of $E$). For both cases, the domain $G$ can be mapped conformally onto an annulus $r<|w|<1/r$. Then the elliptic capacity is given by~\cite{DuKu}
\[
\capa_e(E)=r.
\]

We can use the method described above to map the domain $G$ onto an annulus $q<|w|<1$ which is conformally equivalent to the annulus $r<|w|<1/r$ with $r=\sqrt{q}$. Thus, we have
\begin{equation}\label{eq:ell-cap2}
\capa_e(E)=\sqrt{q}.
\end{equation}

\begin{remark}
For a closed and connected subset $E$ of the unit disk $\B^2$, Duren and K\"uhnau~\cite{DuKu} have proved that
\[
\capa_e(E)\le\capa_h(E),
\]
with equality if and only if $E=-E$. 
\end{remark}

\begin{example}[Elliptic capacity of an ellipse]\label{ex:ell-ellipse}
Let $E$ be the closed region bordered by the ellipse
\[
\eta_1(t)=0.75\cos t-\i 0.5\sin t, \quad 0\le t\le 2\pi,
\]
which is oriented clockwise.
By the symmetry of $E$, the boundary of $E^\ast$ can be parametrized by
\[
\eta_0(t)=\frac{1}{0.75\cos t-\i 0.5\sin t}, \quad 0\le t\le 2\pi,
\]
which is oriented counterclockwise. Let $G$ be the bounded doubly connected domain in the exterior of the curve with parametrization $\eta_1(t)$ and in the interior of the curve with the parametrization $\eta_0(t)$.
Let also $\alpha=0.75\i$ and $z_1=0$. Then, we used the above method to compute the conformal mapping from $G$ onto the annulus $q<|w|<1$, and hence $\capa_e(E)=\sqrt{q}$. The obtained approximate results (with $n=2^{10}$) of $\capa_e(E)$ is $0.634497711721982$.
\end{example}

\begin{remark}
Note that $E=-E$ in Examples~\ref{ex:hyc-ellipse} and~\ref{ex:ell-ellipse} which implies that the hyperbolic and elliptic capacities must be equal. For the numerically computed values in Examples~\ref{ex:hyc-ellipse} and~\ref{ex:ell-ellipse}, note that
\[
|\capa_h(E)-\capa_e(E)|=9.99\times10^{-16},
\]
which could be considered as an estimation of the error in the computed values of the hyperbolic and elliptic capacities.
\end{remark}

\section{Reduced modulus}\label{sec:rm}

\subsection{Conformal mappings of simply connected domains}\label{sec:sim-cm}
Assume that $G$ is a simply connected domain, i.e., $m=0$ in the notation of Section~\ref{sec:bie}. We assume that $G$ is bordered by the piecewise smooth Jordan curve $\Gamma=\Gamma_0$ which is parametrized by the function $\eta(t)=\eta_0(t)$ in~\eqref{eq:mul-eta}. Assume also that the function $A(t)$ is defined by~\eqref{eq:A}, and the kernels $N(s,t)$ and $M(s,t)$ of the integral operators $\bN$ and $\bM$, respectively, are formed with these functions $\eta(t)$ and $A(t)$.
In this subsection, we review a numerical method based on using the integral equation~\eqref{eqn:ie} for computing a conformal mapping $w=\Phi(z)$ from a simply connected domain $G$ onto the unit disk $\B^2$~\cite{Nas-cmft15,nrrvwyz}. Some applications of the reduced
modulus to geometric function theory are discussed
in \cite{akn}.

For a bounded domain $G$, the mapping function $\Phi$ is unique if we assume that
\begin{equation}\label{eq:Phi-cond}
\Phi(\alpha)=0, \quad \Phi'(\alpha)>0
\end{equation}
where $\alpha$ is a given auxiliary point in the domain $G$.
Let 
\begin{equation}\label{eqn:gam}
\gamma(t)=-\log|\eta(t)-\alpha|, \quad t\in[0,2\pi],  
\end{equation}
let $\rho$ be the unique solution of the integral equation~\eqref{eqn:ie}, and let the constant $\nu$ be given by~\eqref{eqn:h}. 
Then, the mapping function $\Phi$ with normalization~\eqref{eq:Phi-cond}
can be written for $z\in G\cup\Gamma$ as
\begin{equation}\label{eq:Phi-b}
\Phi(z)=c(z-\alpha)e^{(z-\alpha)f(z)}
\end{equation}
where the function $f(z)$ is analytic in $G$ with the boundary values~\eqref{eqn:Af} and $c=e^{-\nu}=\Phi'(\alpha)>0$.

For unbounded $G$, the mapping function $\Phi$ is unique if we assume that
\begin{equation}\label{eq:Phi-condu}
\Phi(\infty)=0, \quad \Phi'(\infty)>0.
\end{equation}
Let 
\begin{equation}\label{eqn:gamu}
\gamma(t)= \log|\eta(t)-z_0|, \quad t\in[0,2\pi],  
\end{equation}
where $z_0$ is a given auxiliary point in the bounded domain interior to $\Gamma$,
let $\rho$ be the unique solution of the integral equation~\eqref{eqn:ie}, and let the constant $\nu$ be given by~\eqref{eqn:h}. 
Then, the mapping function $\Phi$ with normalization~\eqref{eq:Phi-condu}
can be written for $z\in G\cup\Gamma$ as
\begin{equation}\label{eq:Phi-u}
\Phi(z)=\frac{c}{z-z_0}e^{f(z)}
\end{equation}
where the function $f(z)$ is analytic in $G$ with $f(\infty)=0$, the boundary values of $f$ are given by~\eqref{eqn:Af}, and
and $c=e^{-\nu}=\Phi'(\infty)>0$.

By computing the function $\rho$ and the constant $\nu$ numerically, we obtain approximations of the boundary values of the analytic function $f(z)$ through~\eqref{eqn:Af} and hence the boundary values $\Phi(\eta(t))$ of the mapping function $\Phi(z)$ can be computed by~\eqref{eq:Phi-b} for bounded $G$ and by~\eqref{eq:Phi-u} for unbounded $G$. The values of the mapping function $w=\Phi(z)$ can be then computed for $z\in G$ using the Cauchy integral formula 
\begin{equation}\label{eqn:Phi-z}
\Phi(z)=\frac{1}{2\pi\i}\int_{\Gamma} \frac{\Phi(\eta)}{\eta-z}d\eta
=\frac{1}{2\pi\i}\int_{0}^{2\pi} \frac{\Phi(\eta(t))}{\eta(t)-z}\eta'(t)dt 
=\frac{1}{2\pi\i}\int_{0}^{2\pi} \frac{\zeta(t)}{\eta(t)-z}\eta'(t)dt.
\end{equation} 

The function $\zeta(t)=\Phi(\eta(t))$, $t\in[0,2\pi]$, is a parametrization of the unit circle.
Thus, for computing the values of the inverse mapping function, we first compute numerically the derivative $\zeta'(t)$ by interpolating both $\Re[\zeta(t)]$ and $\Im[\zeta(t)]$ by trigonometric interpolating polynomials and then differentiating the interpolating polynomials. These polynomials can be computed with FFT~\cite{weg05}.
Then, for a bounded domain $G$, the values of the inverse mapping function $z=\Phi^{-1}(w)$ for $w\in \B^2$ can be  computed using the Cauchy integral formula,
\begin{equation}\label{eqn:Phi-w}
\Phi^{-1}(w)=\frac{1}{2\pi\i}\int_{\partial\B^2} \frac{\Phi^{-1}(\zeta)}{\zeta-w}d\zeta
=\frac{1}{2\pi\i}\int_{0}^{2\pi} \frac{\Phi^{-1}(\zeta(t))}{\zeta(t)-w}\zeta'(t)dt 
=\frac{1}{2\pi\i}\int_{0}^{2\pi} \frac{\eta(t)}{\zeta(t)-w}\zeta'(t)dt.
\end{equation}
For unbounded $G$, note that the function $\Phi^{-1}(w)$ has a simple pole at $w=0$ and the function $g(w)=w\Phi^{-1}(w)$ is analytic in $\B^2$. Thus, the values of the inverse mapping function $z=\Phi^{-1}(w)$ for $w\in \B^2$ can be  computed using the Cauchy integral formula,
\begin{equation}\label{eqn:Phi-wu}
\Phi^{-1}(w)=\frac{1}{w}g(w)=\frac{1}{w}\frac{1}{2\pi\i}\int_{\partial\B^2} \frac{\zeta\Phi^{-1}(\zeta)}{\zeta-w}d\zeta
=\frac{1}{w} \frac{1}{2\pi\i}\int_{0}^{2\pi} \frac{\zeta(t)\eta(t)}{\zeta(t)-w}\zeta'(t)dt.
\end{equation}

\subsection{Reduced modulus of simply connected domains}\label{sec:simb-rm}

If $G$ is a bounded simply connected domain and $w=\Phi(z)$ is the conformal mapping from $G$ onto the unit disk $\B^2$ with the normalization~\eqref{eq:Phi-cond}, then the conformal radius of $G$ with respect to the point $\alpha$ is defined by~\cite{Mit}
\begin{equation}\label{eq:R}
R(G,\alpha)=\frac{1}{\Phi'(\alpha)}.
\end{equation}
The reduced modulus of the domain $G$ with respect to the point $\alpha$ is then given by~\cite[p.~16]{Vas02}
\begin{equation}\label{eq:m}
m(G,\alpha)=\frac{1}{2\pi}\log R(G,\alpha)=-\frac{1}{2\pi}\log\Phi'(\alpha).
\end{equation}

\begin{remark}
Vasil'ev~\cite[Section 2.2.1]{Vas02} assumed that the mapping function $w=\hat\Phi(z)$ from $G$ onto $\B^2$ satisfies the normalization $\hat\Phi(\alpha)=0$ and $\hat\Phi'(\alpha)=1$. Hence $w=\hat\Phi(z)$ maps the domain $G$ onto the disk $|w|<R$ where $R=R(G,\alpha)$ is the conformal radius of $G$ with respect to the point $\alpha$. This is equivalent to the above definition~\eqref{eq:R} since $\hat\Phi(z)=R\,\Phi(z)$.
\end{remark}

For unbounded simply connected domain $G$, the reduced modulus of the domain $G$ with respect to the point $\infty$ is defined by~\cite[p.~17]{Vas02}
\begin{equation}\label{eq:mu}
m(G,\infty)=-\frac{1}{2\pi}\log\Phi'(\infty),
\end{equation}
where $w=\Phi(z)$ is the conformal mapping from $G$ onto the unit disk $\B^2$ with the normalization~\eqref{eq:Phi-condu}. 

For both bounded and unbounded simply connected domain $G$, it follows from \S\ref{sec:sim-cm} that $m(G,\alpha)=\nu/2\pi$ and $m(G,\infty)=\nu/2\pi$ where the constant $\nu$ is given by~\eqref{eqn:h}.

\begin{example}[Domain interior to an ellipse]\label{ex:redm-in}
We consider the simply connected domain $G$ in the interior of the ellipse
\[
\eta(t)=\cosh(r+\i t)=\cosh r\cos t+\i\sinh r\sin t, \quad 0\le t\le 2\pi, \quad 0<r.
\]
Let $w=\Phi(z)$ be the unique conformal mapping from the interior of the ellipse onto the interior of the unit circle with the normalization $\Phi(0)=0$ and $\Phi'(0)>0$. The  exact form of the inverse conformal mapping $z=\Phi^{-1}(w)$ is given in~\cite{KS06}. In particular, it was shown in~\cite{KS06} that $(\Phi^{-1})'(0)=\pi/(2\sqrt{s}\K(s))$ where $s=\mu^{-1}(2r)$. Hence, $\Phi'(0)=2\sqrt{s}\K(s)/\pi$. Thus, $R(G,0)=1/\Phi'(0)=\pi/(2\sqrt{s}\K(s))$ and hence
\[
m(G,0)=\frac{1}{2\pi}\log\frac{\pi}{2\sqrt{s}\K(s)}, \quad s=\mu^{-1}(2r).
\]
Figure~\ref{fig:redm-sim} (left) shows
the relative error in the numerically computed values using the integral equation method presented in~\S\ref{sec:sim-cm} with $n=2^8\,.$ 
To study the effect of the location of $\alpha$ on the values of $m(g,\alpha)$, we define the function $u(x,y)$ for all $x$ and $y$ such that $x+\i y\in G$ by
\[
u(x,y) = m(G,x+\i y).
\]
We use the integral equation method with $n=2^{12}$ to compute the values of the function $u(x,y)$. The level curves for the function $u(x,y)$ corresponding to the several levels are shown in Figure~\ref{fig:redm-simL}. Each of these level curves describes the locations of the point $\alpha$ for which the values $m(G,\alpha)$ are a constant. The maximum value of $m(G,\alpha)$ occurs when $\alpha=0$.
\end{example}

\begin{example}[Domain exterior to an ellipse]\label{ex:redm-ex}
Consider the simply connected domain $G$ in the exterior of the ellipse
\[
\eta(t)=\cosh(r-\i t), \quad 0\le t\le 2\pi, \quad 0<r < 1.
\]
We can easily show that the function
\[
z=\Psi(w)=\frac{1}{2}\left(e^{-r} w+\frac{1}{e^{-r} w}\right),
\]
maps the unit disk onto the domain exterior of the ellipse. Hence, the inverse mapping
\begin{equation}\label{eq:map-to-ellipse}
w=\Phi(z)=\frac{e^r}{z\left(1+\sqrt{1-1/z^2}\right)},
\end{equation}
maps the domain exterior to the ellipse onto the unit disk, where the branch of the square root is chosen such that $\sqrt{1}=1$. 
The function $\Psi(w)$ is a simple modification of the Joukowski map~\cite{hen}.
It is clear that the function $\Phi$ satisfies $\Phi(\infty)=0$ and $\Phi'(\infty)=e^r/2$. Hence, 
\[
m(G,\infty)=-\frac{1}{2\pi}\log\frac{e^r}{2}=\frac{1}{2\pi}\left(\log2-r\right).
\]
The relative error in the numerically computed values using the method presented in~\S\ref{sec:sim-cm} with $n=2^8$ is presented in Figure~\ref{fig:redm-sim} (right). 
\end{example}

\begin{figure}[ht] %
\centerline{
		\scalebox{0.6}{\includegraphics[trim=0cm 0.5cm 0cm 0.5cm,clip]{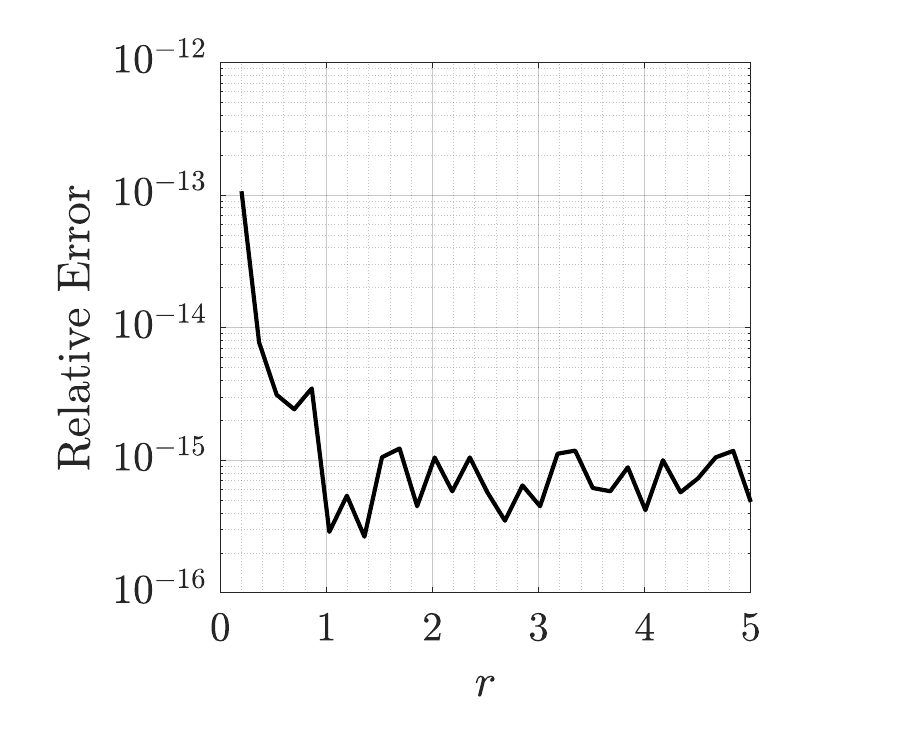}}
		\hfill
		\scalebox{0.6}{\includegraphics[trim=0cm 0.5cm 0cm 0.5cm,clip]{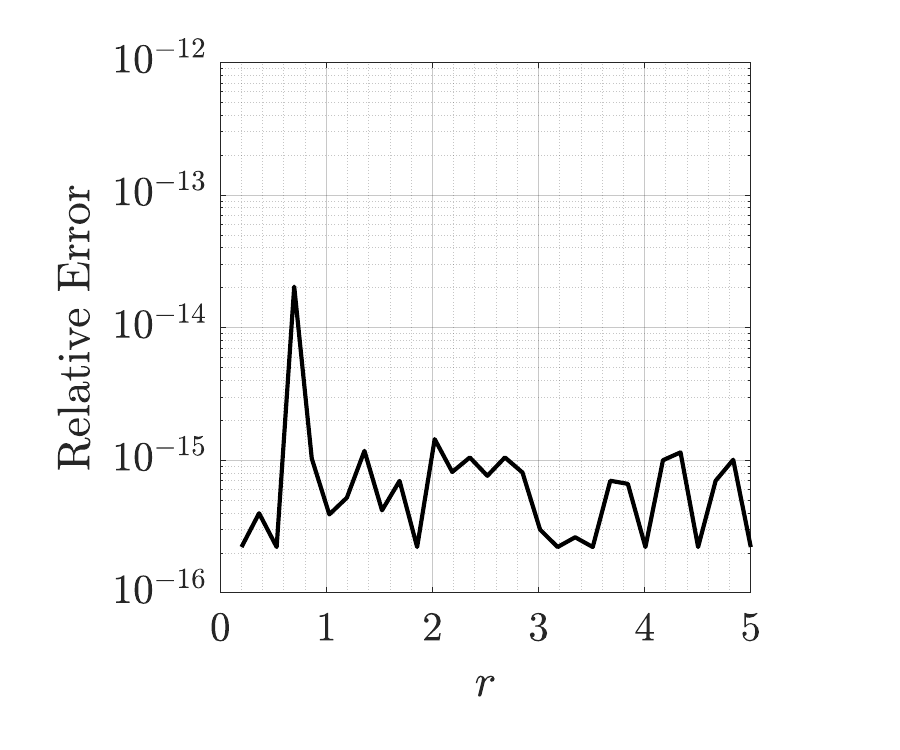}}		
	}
\caption{The relative error in the computed values of the reduce modulus vs. $r$ for Example~\ref{ex:redm-in} (left) and Example~\ref{ex:redm-ex} (right) obtained with $n=2^8$ where $n$ is the number of mesh points on the boundary of $G$.}
	\label{fig:redm-sim}
\end{figure}

\begin{figure}[ht] %
\centerline{
		\scalebox{0.6}{\includegraphics[trim=0cm 0.5cm 0cm 0.5cm,clip]{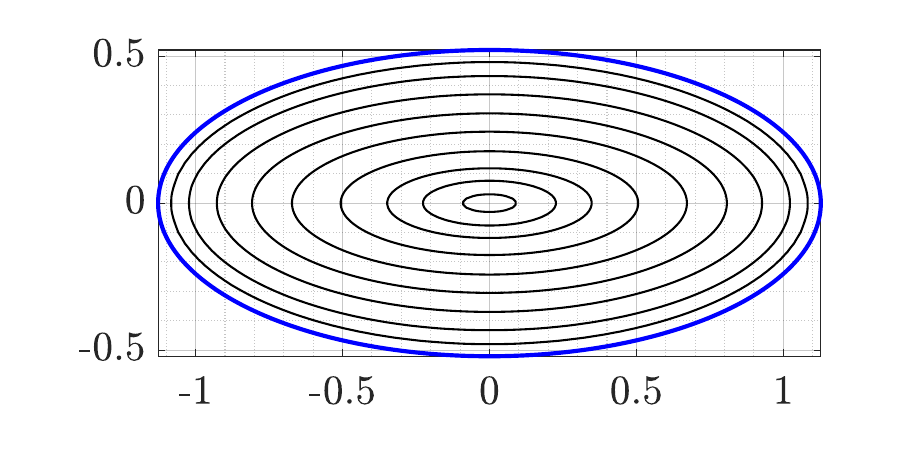}}		
	}
\caption{The contour lines of the function $u(x,y)$ in Example~\ref{ex:redm-in}.}
	\label{fig:redm-simL}
\end{figure}

\subsection{Generalized reduced modulus}\label{sec:mul-rm}

A generalization of the reduced modulus to multiply connected domains has been proposed by Mityuk~\cite{Mit}. For multiply connected domains, several canonical domains are available and two of these canonical domains have been considered in~\cite{Mit}. The boundary integral equation with the generalized Neumann kernel has been used in~\cite{knv} to compute the generalized reduced module for these two canonical domains. Here, we consider one of these canonical domains, namely, the unit disk with circular slits.

Let $G$ be a given bounded multiply connected domain of connectivity $m+1$, let $\Gamma=\partial G=\cup_{k=0}^m \Gamma_k$ be parametrized by the function $\eta(t)$ in~\eqref{eq:mul-eta}, and let the function $A(t)$ be defined by~\eqref{eq:A}, i.e., $A(t)=\eta(t)-\alpha$. Let also the kernels $N(s,t)$ and $M(s,t)$ of the integral operators $\bN$ and $\bM$, respectively, be formed with these functions $\eta(t)$ and $A(t)$.

For the given domain $G$, there exists a conformal mapping $w=\Phi(z)$ from the domain $G$ onto the canonical domain $D$ obtained by removing $m$ concentric circular slits from the unit disk. We assume that these slits are
subarcs of circles centered at $0$
with radii $R_1,\ldots,R_m$ which are undetermined real constants.
With the normalization
\begin{equation}\label{eq:Phi-cond-m}
\Phi(\alpha)=0, \quad \Phi'(\alpha)>0,
\end{equation}
this conformal mapping is unique. 
Hence, the definition~\eqref{eq:m} of the reduced modulus of bounded simply connected domains can be generalized to  the bounded multiply connected domain $G$. That is, the generalized reduced modulus of the bounded multiply connected domain $G$ with respect to the point $\alpha$ and the canonical domain $D$ can be defined by
\begin{equation}\label{eq:m-m}
m(G,\alpha)=-\frac{1}{2\pi}\log\Phi'(\alpha).
\end{equation}

Let 
\begin{equation}\label{eqn:gam-m}
\gamma(t)=-\log|\eta(t)-\alpha|, \quad t\in[0,2\pi],  
\end{equation}
let $\rho$ be the unique solution of the integral equation~\eqref{eqn:ie} and let the piecewise constant function $\nu(t)=(\nu_0,\nu_1,\ldots,\nu_m)$ be given by~\eqref{eqn:h}. 
Then, the mapping function $\Phi$ with normalization~\eqref{eq:Phi-cond-m} can be written for $z\in G\cup\Gamma$ as
\begin{equation}\label{eq:Phi-b-m}
\Phi(z)=c(z-\alpha)e^{(z-\alpha)f(z)}
\end{equation}
where the function $f(z)$ is analytic in $G$ with the boundary values~\eqref{eqn:Af}, $c=e^{-\nu_0}=\Phi'(\alpha)>0$, and $R_k=e^{\nu_k-\nu_0}$ for $k=1,\ldots,m$.
Then, it follows from~\eqref{eq:Phi-cond-m} and~\eqref{eq:m-m} that 
\[
m(G,\alpha)=\nu_0/2\pi.
\]
See~\cite[\S4.2]{Nas-Siam1} for details. 
Several numerical examples using this integral equation method are presented in~\cite{knv}.

\begin{example}[Circular domain]\label{ex:redm-mul}
We consider the multiply connected domain $G$ of connectivity $5$ in the interior of the unit circle and in the exterior of the circles $|z|=0.25$, $|z-0.6|=0.2$, $|z-(-0.3+0.5\i)|=0.2$, and $|z-(-0.3-0.5\i)|=0.2$.
We define the function $v(x,y)$ for all $x$ and $y$ such that $x+\i y\in G$ by
\[
v(x,y) = m(G,x+\i y).
\]
We use the integral equation method with $n=2^{11}$ to compute the values of the function $v(x,y)$. The level curves for the function $v(x,y)$ corresponding to the several levels are shown in Figure~\ref{fig:redm-mulL}. Each of these level curves describes the locations of the point $\alpha$ for which the values $m(G,\alpha)$ are a constant. It is clear from this figure that there are three locations of the point $\alpha$ at which $m(G,\alpha)$ has local maximum.
\end{example}

\begin{figure}[ht] %
\centerline{
		\scalebox{0.6}{\includegraphics[trim=0cm 0.5cm 0cm 0.5cm,clip]{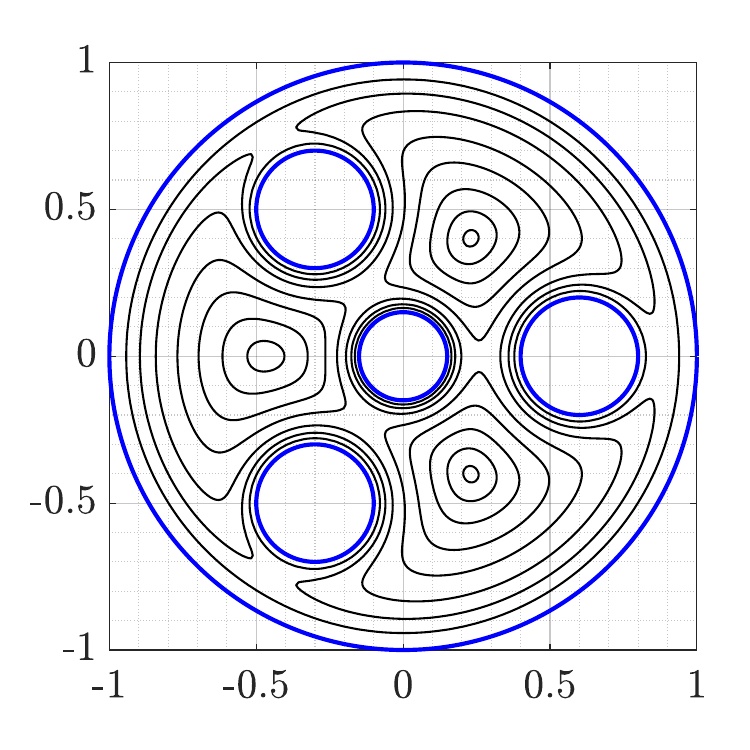}}		
	}
\caption{The contour lines of the function $v(x,y)$ in Example~\ref{ex:redm-mul}.}
	\label{fig:redm-mulL}
\end{figure}


\section{Moduli of quadrilaterals}\label{sec:sim-quad}

\subsection{Quadrilaterals}

For a given bounded simply connected domain $G$ and for a quadruple $ \{z_1 ,z_2 ,z_3 ,z_4 \}$ of its boundary points, we call $(G;z_1 ,z_2 ,z_3 ,z_4 )$ a
\textit{quadrilateral} if the points $z_1 ,z_2 ,z_3 ,z_4$ occur in this order when the boundary curve is traversed
in the positive direction. The points $z_1 ,z_2 ,z_3 ,z_4$ are called the vertices, and the
part of the oriented boundary between two successive vertices such as $z_1$ and $z_2$ is
called a boundary arc  and denoted $(z_1 ,z_2 ).$ 
By Riemann's mapping theorem, there is a conformal mapping $w=\psi(z)$ of $G$ onto a rectangle $R$ with vertices $0,1,1+h\i,h\i$, $h>0$, such that~\cite[p.52]{ps10},~\cite{lv}
\begin{equation}\label{eq:quad-bd}
\psi(z_1) = 0, \quad \psi(z_2)=1, \quad \psi(z_3)=1+h\i, \quad \psi(z_4)=h\i.
\end{equation}
Then the value $h$ is called the \emph{conformal modulus} of $G$:
\[
h =\moda(G;z_1 ,z_2 ,z_3 ,z_4 ) \equiv \M(\Delta([0,1],[h\i, 1+h\i]; G)) .
\]

An alternative method to find the modulus is to solve  the 
following Dirichlet-Neumann boundary value problem for the Laplace equation \cite{ah}. Suppose that $\partial G = \cup_{k=1}^4\partial G_k$; all the four boundary arcs $\partial G_k$ between vertices are
assumed to be non-degenerate. This problem is
\begin{equation}\label{eqn:dirichlet}
\left\{\begin{matrix}
    \Delta u& =&\ 0,& \text{on}\ &{\ G,} \\
    \partial u/\partial n&  =&\ 0,& \text{on}\ &{\partial G_1 = (z_1,z_2),}\\
    u& =&\ 1,&  \text{on}\ &{\partial G_2 = (z_2,z_3),}\\
    \partial u/\partial n&  =&\ 0,& \text{on}\ &{\partial G_3 = (z_3,z_4),}\\
    u& =&\ 0, & \text{on}\ &{\partial G_4 = (z_4,z_1).}\\
\end{matrix}\right.
\end{equation}
In terms of a solution function $u$ to the above problem,  the modulus can be
computed  as 
\begin{equation}\label{eq:q-mod}
h =\iint_G |\nabla u|^2 dm\,.
\end{equation}

It is an obvious fact that
\begin{equation}\label{recip}
\moda(G;z_1 ,z_2 ,z_3 ,z_4 ) \moda(G;z_4,z_1 ,z_2 ,z_3  )=1 \,.
\end{equation}
In the case of numerical computations, the difference
\begin{equation}\label{recError}
|\moda(G;z_1 ,z_2 ,z_3 ,z_4 ) \moda(G;z_4,z_1 ,z_2 ,z_3  )-1|.
\end{equation}
can be used as an experimental error characteristic if no
other error estimates are available. This method was used
in \cite{hvv} and as later work in \cite{hrv1,hrv2,hrv3} showed
it is often compatible with other error estimates.

\subsection{The numerical method}

To compute the modulus $h$ of the quadrilateral $(G;z_1,z_2,z_3,z_4)$, we first compute the conformal mapping $\hat w=\Phi_1(z)$ from the simply connected domain $G$ onto the unit disk $\B^2$ such that $\Phi_1(\hat\alpha)=0$ and $\Phi'_1(\hat\alpha)>0$ for some point $\hat\alpha$ in $G$. 
The mapping function $\hat w=\Phi_1(z)$ maps the positively oriented points $z_1,z_2,z_3,z_4$ on $\partial G$ onto positively oriented points $\hat w_1,\hat w_2,\hat w_3,\hat w_4$ on $\partial\B^2$. An exact formula for computing the modulus of the quadrilateral $(\B^2;\hat w_1,\hat w_2,\hat w_3,\hat w_4)$ is known in literature~\cite[(2.6.1)]{ps10}. Thus, by the conformal invariance of the modulus, we have
\[
h=\moda(G;z_1,z_2,z_3,z_4) = \moda(\B^2;\hat w_1,\hat w_2,\hat w_3,\hat w_4)
= \frac{2}{\pi}\mu\left(1/\sqrt{k}\right)
\]
where $k$ is given by the absolute (cross) ratio
\begin{equation}\label{eq:ration}
k=|\hat w_1,\hat w_2,\hat w_3,\hat w_4| = \frac{|\hat w_1-\hat w_3||\hat w_2-\hat w_4|}{|\hat w_1-\hat w_2||\hat w_2-\hat w_4|}.
\end{equation}

By the definition of $h=\moda(\B^2;\hat w_1,\hat w_2,\hat w_3,\hat w_4)$, there exists a conformal mapping
\[
w=\Phi_2(\hat w)
\]
from the unit disk $\B^2$ onto the rectangle $R$ such that
\[
\Phi_2(\hat w_1)=0,\;\Phi_2(\hat w_2)=1,\;\Phi_2(\hat w_3)=1+h\i,\;
\Phi_2(\hat w_4)=h\i.
\]
To compute such a conformal mapping $\Phi_2$, we first compute
the unique conformal mapping
\[
\tilde w=\Psi_1(w)
\]
from the domain $R$ onto the unit disk $\B^2$ with the normalization
\begin{equation}\label{eq:Phi-cond-b}
\Psi_1(\tilde\alpha)=0, \quad \Psi_1'(\tilde\alpha)>0
\end{equation}
where $\tilde\alpha$ is an auxiliary point in $R$, say $\tilde\alpha=(1+\i h)/2$. This conformal mapping $\Psi_1$ can be computed by the method presented in \S\ref{sec:sim-cm}. The mapping function $\tilde w=\Psi_1(w)$ maps the vertices $0,1,1+\i h,\i h$ of $\partial R$ onto four points  $\tilde w_1,\tilde w_2,\tilde w_3,\tilde w_4\in\partial\B^2$. These points are in general different from the points
$\hat w_1,\hat w_2,\hat w_3,\hat w_4$. Let
\[
\hat w=\Psi_2(\tilde w)= \hat w_3+
\frac{(\hat w_3-\hat w_1)(\hat w_2-\hat w_3)(\tilde w_2-\tilde w_1)(\tilde w-\tilde w_3)}{(\hat w_2-\hat w_1)(\tilde w_2-\tilde w_3)(\tilde w-\tilde w_1)-(\hat w_2-\hat w_3)(\tilde w_2-\tilde w_1)(\tilde w-\tilde w_3)}\,,
\]
then $\Psi_2$ maps the unit disk $\B^2$ onto itself such
that $\Psi_2(\tilde w_1)=\hat w_1$,
$\Psi_2(\tilde w_2)=\hat w_2$, and $\Psi_2(\tilde w_3)=\hat w_3$.
Thus, the function
\[
w=(\Psi_1^{-1}\circ\Psi_2^{-1})(\hat w)
\]
maps the unit disk $\B^2$ onto the rectangle $R$ and takes
the three points $\hat w_1$, $\hat w_2$, $\hat w_3$ to the three points $0$, $1$, $1+h\i$, respectively.
Since the function $\Phi_2$ is also a conformal mapping from the unit disk $\B^2$ onto the rectangle $R$ and maps the three points $\hat w_1,\hat w_2,\hat w_3$ to the three points $0$, $1$, $1+h\i$, respectively, then we have
\[
\Phi_2 = \Psi_1^{-1}\circ\Psi_2^{-1}.
\]
This is due to the fact of the uniqueness of the conformal mapping that maps the unit disk $\B^2$ onto the domain $R$ and maps three points on $\partial \B^2$ to three points on $\partial R$ when $h$ is fixed.
Thus, the function
\[
w=\psi(z)=(\Phi_2\circ\Phi_1)(z)=(\Psi_1^{-1}\circ\Psi_2^{-1}\circ\Phi_1)(z)
\]
is the required unique conformal mapping from the simply connected domain $G$ onto the rectangle $R=\{w\,:\, 0<\Re w<1,\,0<\Im w<h\}$ which satisfies~\eqref{eq:quad-bd} and the harmonic function $u(z)=\Re[\psi(z)]$ is then the unique solution of the boundary value problem~\eqref{eqn:dirichlet}.

\begin{example}[Trapezoid:~{\cite{hnv0,ps10}}]\label{ex:mod-trap}
Consider the trapezoid $T$ with the vertices $z_1=0$, $z_2=1$, $z_3=1+\i L$, $z_4=\i(L-1)$. The exact value of the modulus $\moda(T; z_1,z_2,z_3,z_4)$ is given for $L>1$ by~\cite[p.~82]{ps10}
\begin{equation}\label{eqn:mod-trap}
\moda(T; z_1,z_2,z_3,z_4)=\frac{\pi}{2\mu({ k})},  
\end{equation}
where
\[
{ k}=\frac{1-2\lambda\lambda'}{1+2\lambda\lambda'}, \quad
\lambda=\mu^{-1}\left(\frac{\pi}{2(2L-1)}\right),\quad \lambda'=\sqrt{1-\lambda^2}.
\]

The approximate values of the modulus $\moda(T; z_1,z_2,z_3,z_4)$ has been computed in~\cite{hnv0} for several values of $L$. 
The relative error in the approximate value of $\moda(T; z_1,z_2,z_3,z_4)$ computed using the above method for $L=1.5$ with $n=2^{12}$ is $5.47\times 10^{-14}$. 
The level curves of the potential function $u$ are shown in Figure~\ref{fig:trap-gear} (left).
\end{example}

\begin{figure}[ht] %
\centerline{
		\scalebox{0.6}{\includegraphics[trim=0cm 0.5cm 0cm 0.5cm,clip]{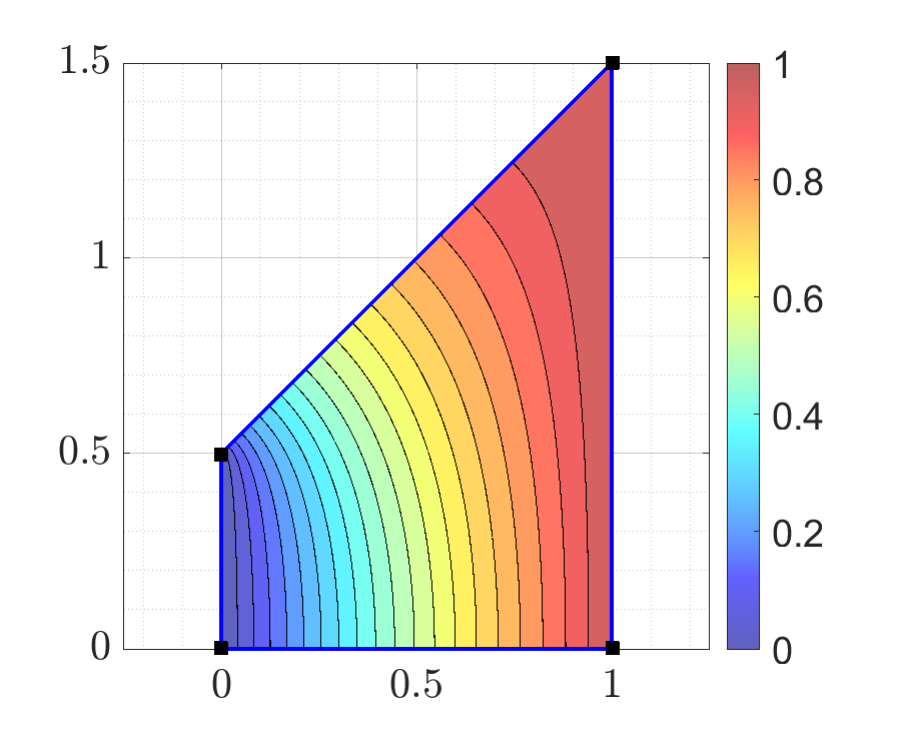}}
		\hfill
		\scalebox{0.6}{\includegraphics[trim=0cm 0.5cm 0cm 0.5cm,clip]{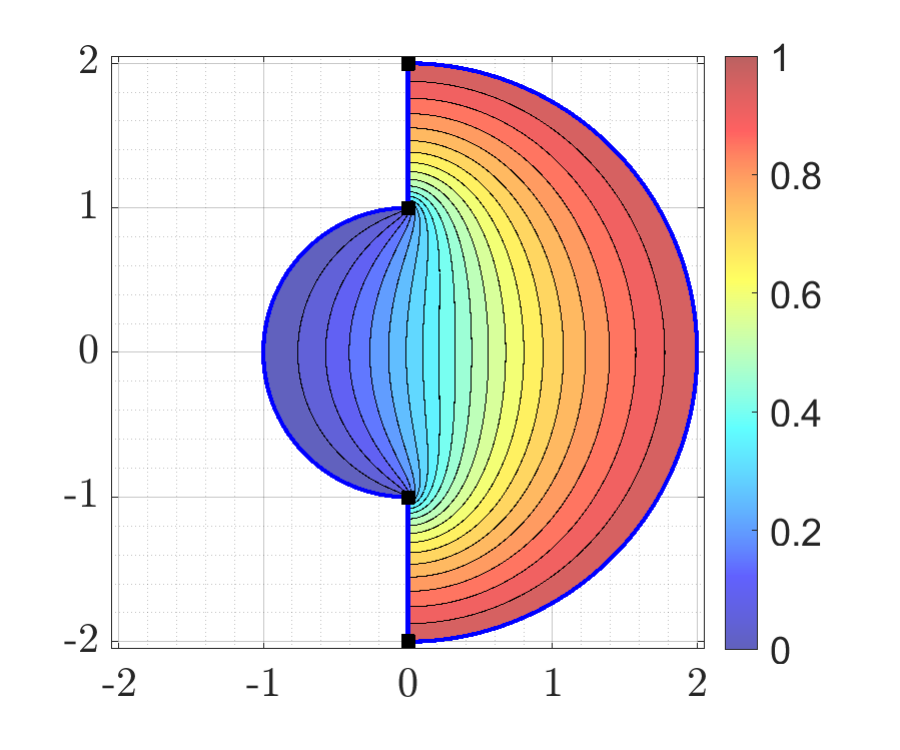}}		
	}
\caption{The level curves of the potential function $u$ for the trapezoid in Example~\ref{ex:mod-trap} with $L=1.5$ (left) and the gear domain in Example~\ref{ex:mod-gear} (right).}
	\label{fig:trap-gear}
\end{figure}

\begin{example}[Gear domain]\label{ex:mod-gear}
We consider a gear domain $G$ with one tooth which is a polycircular domain whose boundary consists of the segment $[-\i,-2\i]$, the semicircular arc connecting the points $-2\i$ and $2\i$ on the right half  plane, the segment $[2\i,\i]$, and the semicircular arc connecting the points $\i$ and $-\i$ on the left half  plane~\cite{hnv0}.
The approximate values of the modulus $\moda(G; -\i,-2\i,2\i,\i)$ are computed using the above method with $n=2^{12}$ and the obtained approximate value is $1.76445071147738$. 
The level curves of the potential function $u$ are shown in Figure~\ref{fig:trap-gear} (right).
\end{example}

\begin{example}[{Amoeba-shaped domain}]\label{ex:mod-amoeba}
Consider the simply connected domain $G$ in the interior of the curve $\Gamma$ (amoeba-shaped boundary) with the parametrization
\[
\eta(t) = \left(e^{\cos t}\cos^2(2t)+e^{\sin t}\sin^2(2t)\right)e^{\i t}, \quad 0\le t\le 2\pi.
\]
The approximate values of the modulus $\moda(G; e,\i,-e^{-1},-\i)$ are computed using the above method with $n=2^{10}$ and the obtained approximate value is $1.20089247845316$. 
The level curves of the potential function $u$ are shown in Figure~\ref{fig:amoeba}.
\end{example}

\begin{figure}[ht] %
\centerline{
		\scalebox{0.6}{\includegraphics[trim=0cm 0.5cm 0cm 0.5cm,clip]{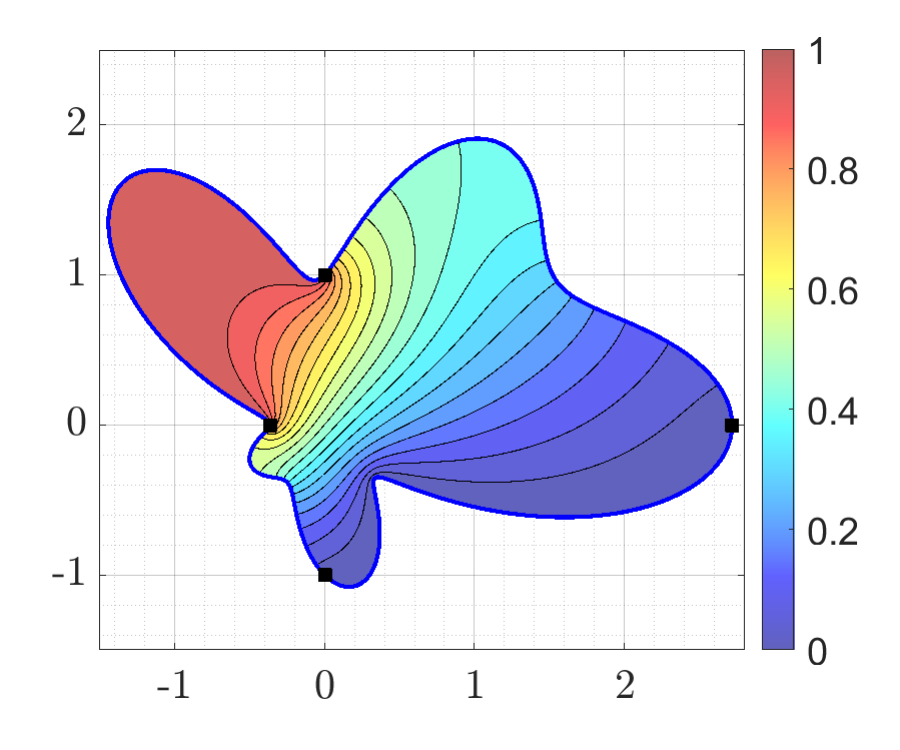}}
	}
\caption{The level curves of the potential function $u$ for the amoeba-shaped domain in Example~\ref{ex:mod-amoeba}.}
	\label{fig:amoeba}
\end{figure}

\subsection{Convex polygonal quadrilateral}

A look at the literature \cite{ky} shows that there are very few non-trivial domains for which a quadrilateral can be analytically handled. In the following theorem, an analytic formula  is given for the modulus of a bounded convex polygonal quadrilateral $G$ in the upper-half plane with vertices $0,1,A,B$. 

\begin{theorem}[\cite{hvv}]\label{thm:md-conv} 
Choose $a,b,c$ such that $0<a,b<1$ and $\max\{a+b,1\}\leq c\leq 1+\min\{a,b\}$. Let $G$ be a polygonal quadrilateral in the upper half-plane with interior angles $b\pi$, $(c-b)\pi$, $(1-a)\pi$ and $(1+a-c)\pi$ at the vertices $0,1,A,B$, respectively. Then the conformal modulus of $(G;0,1,A,B)$ is given by
\begin{align*}
\moda(G;0,1,A,B)=\frac{2}{\pi}\mu(r),    
\end{align*}
where $0<r<1$ fulfills the equation
\begin{align}
A-1=\frac{L(1-r^2)^{c-a-b}F(c-a,c-b;c+1-a-b;1-r^2)}{F(a,b;c;r^2)}
\end{align}
with
\begin{equation}\label{eq:L}
L=\frac{B(c-b,1-a)}{B(b,c-b)}e^{(b+1-c)\pi i} 
\end{equation}
where $B(z,w)= \Gamma(z) \Gamma(w)/\Gamma(z+w)$ is the beta function.
\end{theorem}

A Mathematica code for the implementation of the method in Theorem  \ref{thm:md-conv} is presented in~\cite{hvv} and a MATLAB code is presented in~\cite{nrv}.

\begin{example}[Convex domain]\label{ex:mod-conv}
Consider the polygonal domain $G$ with the vertices $0,1,x+\i y,\i$ with $x>0$, $y>0$, and $x+y>1$. We define the function 
\[
v(x,y) = \moda(G; 0,1,x+\i y,\i)
\]
for $0.1\le x\le 3$, $0.1\le y\le 3$, and $x+y>1$. 
We use the presented numerical method with $n=2^{12}$ to compute approximate values $v_n(x,y)$ of the function $v(x,y)$. 
The level curves of the function $v_n(x,y)$ are shown in Figure~\ref{fig:md-conv} (left). The relative error $|v(x,y)-v_n(x,y)|/|v(x,y)|$ in the computed values $v_n(x,y)$ is shown in Figure~\ref{fig:md-conv} (right) where $v(x,y)$ are the exact values computed using the method presented in Theorem~\ref{thm:md-conv}. It is clear from Figure~\ref{fig:md-conv} that $v(x,y)=1$ for $y=x$, $v(x,y)<1$ for $y<x$, and $v(x,y)>1$ for $y>x$. 
\end{example}

\begin{figure}[ht] %
\centerline{
		\scalebox{0.6}{\includegraphics[trim=0cm 0.5cm 0cm 0.5cm,clip]{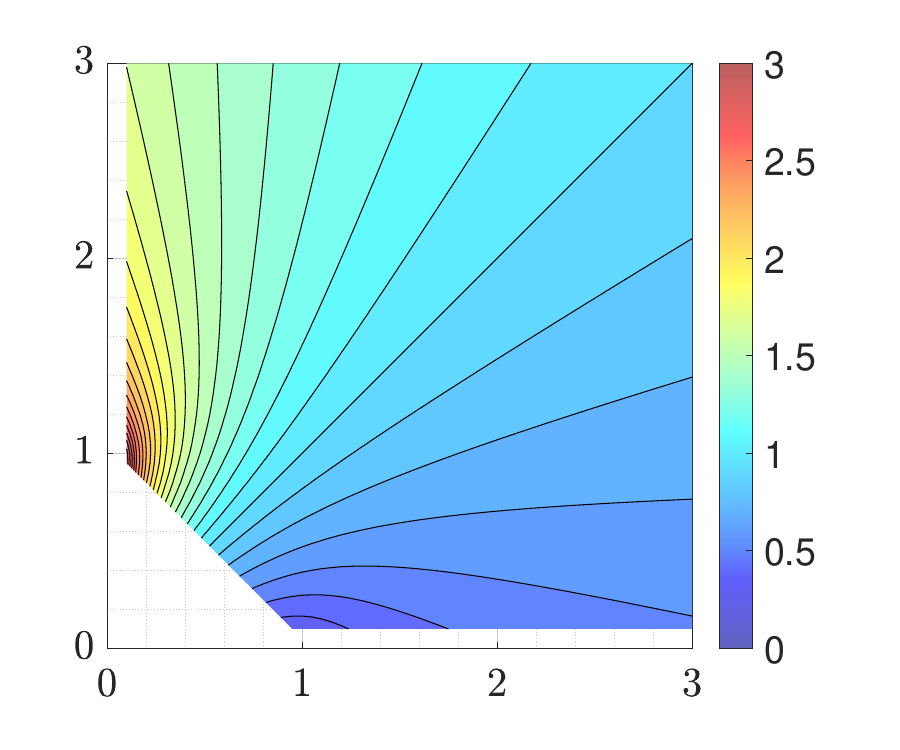}}
		\hfill
		\scalebox{0.6}{\includegraphics[trim=0cm 0.5cm 0cm 0.5cm,clip]{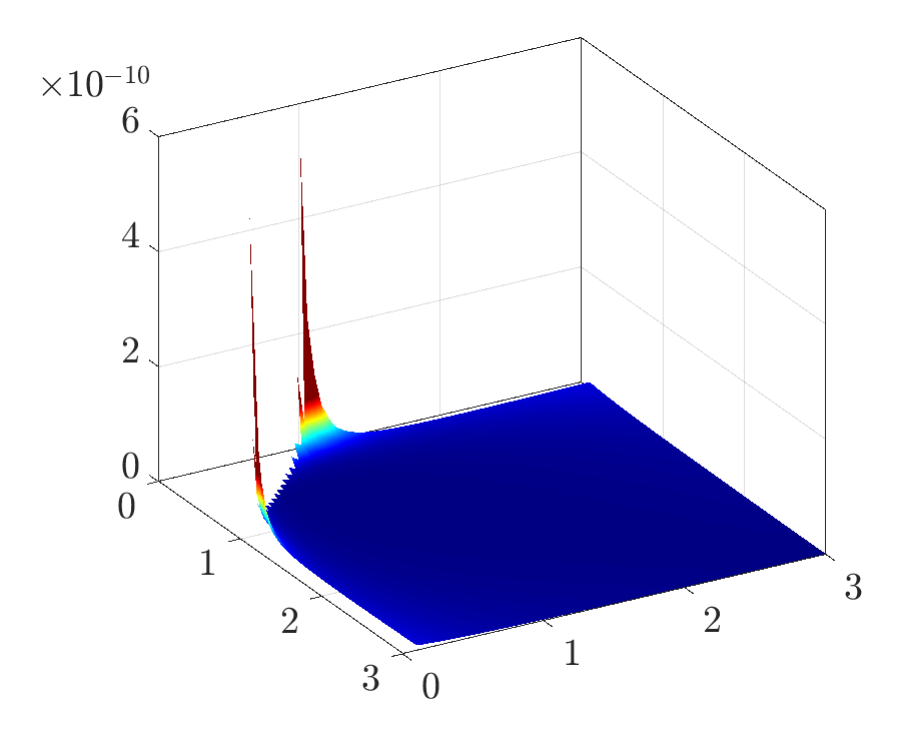}}		
	}
\caption{The level curves of the function $v$ (left) and the relative error (right) in Example~\ref{ex:mod-conv}.}
	\label{fig:md-conv}
\end{figure}

\subsection{Unbounded quadrilaterals}\label{sec:mod-ext}

The above presented method can be extended to compute the modulus of unbounded quadrilaterals, i.e., when $G$ is an unbounded simply connected domain~\cite{nnv}. Because harmonic
functions satisfy the maximum principle, there is a unique level set of $u$ corresponding to the level value $u(\infty)$,
which passes through $\infty$.
Numerical methods for computing the value of the potential function at $\infty$ are given in \cite{hrv2,nnv}.
In  \cite{nnv},  for the polygonal quadrilateral of Theorem \ref{thm:md-conv}, 
the value of the potential function $u(\infty)$ was found also in terms of an algorithm for the exterior modulus, implemented in~\cite{nsv}.

In the following example, we consider an unbounded quadrilateral with known exact modulus.
For the numerical implementation of the boundary integral equation method to compute the modulus of unbounded quadrilateral, see~\cite{nnv,nrrvwyz}.

\begin{example}[The exterior of a rectangle]
Let $R$ be the unbounded simply connected domain in the exterior of the rectangle with the vertices $1$, $0$, $\i k$, $1 +\i k$, $k > 0$. 
Then the exterior modulus $h$ of the quadrilateral $(R;1,0,\i k,1+\i k)$ can be expressed as
\[
h=\frac{1}{\pi}\mu(\psi^{-1}(1/k)),
\]
where the function $\psi(r)$ is defined by Duren and Pfaltzgraff~\cite{dp} as
\[ 
\psi(r) = \frac{2( \E(r)-(1- r) \K(r))}{
\E(r')- r \K(r')},\quad r'=\sqrt{1- r^2}.
\]

Using the  numerical method   presented above with $n=2^{12}$, the obtained approximate value of the modulus $\moda(R; 1,0,\i k,1+\i k)$ for $k=2$ is $1.154924858699863$ and the relative error in this approximate value is $1.31\times 10^{-13}$.
\end{example}

\section{Harmonic measure}\label{sec:hm}

\subsection{Harmonic measure for simply connected domains}\label{sec:sim-hm}

The harmonic measure is one of the key notions of potential theory and it has numerous applications to geometric function theory \cite{garmar}.
Assume that $G$ is a simply connected domain, i.e., $m=0$ in Section~\ref{sec:bie}. We assume that $G$ is bordered by the piecewise smooth Jordan curve $\Gamma=\Gamma_0$ which is parametrized by the function $\eta(t)=\eta_0(t)$ in~\eqref{eq:mul-eta}. 
Let $L$ be a nonempty subset of $\Gamma$ such that $\Gamma\backslash L\ne\emptyset$. 
The harmonic measure of $L$ with respect to $G$ is the $C^2(G)$ function $u: G \to (0,1)$ satisfying the Laplace equation
\[
\Delta u = 0
\]
in $G$ and the boundary conditions $u(z)=1$ when $z \in L $ and $u(z)=0$ when $z\in\Gamma\setminus L$. For unbounded $G$, we assume that $u(z)$ is bounded at $\infty$.
The harmonic measure of $L$ with respect to $G$ will be denoted by $\omega(z,L,G)$ (see e.g.,~\cite[p.~123]{avv}, \cite[Ch I]{garmar}, and~\cite[p.~111]{tsu}).

We assume that $L$ is a union of $\ell$ nonempty connected and disjoint arcs $L_1,\ldots,L_\ell$ such that the two end points of $L_j$ are $z_{2j-1}$ and $z_{2j}$, $j=1,\ldots,\ell$. 
That is, the value of the harmonic measure $\omega(z,L,G)$ is $1$ when $z$ is on the arc between the two points $z_{2j-1}$ and $z_{2j}$, and $0$ when $z$ is on the arc between the two points $z_{2j}$ and $z_{2j+1}$, $j=1,\ldots,\ell$, where $z_{2\ell+1}=z_{1}$.
We assume that $z_1,z_2,\ldots,z_{2\ell}$ are oriented in the direction of $\Gamma$.
Let $w=\Phi(z)$ be the unique conformal mapping from $G$ onto the unit disk $\B^2$ (with the normalization~\eqref{eq:Phi-cond} or~\eqref{eq:Phi-condu}). Then $\Phi$ maps the boundary $\Gamma$ onto the unit circle and maps the point $z_j$ on $\Gamma$ onto the point $w_j=\Phi(z_j)$ on the unit circle, $j=1,\ldots,2\ell$. The points $w_1,w_2,\ldots,w_{2\ell}$ are oriented counterclockwise. 

Let $\theta_j=\arg(w_j)$ such that $\beta<\theta_1<\theta_2<\cdots<\theta_{2\ell}<\beta+2\pi$ for some real number $\beta$. 
Then the Mobi\"us transformation 
\[
\zeta = \Psi(w)=\i\,\frac{e^{\i\beta}+w}{e^{\i\beta}-w}
\]
maps the unit disk $|w|<1$ onto the upper half plane $\Im\,\zeta>0$ and maps the unit circle onto the real line such that $\Psi(e^{\i\beta})=\infty$. It maps also the points $w_1,w_2,\ldots,w_{2\ell}$ on the unit circle onto the points $x_1,x_2,\ldots,x_{2\ell}$ on the real line such that $-\infty<x_1<x_2<\cdots<x_{2\ell-1}<x_{2\ell}<\infty$. It is clear that the function~\cite[p.~421]{zs}
\begin{equation}\label{eq:hm-ps}
\psi(\zeta) = \frac{1}{\pi}\sum_{j=1}^{2\ell}(-1)^{j}\arg(\zeta-x_j)
\end{equation}
is harmonic in the upper half plane $\Im\,\zeta>0$, $\psi(x)=1$ when $x$ is real with $x_{2j-1}<x<x_{2j}$ for $j=1,2,\ldots,\ell$, and $\psi(x)=0$ when $x$ is real with $x<x_1$, $x>x_{2\ell}$, or $x_{2j}<x<x_{2j+1}$ for $j=1,2,\ldots,\ell-1$. The branch of $\arg$ is chosen such that $\arg(1)=0$.

Since $w_j=e^{\i\theta_j}$ for $j=1,2,\ldots,2\ell$, then
\[
x_j=\Psi(w_j)=\i\,\frac{e^{\i\beta}+e^{\i\theta_j}}{e^{\i\beta}-e^{\i\theta_j}}
=\i\,\frac{e^{\i(\beta-\theta_j)/2}+e^{-\i(\beta-\theta_j)/2}}{e^{\i(\beta-\theta_j)/2}-e^{-\i(\beta-\theta_j)/2}}
=\cot\left(\frac{\beta-\theta_j}{2}\right).
\]
Note also that
\[
\zeta = \i\,\frac{e^{\i\beta}+w}{e^{\i\beta}-w} 
= \i\,\frac{e^{\i\beta}+\Phi(z)}{e^{\i\beta}-\Phi(z)}
= \i\,\frac{1+e^{-\i\beta}\Phi(z)}{1-e^{-\i\beta}\Phi(z)}.
\]
Then, it follows from~\eqref{eq:hm-ps} that the harmonic measure $\omega(z,L,G)$ is given by
\begin{equation}\label{eq:sim-hm-d}
\omega(z,L,G) = \frac{1}{\pi}\sum_{j=1}^{2\ell}(-1)^{j}\arg\left(\i\frac{1+e^{-\i\beta}\Phi(z)}{1-e^{-\i\beta}\Phi(z)}+\cot\left(\frac{\theta_j-\beta}{2}\right)\right).
\end{equation} 
The conformal mapping $\Phi(z)$ will be computed by the method presented in~\S\ref{sec:sim-cm}. 

\subsection{Harmonic-measure distribution function}

For a given basepoint $z_0$ in a given simply connected domain $G$, the harmonic-measure distribution function $h(r)$, or the $h$-function, is the piecewise smooth, non-decreasing function, $h:[0,\infty)\to[0,1]$, defined by
\[
h(r)=\omega(z_0,\Gamma\cap\overline{B(z_0,r)},G),
\]
i.e., $h(r)$ is the value of the harmonic measure of $\Gamma\cap\overline{B(z_0,r)}$ with respect to $G$ at the point $z_0$.
The $h$-function was first studied in depth by Walden \& Ward~\cite{WaWa}. An overview of the main properties of $h$-functions is given in the survey paper~\cite{SnWa}.

The value of $h(r)$ is the probability that a Brownian particle reaches a boundary component of $G$ within a certain distance $r$ from its point of release $z_0$.
The properties of two-dimensional Brownian motions were investigated extensively by Kakutani~\cite{kak} who found a deep connection between Brownian motion and harmonic functions (see also~\cite{SnWa}).  Han-Rasila-Sottinen \cite{hrs} also studied harmonic measure using random walk simulations.

Explicit formulas of $h$-functions for several simply connected domains are given in~\cite{AMth,AMta,SnWa}. 
In~\cite{GSWC}, explicit formulas of $h$-functions for a class of multiply connected symmetrical slit domains were derived in terms of the Schottky-Klein prime function~\cite{c,CKGN}. A numerical method for computing the $h$-functions for such symmetric slit domains with high connectivity is presented in~\cite{GrNa24} (see also~\cite{GrNag}). The method is based on using the boundary integral equation with the generalized Neumann kernel.

\begin{example}[Domain interior to an ellipse]\label{ex:hfun-in}
We consider the simply connected domain $G$ in the interior of the ellipse
\[
\eta(t)=\cosh(\tau+\i t), \quad 0\le t\le 2\pi, 
\]
for $\tau=0.5$ and $z_0=0$.

The values of $h(r)$ are computed using the above method with $n=2^{12}$ and the obtained results are presented in Figure~\ref{fig:hfunction} (left). It is clear that $h(r)=0$ for $r\le\sinh(\tau)\approx0.5211$ and $h(r)=1$ for $r\ge\cosh(\tau)\approx1.1276$.  
\end{example}

\begin{example}[Domain exterior to an ellipse]\label{ex:hfun-ex}
Consider the simply connected domain $G$ in the exterior of the ellipse
\[
\eta(t)=\cosh(\tau-\i t), \quad 0\le t\le 2\pi, 
\]
for $\tau=0.9$ and $z_0=1+\cosh \tau$.

The values of $h(r)$ is computed using the above method with $n=2^{12}$ and the obtained results are presented in Figure~\ref{fig:hfunction} (right). It is clear that $h(r)=0$ for $r\le1$ and $h(r)=1$ for $r\ge1+\cosh(\tau)\approx3.8662$ (see also~\cite[Figure~5.7]{AMth}).
\end{example}

\begin{figure}[ht] %
\centerline{
		\scalebox{0.6}{\includegraphics[trim=0cm 0.5cm 0cm 0.5cm,clip]{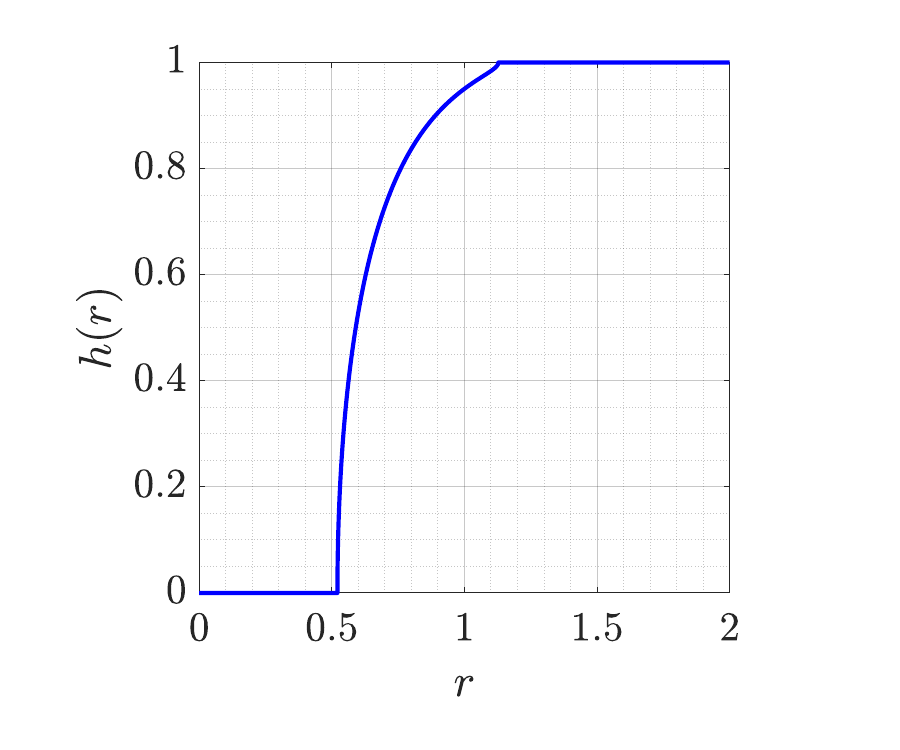}}
		\hfill
		\scalebox{0.6}{\includegraphics[trim=0cm 0.5cm 0cm 0.5cm,clip]{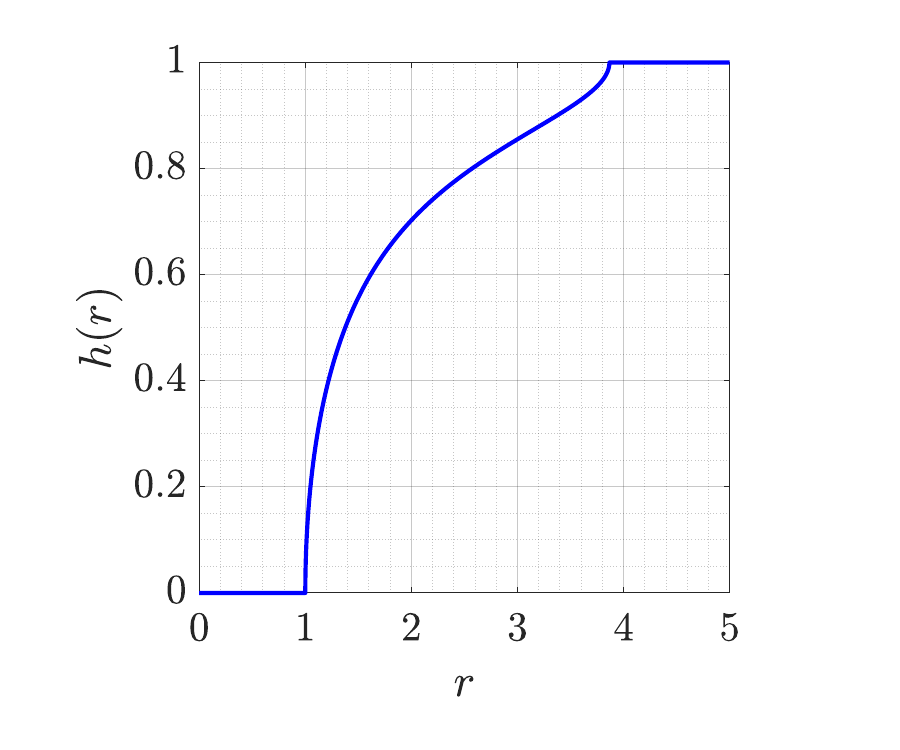}}		
	}
\caption{The graphs of the $h$-function $h(r)$ for Example~\ref{ex:hfun-in} (left) and Example~\ref{ex:hfun-ex} (right).}
	\label{fig:hfunction}
\end{figure}

\subsection{Harmonic measure for multiply connected domains}\label{sec:mul-hm}

Let $G$ be a given bounded multiply connected domain of connectivity $m+1$, let $\Gamma=\partial G=\cup_{k=0}^m\Gamma_k$ be parametrized by the function $\eta(t)$ in~\eqref{eq:mul-eta}, and let the function $A(t)$ be defined by~\eqref{eq:A}, i.e., $A(t)=\eta(t)-\alpha$. Let also the kernels $N(s,t)$ and $M(s,t)$ of the integral operators $\bN$ and $\bM$, respectively, be formed with these functions $\eta(t)$ and $A(t)$.

For $k=0,1,\ldots,m$, let $\omega(z,\Gamma_k,G)$ be the harmonic measure of $\Gamma_k$ with respect to $G$, i.e., $\omega(z,\Gamma_k,G)=\sigma_k(z)$ where $\sigma_k(z)$ is the unique solution of the Dirichlet problem:
\begin{subequations}\label{eq:bdv-sig}
	\begin{align}
		\label{eq:sig-Lap}
		\nabla^2 \sigma_k(z) &= 0 \quad\quad \mbox{if }z\in G, \\
		\label{eq:sig-j}
		\sigma_k(z)&= \delta_{k,j} \quad \mbox{if }z\in \Gamma_j, \quad j=0,1,\ldots,m, 
	\end{align}
\end{subequations}
where $\delta_{k,j}$ is the Kronecker delta function. It is clear that $\sum_{k=0}^m\sigma_k(z)=1$ and hence
\[
\sigma_0(z)=1-\sum_{k=1}^m\sigma_k(z).
\]
Thus, it is enough to compute the functions $\sigma_1(z),\ldots,\sigma_m(z)$.

For each $j=1,2,\ldots,m$, let the function $\gamma_j(t)$ be defined by~\eqref{eq:gam-k}, let $\rho_j(t)$ be the unique solution of the integral equation~\eqref{eq:ie}, and let the piecewise constant function $\nu_j(t)=(\nu_{0,j}, \nu_{1,j}, \ldots,\nu_{m,j})$ be given by~\eqref{eq:h}. Then, it follows from~\S\ref{sec:cgc} that the function $f_j(z)$ with the boundary values
\[
A(t)f_j(\eta(t))=\gamma_j(t)+\nu_j(t)+\i\mu_j(t), \quad t\in J,
\]
is analytic in $G$. Then the function $\sigma_k(z)$ is given for $z\in G$ by
\begin{equation}\label{eq:F-sigm-k}
\sigma_k(z)=\Re[(z-\alpha)g_k(z)]+c_k-\sum_{j=1}^{m} a_{kj}\log|z-\alpha_j|, \quad k=1,\ldots,m,
\end{equation}
where 
\[
g_k(z) = \sum_{j=1}^{m} a_{kj} f_j(z)
\]
and the real constants $c_k$ and $a_{k,1},\ldots,a_{k,m}$ are computed by solving the uniquely solvable linear system
\begin{equation}\label{eq:sys-sigma}
	\left[\begin{array}{ccccc}
		\nu_{0,1}    &\nu_{0,2}    &\cdots &\nu_{0,m}      &1       \\
		\nu_{1,1}    &\nu_{1,2}    &\cdots &\nu_{1,m}      &1       \\
		\vdots       &\vdots       &\ddots &\vdots       &\vdots  \\
		\nu_{m,1}    &\nu_{m,2}    &\cdots &\nu_{m,m}      &1       \\
	\end{array}\right]
	\left[\begin{array}{c}
		a_{k,1}    \\a_{k,2}    \\ \vdots \\ a_{k,m} \\  c_k 
	\end{array}\right]
	= \left[\begin{array}{c}
		\delta_{k,0} \\  \delta_{k,1} \\  \vdots \\ \delta_{k,m}  
	\end{array}\right].
\end{equation}

By computing the functions $\rho_j$ and $\nu_j$, and the constants $c_k$ and $a_{k,j}$ for $j,k=1,\ldots,m$, we obtain approximations of the boundary values of the analytic function $g_k(z)$ by
\[
A(t)g_k(\eta(t)) = \sum_{j=1}^{m} a_{kj}(\gamma_j(t)+\nu_j(t)+\i\rho_j(t)), \quad t\in J.
\]
Then, the values of $g_k(z)$ for $z\in G$ can be computed by the Cauchy integral formula, and hence the values of $\sigma_k(z)$ can be computed for $z\in G$ by~\eqref{eq:F-sigm-k}.

\begin{example}[Circular domain]\label{ex:hm-mul}
We consider the multiply connected domain $G$ of connectivity $5$ in the interior of the unit circle and in the exterior of the circles $|z|=0.25$, $|z-0.6|=0.2$, $|z-(-0.3+0.5\i)|=0.2$, and $|z-(-0.3-0.5\i)|=0.2$.
We use the integral equation method with $n=2^{11}$ to compute the values of the functions $\sigma_1(z)=\omega(z,\Gamma_1,G)$ and $\sigma_2(z)=\omega(z,\Gamma_2,G)$. The level curves for these functions are shown in Figure~\ref{fig:hm-mulL}. 
\end{example}

\begin{figure}[ht] %
	\centering{
	\hfill{\includegraphics[width=0.45\textwidth]{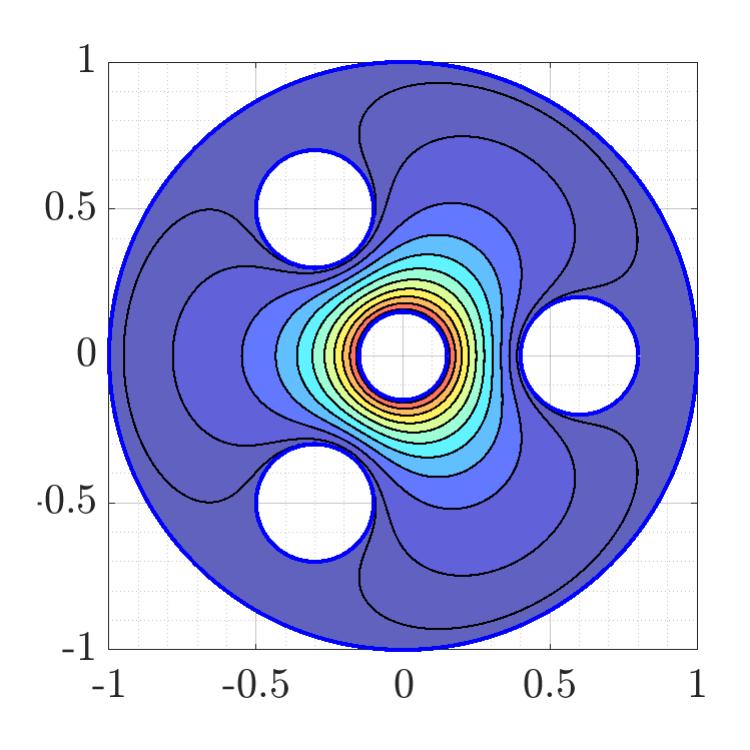}}
	\hfill{\includegraphics[width=0.45\textwidth]{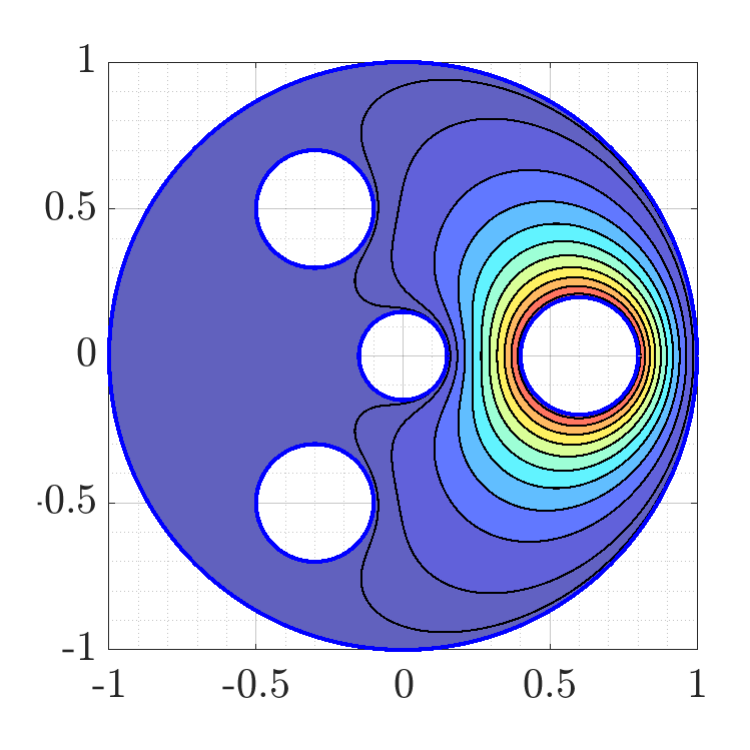}}
	\hfill}
\caption{The contour lines of the functions $\sigma_1(z)$ (left) and $\sigma_2(z)$ (right) in Example~\ref{ex:hm-mul}.}
	\label{fig:hm-mulL}
\end{figure}

\section{Hyperbolic distance}\label{sec:hd}

Assume that $G$ is a given bounded simply connected domain and $\alpha$ is a given point in $G$.
Let $w=\Phi(z)$ be the conformal mapping from the bounded simply connected domain $G$ onto the unit disk $\B^2$ with the normalization~\eqref{eq:Phi-cond}. 
By \eqref{rhodef} one can define the hyperbolic metric on the Jordan domain $G$ in terms of the conformal mapping function $w=\Phi(z)$  as follows:
 \[
 \rho_G(x,y) = \rho_{\B^2}(\Phi(x),\Phi(y))\,.
 \]
In hyperbolic geometry the boundary $\partial G$ has the same role as the point of $\{\infty\}$ in Euclidean geometry: both are like ``horizon's'', we cannot see beyond the
 horizon.
Due to conformal invariance, the hyperbolic geometry is more useful than the Euclidean geometry when studying the inner geometry of domains in geometric function theory.

\begin{example}\label{ex:hyp-disk}
We consider the simply connected domain $G$ in the interior of the curve $\Gamma$ with the parametrization
\[
\eta(t)=(2+\cos(5t))e^{\i t},\quad 0\le t\le 2\pi.
\]
Then, for a given point $z_0$ in $G$, we define the function $u(x,y)$ for all $x$ and $y$ such that $x+\i y\in G$ by
\[
u(x,y) = \rho_G(z_0,x+\i y).
\]
We use the method presented in~\S\ref{sec:sim-cm} with $n=2^{12}$ and $\alpha=z_0$ to compute the the conformal mapping $w=\Phi(z)$ with the normalization~\eqref{eq:Phi-cond}, and hence to compute the values of the function $u(x,y)$ in the domain $G$. We plot the contour lines for the function $u(x,y)$ corresponding to several levels in Figure~\ref{fig:hyp-dis} for $z_0=0$ (left) and $z_0=1$ (right). Each of these level curves is a hyperbolic circle in $G$ with center $z_0$.  

\end{example}

\begin{figure}[ht] %
\centerline{
		\scalebox{0.6}{\includegraphics[trim=0cm 0.5cm 0cm 0.5cm,clip]{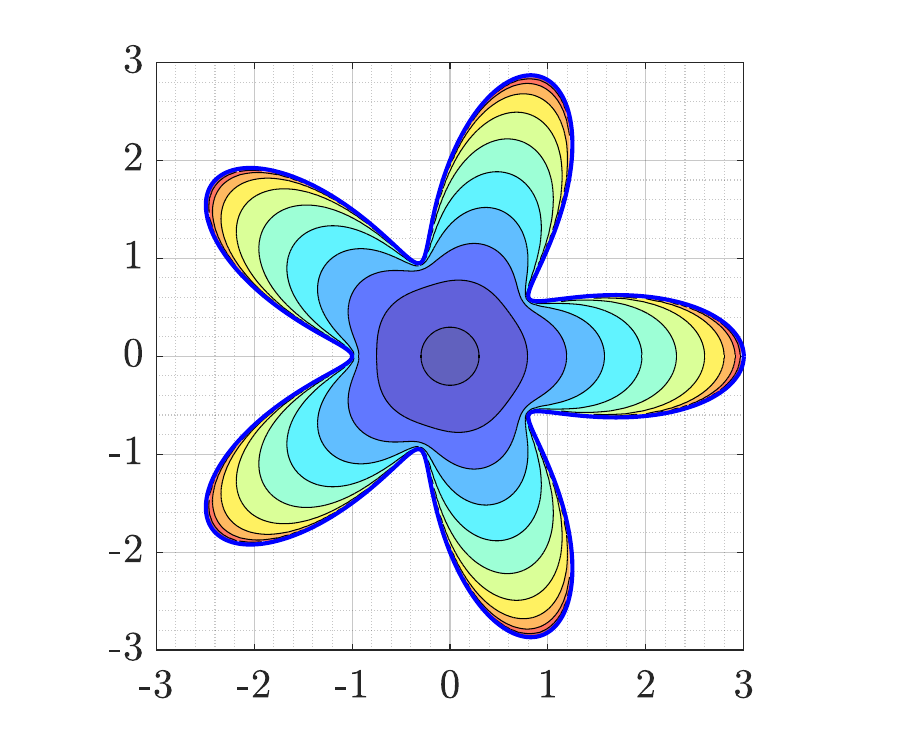}}
		\hfill
		\scalebox{0.6}{\includegraphics[trim=0cm 0.5cm 0cm 0.5cm,clip]{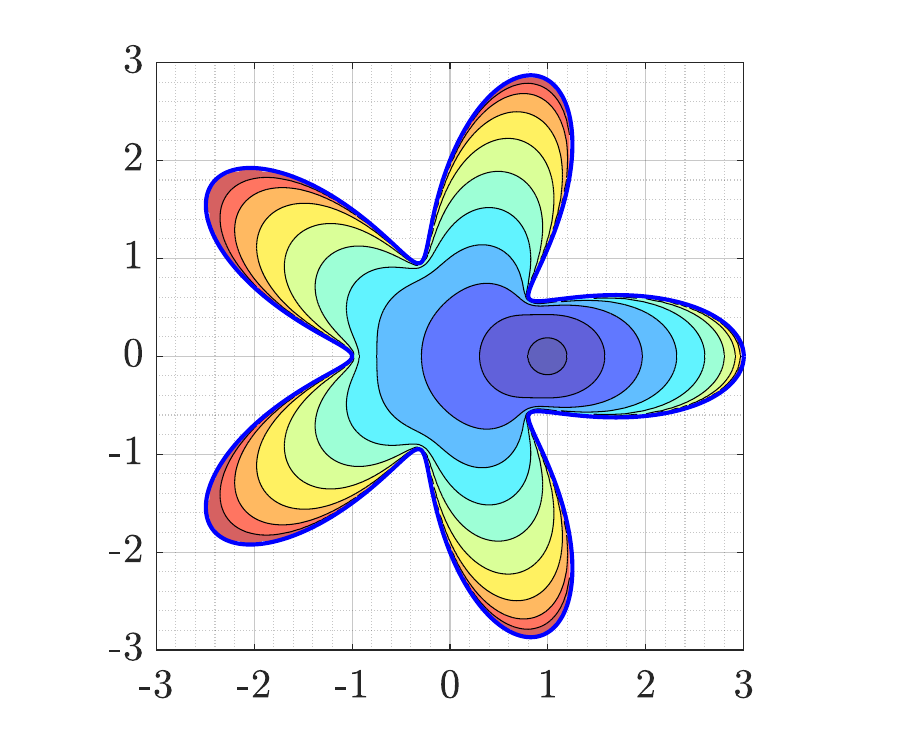}}		
	}
\caption{The contour lines of the function $u(x,y)$ in Example~\ref{ex:hyp-disk} for $z_0=0$ (left) and $z_0=1$ (right).}
	\label{fig:hyp-dis}
\end{figure}

\begin{example}\label{ex:hyp-geod}
We consider two simply connected polygonal domains $G$ as in Figure~\ref{fig:hyp-geod}.
For each of these two polygons, we use the method presented in~\S\ref{sec:sim-cm} with $n=3\times2^{14}$, $\alpha=0.5-1.5\i$ for the graph in the left and $n=7\times2^{13}$, $\alpha=-\i$ for the graph in the right to compute the the conformal mapping $w=\Phi(z)$ with the normalization~\eqref{eq:Phi-cond}. Then, we use this conformal mapping to compute the hyperbolic geodesics starting at the middle point of each boundary segment of these polygons, passing through the point $\alpha$ and continue to the boundary of the polygon. We plot the computed hyperbolic geodesics in Figure~\ref{fig:hyp-geod} where the the middle point of each boundary segment is marked with a red-dot and the point $\alpha$ is marked with a red square. In Tables~\ref{tab:hdis-L} and~\ref{tab:hdis-R}, we choose five points in the interior of the domain $G$ in Figure~\ref{fig:hyp-geod}(left) and Figure~\ref{fig:hyp-geod}(right), respectively, and compute the hyperbolic distance between all possible pair of points from these five points. 
\end{example}

\begin{table}[ht]
\caption{The hyperbolic distance $\rho_G(z_k,z_j)$ between the five points $z_1=-1+2.5\i$, $z_2=-1.5$, $z_3=0.5-1.5\i$, $z_4=2.5-0.5\i$, and $z_5=1+0.5\i$ in the interior of the domain $G$ in Figure~\ref{fig:hyp-geod}(left).}
\label{tab:hdis-L}
\centering
\begin{tabular}{cccccc} 
\hline
$k\backslash j$ & $1$ & $2$ & $3$ & $4$ & $5$ \\ \hline
$1$  &$0$               &$8.211957028453$  &$17.82130682562$  &$25.86706921787$  &$32.31461859625$ \\
$2$  &$8.211957028453$  &$0$               &$9.609696981382$  &$17.65545939223$  &$24.10300877075$ \\
$3$  &$17.82130682562$  &$9.609696981382$  &$0$               &$8.045817781269$  &$14.49336799127$ \\
$4$  &$25.86706921787$  &$17.65545939223$  &$8.045817781269$  &$0$               &$6.477277067630$ \\
$5$  &$32.31461859625$  &$24.10300877075$  &$14.49336799127$  &$6.477277067630$  &$0$              \\
		\hline
	\end{tabular}
\end{table}

\begin{table}[ht]
\caption{The hyperbolic distance $\rho_G(z_k,z_j)$ between the five points $z_1=-2+\i$, $z_2=\i$, $z_3=-\i$, $z_4=-2\i$, and $z_5=2+\i$ in the interior of the domain $G$ in Figure~\ref{fig:hyp-geod}(right).}
\label{tab:hdis-R}
\centering
\begin{tabular}{cccccc} 
\hline
$k\backslash j$ & $1$ & $2$ & $3$ & $4$ & $5$ \\ \hline
$1$  &$0$               &$7.522231288671$  &$5.153143098454$  &$5.612327759939$  &$10.30445568667$ \\
$2$  &$7.522231288671$  &$0$               &$3.103715005941$  &$4.095139688250$  &$7.522231288628$ \\
$3$  &$5.153143098454$  &$3.103715005941$  &$0$               &$0.991424682309$  &$5.153143098413$ \\
$4$  &$5.612327759939$  &$4.095139688250$  &$0.991424682309$  &$0$               &$5.612327759898$ \\
$5$  &$10.30445568667$  &$7.522231288628$  &$5.153143098413$  &$5.612327759898$  &$0$       \\
\hline
	\end{tabular}
\end{table}


\begin{figure}[ht] %
\centerline{
		\scalebox{0.6}{\includegraphics[trim=0cm 0.5cm 0cm 0.5cm,clip]{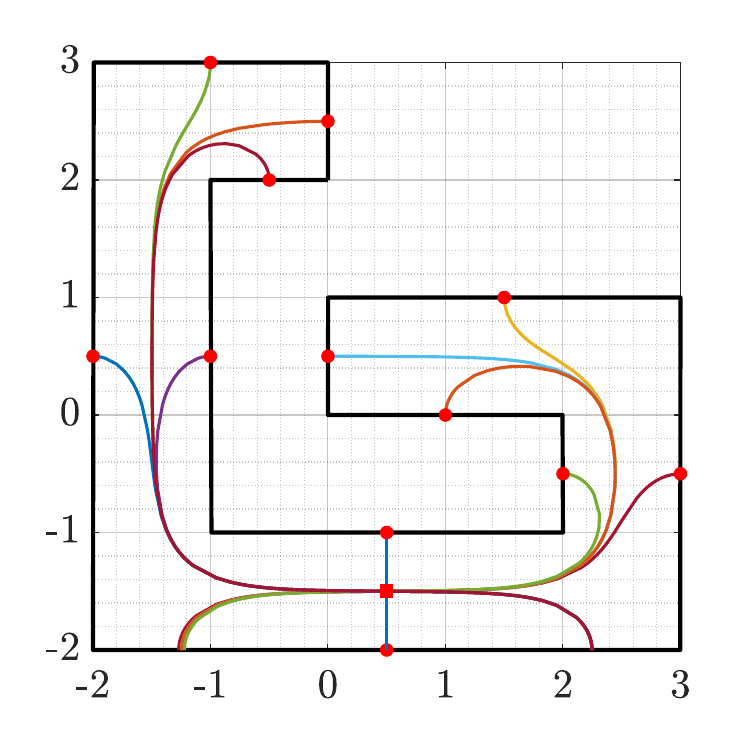}}
		\hfill
		\scalebox{0.6}{\includegraphics[trim=0cm 0.5cm 0cm 0.5cm,clip]{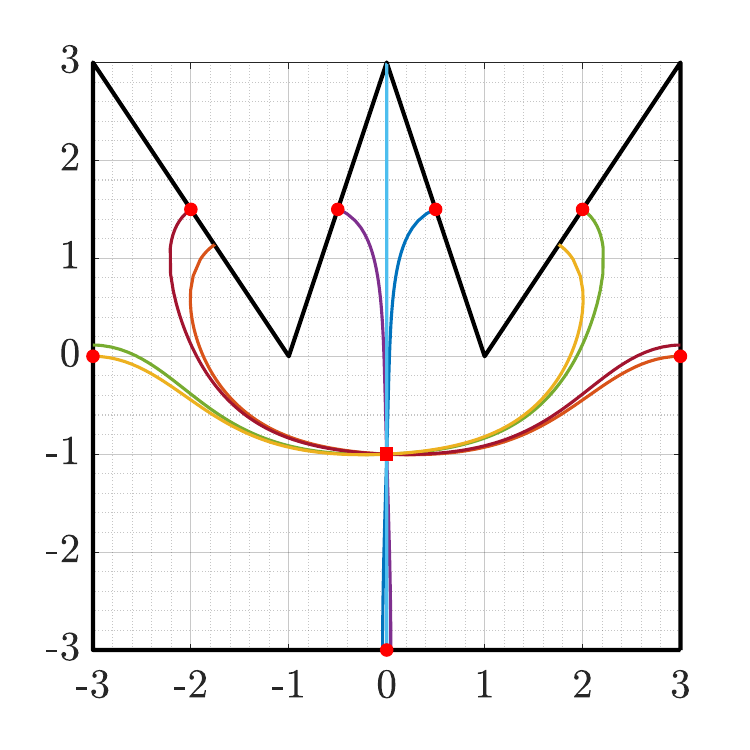}}		
	}
\caption{The hyperbolic geodesics pass through the middle of each boundary segment (marked with a red-dot) and the point $\alpha$ (marked with a red-square) for the two polygons.}
	\label{fig:hyp-geod}
\end{figure}

\section{Conclusions}\label{sec:con}

We have seen in this paper that the boundary integral equation method can be efficiently applied for numerical approximation of a great variety of conformal invariants. Comparison to the other methods we are familiar with suggests that:
\begin{itemize}
\item[•] The accuracy of the results is carefully compared in our papers to
the literature and to the results provided by $hp$-FEM computations. The conclusion is that the same accuracy is attained.
\item[•] The speed of computation in the cases tested is usually significantly faster than other methods.
\item[•] The flexibility to modify the method from case to case makes the method the preferred in many cases.
\end{itemize}

In the case of simply connected polygonal domains, the SCToolbox of Driscoll \cite{dri} seems to be most widely used. For such domains, the boundary integral method presented in this paper produces results with accuracy almost comparable to the popular SCToolbox, see e.g.~\cite{Nas-PlgCir}.

The classic book of P\'{o}lya-Szeg\"{o} \cite{ps} has inspired several generations of mathematicians to study
isoperimetric problems. This indeed is
pioneering work which has also inspired our work in the past. They produced numerous tables with numerical data at a precomputer time which was remarkable and formulated several problems which are open even now. One of these
problems was to find the fundamental frequency of the Laplace operator in polygonal plane domains, discussed
recently in \cite{ind}.

\subsection{Topics for future work} 
In our papers we have mentioned many open problems, specific to the topics discussed.
The following two topics seem to offer many challenges for later research.

As far as we know very little is published
about the concrete estimates of the principal frequency of the Laplace operator in a bounded polygonal domain.
Such an estimate is expected to depend
on the geometry of the domain \cite{ps, ind}.

The classical problems of Gr\"otzsch and Teichm\"uller discussed at the end of Section 2  could be also studied in polygonal planar domains.

\section{Topicwise guides to literature}\label{sec:topic}

Numerical methods for conformal mappings have applications to many
areas of science and technology, from physics and engineering to
computer graphics and fluid dynamics. We give here a list of literature on mathematical aspects of the topic, based on well-known
sources.

The bibliographies of  Gaier \cite{gai} (ca. 300 items),
Wegmann \cite{weg05} (ca. 300 items), Kythe \cite{ky} (ca. 1000 items),
Ivanov-Trubetskov \cite{IT} (ca. 170 items),
Papamichael-Stylianopoulos \cite{ps10} (ca. 200 items),
and the books listed below in Subsection \ref{topic04} and their bibliographies provide 
information about the literature before the year 2000.
See also the bibliography of the most recent book on the
topic by D. Crowdy \cite{c} (ca. 120 items).

\subsection{Books on complex analysis} \label{topic02}
Ahlfors \cite{ah}, Sansone-Gerretsen \cite{sg}, Ablowitz-Fokas \cite{af}, Goluzin \cite{gol},  Volkovyski\v{i},  Lunts, and  Aramanovich \cite{vla}, Nehari \cite{neh}, Tsuji \cite{tsu}, Shaw \cite{shaw}, Dovbush-Krantz \cite{dk}.

\subsection{Books on numerical conformal mappings} \label{topic08}
Gaier \cite{gai}, Henrici \cite{hen}, Driscoll-Trefethen \cite{dt}, Crowdy \cite{c}, Papamichael-Stylianopoulos \cite{ps10},  Kythe \cite{ky}, Wen \cite{wen}.

\subsection{Formulas for conformal mappings,  tables} \label{topic01}
Abramowitz-Stegun \cite{ab},   Kythe \cite{ky},
von Koppenfels-Stallmann \cite{vs}, Lavrik-Savenkov \cite{LaSa}, Dalischau \cite{da}, Schintzinger-Laura \cite{ScLa},
Ivanov-Trubetskov \cite{IT}.

\subsection{Elliptic functions and integrals} \label{topic03}
Akhiezer \cite{akh}, Lawden \cite{law},  Anderson-Vamanamurthy-Vuorinen \cite{avv}, Borwein-Borwein \cite{bb}.

\subsection{Collections of papers, computational complex analysis} \label{topic04}
Beckenbach \cite{beck}, Todd \cite{to}, Trefethen \cite{tr}, Papamichael-Saff \cite{pasa},  K\"uhnau \cite{kuhHB, kuhbib}.

\subsection{Canonical domains} \label{topic05}
Koebe \cite{koe}, K\"uhnau \cite{kuh}, Tsuji \cite{tsu}, Nehari \cite{neh}, Wen \cite{wen}, Bergman \cite{ber}.

\subsection{Potential theory, capacity} \label{topic06}

\begin{itemize}
\item[-] Popular surveys: Hardin-Saff \cite{hs}.
\item[-] Capacity of condensers: P\'{o}lya-Szeg\"{o} \cite{ps}, Landkof \cite{lan},   Ransford \cite{r}, Tsuji \cite{tsu}, Garnett-Marshall \cite{garmar},
Gol'dshte\v{i}n-Reshetnyak \cite{gor}, Heinonen-Kilpel\"ainen-Martio \cite{hkm}, Kirsch \cite{ki}, Ohtsuka \cite{o}.
\item[-] Generalized condenser: Dubinin \cite{du}.
\item[-] Logarithmic capacity:  Saff-Totik \cite{st},  Borodachov,  Hardin, and  Saff \cite{bhs}, Liesen-S\`ete-Nasser \cite{LSN17, LSN23}, Ransford-Rostand \cite{rr}.
\item[-] Analytic capacity:  Ransford-Younsi \cite{yr},   Nasser-Green-Vuorinen \cite{ngv}.  
\end{itemize}

\subsection{Mathematical analysis, tools for constructive complex analysis} \label{topic07} Atkinson \cite{atk},
Gakhov \cite{Gak}, Mikhlin \cite{Mik64}, Kress \cite{kre14}, Wen \cite{wen},
Crowdy \cite{c}, Baernstein \cite{bae},  Trefethen  \cite{Tr2,Tr3}.

\subsection{Surveys and comparisons of methods for numerical conformal mappings} \label{topic09} Gaier \cite{gai},
Papamichael-Stylianopoulos \cite{ps10}, \cite[pp.13-16]{ps}, Papamichael \cite{p1},
Trefethen \cite{Tr2,Tr3,tre5}, 
Trefethen-Driscoll  \cite{TrD},   Kythe \cite[pp.3-12]{ky}, 
Papamichael \cite{p3},
 Porter \cite{por07,por},  Wegmann \cite{weg05}, Gutknecht  \cite{gut},
Driscoll-Trefethen \cite{dt}, 
Badreddine,  DeLillo,  and  Sahraei \cite{bds}, Ivanov-Trubetskov \cite[pp.1-4]{IT}.
\begin{itemize}
\item[-] Popular surveys: Bishop \cite{bis2}, Crowdy \cite{cr08}, Porter \cite{por07}, Trefethen \cite{tre5},
Trefethen-Driscoll  \cite{TrD}, Gu, Luo, and Yau \cite{gly}.
\item[-] Koebe method: Nasser et al. \cite{Nas-cmft15, lmny}.
\item[-]
SC methods: Driscoll-Trefethen \cite{dt},
DeLillo, Elcrat,  Kropf, and  Pfaltzgraff \cite{dekp}.
\item[-]
Circular domain: Badreddine, DeLillo and Sahraei~\cite{bds}, Benchama et. al~\cite{bdhw}, DeLillo, Horn and Pfaltzgraff~\cite{dhp}, Nasser~\cite{Nas-cmft15,Nas-PlgCir}, Wegmann~\cite{weg01}. 
\item[-]
Integral equation methods: Wegmann-Nasser \cite{wn}, Nasser~\cite{Nas-Siam1,Nas-jmaa11,Nas-jmaa13}, Nasser and Al-Shihri~\cite{Nas-Siam2}, Kerzman-Stein \cite{ks}, Kerzman-Trummer \cite{kt}, Henrici \cite[p. 560]{hen}, Bell \cite{bell}, Abzalilov-Ivanshin-Shirokova~\cite{AbISh}, Abzalilov-Shirokova~\cite{AbSh}, Razali-Nashed-Murid~\cite{RNMu}, Murid-Nashed-Razali~\cite{MuNR}.
\item[-]
Charge simulation method: Amano et al. \cite{ama,aoos},
Liesen, S\`ete, and Nasser \cite{LSN17}.
\item[-]
Zipper method: K\"uhnau \cite{kuh83}, Marshall-Rohde \cite{mr07}.
\item[-]
Conjugate function method: Hakula-Rasila \cite{hr}, Hakula-Rasila-Zheng \cite{hrz}.
\item[-]
FEM methods (AFEM, $hp$-FEM): Samuelsson-Vuorinen \cite{sv}, Betsakos-Samuelsson-Vuorinen \cite{bsv}, Hakula-Rasila-Vuorinen 
\cite{hrv1,hrv2,hrv3}, Hakula-Rasila  \cite{hr}.
\item[-] Complexity of numerical conformal mapping: Bishop \cite{bis1}.
\item[-] Circle packing: Stephenson \cite{step}
\end{itemize}

\subsection{Conformal invariants, applications to geometric function theory} \label{topic11}
Ahlfors \cite{ah}, Bishop \cite{bis2}, Garnett-Marshall \cite{garmar}, Lehto-Virtanen \cite{lv},
Anderson-Vamanamurthy-Vuorinen \cite{avv}, Hariri-Kl\'en-Vuorinen \cite{hkv}, Kuz\'{}mina \cite{kuz}, Dubinin \cite{du}, V\"ais\"al\"a \cite{vais},  Gehring-Martin-Palka
\cite{gmp},  Avkhadiev \cite{avk}, Avkhadiev-Kayumov-Nasyrov \cite{akn}.

\subsection{Applications (conformal mappings on Riemann surfaces, image processing)} \label{topic12} 
 Kropf-Yin-Yau-Gu \cite{kyyg},
Choi \cite{choi}, Hakula-Rasila \cite{hr}, Hakula-Rasila-Zheng \cite{hrz},
Rainio-Nasser-Vuorinen-Kl\'en \cite{rnvk}.

\subsection{Software for computational complex analysis} \label{topic13}
Driscoll \cite{dri}, Greengard-Rokhlin \cite{gro, gg}, Nasser \cite{Nas-PlgCir}, Trefethen   \cite{Tr2,Tr3}.

\vspace{1cm}


{\bf Acknowledgements.} The second
author has very pleasant memories about the memorable hospitality of Prof. Vladimir Gol'dshte\v{i}n and his wife Olga during his visits to Sobolev Institute of Mathematics, Academgorodok in 1981 and 1984. He
is also indebted to the first author
for the excellent arrangements for
his research visit to Wichita State University in 2025.

\medskip

\medskip

\noindent{\bf Funding:} No funding.

\noindent{\bf Data availability:}  The MATLAB codes for the figures and numerical results presented in this paper are available at \url{https://github.com/mmsnasser/cap}.

\medskip

\noindent{\bf Declarations}\\
\noindent{\bf Ethics approval and consent to participate:} Not applicable

\noindent{\bf Consent for publication:} Not required

\noindent{\bf Conflict of interest:} The author declares no competing interests.

\noindent{\bf Authors Contribution:} The authors contributed equally.


\medskip


\begin{thebibliography}{333}


\bibitem{af}{\small \textsc{M. J. Ablowitz and A.S. Fokas,} \emph{Complex variables: 
introduction and applications.} Second edition. Cambridge Texts in Applied Mathematics. 
Cambridge University Press, Cambridge, 2003. xii+647 pp.}

%
\bibitem{as}{\small \textsc{ M.  Abramowitz and I. A. Stegun}, \emph{
Handbook of mathematical functions with formulas, graphs, and mathematical tables.} Reprint of the 1972 edition. Dover Publications, Inc., New York, 1992. xiv+1046 pp.} 
%


\bibitem{AbISh}
{\small \textsc{D.F. Abzalilov, P. N. Ivanshin and E.A. Shirokova,} Solution the Dirichlet problem for multiply connected domain using numerical conformal mapping, Complex Anal. Oper. Theory 13 (2019) 1419--1429.}

\bibitem{AbSh}
{\small \textsc{D.F. Abzalilov and E. A. Shirokova,} The Approximate conformal mapping onto multiply connected domains, Probl. Anal. Issues Anal. 8 (26) (2019) 3--16.}

\bibitem{ah}{\small \textsc{L.V. Ahlfors,} \emph{ Conformal invariants: topics 
in geometric function theory.} McGraw-Hill Series in Higher Mathematics.
 McGraw-Hill Book Co.,  1973. ix+157 pp. }


\bibitem{ab}
{\small \textsc{
L. Ahlfors and A. Beurling, } Conformal invariants and function-theoretic null-sets. Acta Math. 83 (1950), 101--129.}



\bibitem{akh}{\small \textsc{N.I.  Akhiezer,}  \emph{ Elements of the theory of elliptic functions.} Translated from the second Russian edition by H. H. McFaden. Translations of Mathematical Monographs, 79. American Mathematical Society, Providence, RI, 1990. viii+237 pp.} 

\bibitem{ama}{\small \textsc{
K. Amano,} A Charge Simulation Method for the Numerical Conformal
Mapping of Interior, Exterior and Doubly Connected Domains. Journal of
Computational and Applied Mathematics. (1994) 53: 353--370.}


\bibitem{aoos}{\small \textsc{K.
Amano, D. Okano, H. Ogata, and M. Sugihara,} Numerical conformal mappings onto the linear slit domain. Jpn. J. Ind. Appl. Math. 29 (2012), no. 2, 165--186.}
\bibitem{avv}{\small \textsc{G.D. Anderson,  M.K.
Vamanamurthy,   and M.K.  Vuorinen, }  \emph{Conformal invariants, inequalities, and quasiconformal maps. } With 1 IBM-PC floppy disk (3.5 inch; HD). Canadian Mathematical Society Series of Monographs and Advanced Texts. A Wiley-Interscience Publication. John Wiley  \& Sons, Inc., New York, 1997. xxviii+505 pp.
}


\bibitem{atk}{\small \textsc{K.E. Atkinson,} \emph{The numerical solution of integral
equations of the second kind.} Cambridge University Press, Cambridge, 1997.}


\bibitem{avk}{\small \textsc{
F.G.~Avkhadiev,} {Conformally invariant inequalities}.  (Russian), Kazan 2020.}


\bibitem{akn}{\small \textsc{F.G. Avkhadiev, I.R. Kayumov, and S.R. Nasyrov,} Extremal problems in geometric function theory. (Russian) Uspekhi Mat. Nauk 78 (2023), no. 2(470), 3--70; translation in Russian Math. Surveys 78 (2023), no. 2, 211--271.}



\bibitem{bds}{\small \textsc{M.
Badreddine, T.K. DeLillo,  and S. Sahraei,} A comparison of some numerical conformal 
mapping methods for simply and multiply connected domains. 
Discrete Contin. Dyn. Syst. Ser. B 24 (2019), no. 1, 55--82.}


\bibitem{bae}{\small \textsc{A.
Baernstein,II}, \emph{Symmetrization in analysis.} With David Drasin and Richard S. Laugesen. 
With a foreword by Walter Hayman. New Mathematical Monographs, 36. 
Cambridge University Press, Cambridge, 2019. xviii+473 pp.}

	\bibitem{b}{\small \textsc{A.\,F. Beardon},  \emph{ The Geometry of Discrete Groups,}
		Springer-Verlag, New York, 1983.}



	


\bibitem{beck} {\small \textsc{E.F. Beckenbach, ed.,} \emph{
Construction and applications of conformal maps.} Proceedings of a symposium.
National Bureau of Standards Applied Mathematics Series, No. 18. 
U.S. Government Printing Office, Washington, DC, 1952.}




\bibitem{bell} {\small \textsc{S. Bell,}
\emph{The Cauchy transform, potential theory and conformal mapping. }Second edition. Chapman \& Hall/CRC, Boca Raton, FL, 2016. xii+209 pp. 
 }


\bibitem{bdhw}{\small \textsc{N. Benchama, T.K. DeLillo, T. Hrycak and L. Wang,} A simplified Fornberg-like method for the conformal mapping of multiply connected regions—comparisons and crowding. 
J. Comput. Appl. Math. 209 (2007) 1--21.}




\bibitem{ber}{\small \textsc{S. Bergman, } \emph{The Kernel Function and 
Conformal Mapping.} AMS, Providence, 1970.}



\bibitem{bsv}{\small \textsc{D. Betsakos, K. Samuelsson,  and M. Vuorinen,} 
The computation of capacity of planar condensers. 
Publ. Inst. Math. (Beograd) (N.S.) 75(89) (2004), 233--252.}
%
\bibitem{bis1}{\small \textsc{Ch. J.
Bishop, } Conformal mapping in linear time. Discrete Comput. Geom. 44 (2010), no. 2, 330--428.}
%
\bibitem{bis2}{\small \textsc{Ch. J.
Bishop,} Harmonic measure: algorithms and applications. Proceedings of the International Congress of Mathematicians--Rio de Janeiro 2018. Vol. III. Invited lectures, 1511--1537, World Sci. Publ., Hackensack, NJ, 2018. }




\bibitem{bhs}{\small \textsc{S.V. Borodachov, D.P. Hardin, and E.B. Saff,}  \emph{Discrete energy on rectifiable sets.  }Springer Monographs in Mathematics. Springer, New York, [2019], 2019. xviii+666 pp.}


\bibitem{bb}{\small \textsc{J.M. Borwein and P.B. Borwein,}
 \emph{Pi and the AGM.
A study in analytic number theory and computational complexity. } Canadian Mathematical Society Series of Monographs and Advanced Texts. A Wiley-Interscience Publication. John Wiley  \& Sons, Inc., New York, 1987. xvi+414 pp.}




\bibitem{choi} {\small \textsc{G.P.T. Choi,} Efficient Conformal Parameterization of 
Multiply-Connected Surfaces Using Quasi-Conformal Theory. J. Sci. Comput. 87 (2021), 70.}


\bibitem{fft} 
{\small \textsc{J.W. Cooley and J.W. Tukey,} 
An algorithm for the machine calculation of complex Fourier series. 
Math. Comput.  19 (1965), 297--301.}
%
\bibitem{cr08}{\small \textsc{D.
Crowdy,} Geometric function theory: a modern view of a classical subject. Nonlinearity 21 (2008), no. 10, T205--T219.}


\bibitem{c-slit} {\small \textsc{D. Crowdy,} Conformal slit maps in applied mathematics. 
ANZIAM J.  53 (2012), 171--189.}

\bibitem{c}{\small \textsc{D. Crowdy,} \emph{ Solving problems in multiply 
connected domains.} SIAM, Philadelphia, PA, 2020.}

\bibitem{CKGN} {\small \textsc{D.G. Crowdy, E.H. Kropf, C.C. Green, and M.M.S. Nasser,} 
The Schottky-Klein prime function: a theoretical and computational tool for applications. 
IMA J. Appl. Math. 81 (2016) 589--628.}


\bibitem{da}{\small \textsc{H. Dalichau,} Conformal mapping and elliptic functions, M\"unchen, 2006  \hfill \newline
{\tt http://dateiena.harald-dalichau.de/spcm/pref11.pdf}}
%



\bibitem{dekp}{\small \textsc{T.K.
DeLillo, A.R. Elcrat,  E.H. Kropf, and J.A. Pfaltzgraff,}
Efficient calculation of Schwarz-Christoffel transformations for multiply connected domains using Laurent series. 
Comput. Methods Funct. Theory 13 (2013), no. 2, 307--336.}





\bibitem{dhp}{\small \textsc{T.K. DeLillo, M.A. Horn, and J.A. Pfaltzgraff,}
Numerical conformal mapping of multiply connected regions by Fornberg-like methods.
Numer. Math. 83 (1999) 205--230.}



\bibitem{dk}{\small \textsc{P.V.
Dovbush and S.G. Krantz,} \emph{The geometric theory of complex variables.} Springer, Cham, 2025. xi+540 pp.} 

\bibitem{dri}{\small \textsc{T.A. Driscoll,} Schwarz--Christoffel Toolbox for MATLAB, 
\url{https://tobydriscoll.net/project/sc-toolbox/}. Accessed 11 May 2021.}


\bibitem{dt}{\small \textsc{T.A. Driscoll and L.N. Trefethen,} \emph{Schwarz-Christoffel mapping.} 
Cambridge Monographs on Applied and Computational Mathematics, 8. 
Cambridge University Press, Cambridge, 2002. xvi+132 pp.}



\bibitem{du}{\small \textsc{V.N.  Dubinin, }  \emph{Condenser Capacities and 
Symmetrization in Geometric Function Theory,} Birkh\"auser, 2014.}


\bibitem{DuKu}{\small \textsc{P. Duren and R. K{\"u}hnau,} 
Elliptic capacity and its distortion under conformal mapping. 
J. Anal. Math. 89(1) (2003), 317-335. }


\bibitem{dp}{\small \textsc{P.
Duren and J. Pfaltzgraff,} Robin capacity and extremal length. J. Math. Anal. Appl. 179 (1993), no. 1, 110--119. }





\bibitem{f}{\small \textsc{B. Fuglede,} Extremal length and functional completion. Acta Math. 98 (1957), 171--219.}








%
\bibitem{gai}{\small \textsc{D. Gaier,} \emph{
Konstruktive Methoden der konformen Abbildung.} (German)
Springer Tracts in Natural Philosophy, Vol. 3. Springer--Verlag, Berlin, 1964. xiii+294 pp.}


\bibitem{Ga}{\small \textsc{D. Gaier,} 
Ermittlung des konformen Moduls von Vierecken mit Differenzenmethoden. 
Numer. Math. 19 (1972), 179--194.}

\bibitem{Gak}{\small \textsc{F.D. Gakhov,}
\emph{Boundary value problems.} Pergamon Press, Oxford, 1966.}


\bibitem{garmar}{\small \textsc{J.B. Garnett and D.E. Marshall,}
\emph{Harmonic measure.} Cambridge University Press, Cambridge, 2008.}



\bibitem{gmp}{\small \textsc{F.W. Gehring, G.J. Martin, and B.P. Palka,}   
\emph{ An introduction to the theory of higher-dimensional quasiconformal mappings.} 
Mathematical Surveys and Monographs, 216. American Mathematical Society, Providence, RI, 2017. ix+430 pp.}
\bibitem{gor}{\small \textsc{V.M.
Gol'dshte\v{i}n and Yu.G. Reshetnyak,} \emph{ Quasiconformal mappings and Sobolev spaces.} Kluwer Academic Publishers Group, Dordrecht, 1990.  }
\bibitem{gol}{\small \textsc{G.M. Goluzin,}
\emph{Geometric Theory of Functions of a Complex Variable.} Amer. Math. Soc., Providence, 1969.}


\bibitem{GrNa24}
{\small \textsc{C.C. Green and M.M.S. Nasser, }
 Towards computing the harmonic-measure distribution function for the middle-thirds Cantor set. 
J. Comput. Appl. Math. 448 (2024) 115903.}

\bibitem{GrNag}
{\small \textsc{C.C. Green and M.M.S. Nasser, }
Numerical computation of Stephenson's $g$-functions in multiply connected domains.
J. Math. Anal. Appl. 554 (2026) 130010.}



\bibitem{GSWC}
{\small \textsc{C.C. Green, M.A. Snipes, L.A. Ward and D.G. Crowdy, }
Harmonic-measure distribution functions for a class of multiply connected symmetrical slit domains. 
Proc. R. Soc. A 478 (2259) (2022) 20210832.}





\bibitem{gg}{\small \textsc{L. Greengard and Z. Gimbutas,} FMMLIB2D: A {MATLAB} toolbox for fast multipole method in two dimensions, version 1.2. 2019, 
\url{www.cims.nyu.edu/cmcl/fmm2dlib/fmm2dlib.html}. Accessed 6 Nov 2020.}




\bibitem{gro}
{\small \textsc{L. Greengard and V. Rokhlin, }
 A fast algorithm for particle simulations. J. Comput. Phys. 73 (1987), no. 2, 325--348.}
 
 

\bibitem{gly}{\small \textsc{ X.
Gu, F. Luo,  and S.-T. Yau,} Computational conformal geometry behind modern technologies. Notices Amer. Math. Soc. 67 (2020), no. 10, 1509--1525.} 


\bibitem{gut}{\small \textsc{M.H.
Gutknecht,}  Numerical conformal mapping methods based on function conjugation. Special issue on numerical conformal mapping. J. Comput. Appl. Math. 14 (1986), no. 1--2, 31--77.}


\bibitem{hnv0}  {\small \textsc{H. Hakula, M.M.S. Nasser, and  M. Vuorinen,} 
{Conformal capacity and polycircular domains}.
\href{https://doi.org/10.1016/j.cam.2022.114802}{J. Comput. Appl. Math. 420 (2023), 114802}, arXiv:2202.12922.}


 
\bibitem{hnv1}  {\small \textsc{H. Hakula, M.M.S. Nasser, and  M. Vuorinen,} 
{Mobile disks in hyperbolic space and minimization of conformal capacity}.
Electron. Trans. Numer. Anal.  60, (2024),   1--19, 
arXiv:2303.00145.}


 \bibitem{hnv2}{ \small \textsc{H. Hakula, M.M.S. Nasser, and  M. Vuorinen,}
{Constrained maximization of conformal capacity}. 
 \href{https://doi.org/10.1016/j.camwa.2025.05.021}{Comput.  Math. Appl.  193 (2025), 132--156}, arXiv:2404.19663.}


 \bibitem{hr}{ \small \textsc{H. Hakula and A. Rasila,} Laplace-Beltrami equations and numerical conformal mappings on surfaces. SIAM J. Sci. Comput. 47 (2025), no. 1, A325--A342.}
 

 
\bibitem{hrv1}{\small \textsc{H. Hakula, A. Rasila, and M. Vuorinen,} 
On moduli of rings and quadrilaterals: algorithms and experiments.
 SIAM J. Sci. Comput. 33 (2011), no. 1, 279--302.}

\bibitem{hrv2}{\small \textsc{H. Hakula,  A. Rasila,  and M. Vuorinen,} 
Computation of exterior moduli of quadrilaterals. 
Electron. Trans. Numer. Anal. 40, 1--16, 2013, ISSN 1068-9613.}

\bibitem{hrv3}{\small \textsc{H. Hakula,  A. Rasila,  and M. Vuorinen,} 
{Conformal modulus on domains with strong singularities and cusps}.
{ Electron. Trans. Numer. Anal. 48, 462--478, 2018,} DOI: 10.1553/etna\_vol48s462.
{ arXiv:1501.06765.}
}


 \bibitem{hrz}{ \small \textsc{H. Hakula, A. Rasila, and Y. Zheng,}
The Conjugate Function Method for Surfaces with Elaborate Topological Types arXiv:2509.01978.}


\bibitem{hrs}{\small \textsc{Q. Han,  A. Rasila,  and T. Sottinen,} Efficient simulation of mixed boundary value problems and conformal mappings. Appl. Math. Comput. 488 (2025), Paper No. 129119, 14 pp.}

\bibitem{hs}{\small \textsc{D.P.
Hardin and E.B. Saff,} Discretizing manifolds via minimum energy points. Notices Amer. Math. Soc. 51 (2004), no. 10, 1186--1194.}

\bibitem{hkv}{\small \textsc{P.~Hariri, R.~Kl\'en, and M.~Vuorinen,} 
 \emph{Conformally {I}nvariant {M}etrics  and {Q}uasiconformal {M}appings}, 
Springer Monographs in Mathematics, Springer, Berlin, 2020.}

\bibitem{hvv} {\small \textsc{V. Heikkala, M. K. Vamanamurthy, and  M. Vuorinen},
{ Generalized elliptic integrals.} Comput. Methods Funct. Theory
9(2009), 75-109,
\href{http://arxiv.org/abs/math/0701436}{arXiv :0701436}.
}

\bibitem{hkm} {\small \textsc{J. Heinonen, T. Kilpel\"ainen, and O. Martio,} \emph{Nonlinear potential theory of degenerate elliptic equations.} Unabridged republication of the 1993 original. Dover Publications, Inc., Mineola, NY, 2006. xii+404 pp.}

\bibitem{hen}{\small \textsc{P. Henrici,}  \emph{Applied and Computational 
Complex Analysis}, Vol. 3, John Wiley \& Sons, New York, 1986.}


%
\bibitem{ind}{\small \textsc{E.
Indrei,} On the first eigenvalue of the Laplacian for polygons. J. Math. Phys. 65 (2024), no. 4, Paper No. 041506, 40 pp.}

\bibitem{IT}{\small  \textsc{V.I. Ivanov and M. K. Trubetskov,}  \emph{
Handbook of conformal mapping with computer-aided visualization.}
With 1 IBM-PC floppy disk (5.25 inch; HD). CRC Press, Boca Raton, FL, 1995. vi+360 pp. ISBN: 0--8493--8936--4.}

\bibitem{jen}{\small \textsc{J.A. Jenkins,}
\emph{Univalent Functions and Conformal Mapping.} Springer, Berlin, 1958.}

\bibitem{kak} {\small \textsc{S. Kakutani, }
Two-dimensional Brownian motion and harmonic functions.
Proc. Imp. Acad. Tokyo 20 (1944) 706--714.}

\bibitem{knv}{\small \textsc{E. Kalmoun, M.M.S. Nasser, and M. Vuorinen,} 
Numerical computation of Mityuk's function and radius for circular$/$radial slit domains. 
J. Math. Anal. Appl. 490 (2020) 124328.}



\bibitem{KS06} {\small \textsc{S. Kanas and T. Sugawa, }
On conformal representations of the interior of an ellipse.
Ann. Acad. Sci. Fenn. Math. 31 (2006), 329--348.}

\bibitem{krv}{\small \textsc{R. Kargar, O. Rainio, and M. Vuorinen,} Landen transformations applied to approximation. Pure Appl. Funct. Anal. 9 (2024), no. 2, 503--517.}


\bibitem{kela}{\small \textsc{L. Keen and N.  Lakic, }
\emph{Hyperbolic geometry from a local viewpoint}. London Mathematical Society Student Texts, 68. Cambridge University Press, Cambridge, 2007.}


\bibitem{ks}{\small \textsc{N. Kerzman and E.M.
 Stein,} The Cauchy kernel, the Szeg\"o kernel, and the Riemann mapping function. Math. Ann. 236 (1978), no. 1, 85--93.}
\bibitem{kt}{\small \textsc{N. Kerzman and M.R. Trummer,} Numerical conformal mapping via the Szeg\"o kernel. Special issue on numerical conformal mapping. J. Comput. Appl. Math. 14 (1986), no. 1--2, 111--123.}
\bibitem{ki}{\small \textsc{S.  Kirsch, }
Transfinite diameter, Chebyshev constant and capacity. (English summary) Handbook of complex analysis: geometric function theory. 
Vol. 2, 243--308, Elsevier Sci. B. V., Amsterdam, 2005. }

%


\bibitem{kre90}{\small \textsc{R. Kress,}  A Nystr\"om method for boundary 
integral equations in domains with corners. Numer. Math. 58(2) (1990), 145--161.}


\bibitem{kre91}{\small \textsc{R. Kress,} Boundary integral equations in 
time--harmonic acoustic scattering. Math. Comput. Modelling 15 (1991), 229--243.}


\bibitem{kre14} {\small \textsc{R. Kress,} \emph{Linear integral equations.} 
Third edition. Applied Mathematical Sciences, 82. Springer, New York, 2014. xvi+412 pp.}



\bibitem{koe10}{\small \textsc{P. Koebe,} \"Uber die Uniformisierung der algebraischen Kurven. II. Math. Ann. 69 (1910), 1--81.}

\bibitem{koe}{\small \textsc{P. Koebe,} Abhandlungen zur Theorie der konformen Abbildung, IV. 
Abbildung mehrfach zusammenh\"angender schlichter Bereiche auf 
Schlitz\-bereiche. Acta Math. 41 (1918), 305--344.}


\bibitem{kyyg} {\small \textsc{E. Kropf, X. Yin,   S.-T. Yau, and X. Gu,}
Conformal parameterization for multiply connected domains: combining finite elements 
and complex analysis, Engineering with Computers (2014) 30(4):441--455, DOI: 10.1007/s00366-013-0348-4.}





\bibitem{kuh83}{\small \textsc{R.
K\"uhnau,} Numerische Realisierung konformer Abbildungen durch "Interpolation''. (German) [Numerical realization of conformal mappings by "interpolation''] Z. Angew. Math. Mech. 63 (1983), no. 12, 631-637.}



\bibitem{kuhHB}{\small \textsc{R. K\"uhnau, ed., }
\emph{Handbook of complex analysis: geometric function theory. Vol. 1--2.}
 Elsevier Science B.V., Amsterdam, 2002-2005.}



\bibitem{kuh}{\small \textsc{R. K\"uhnau, } {
Canonical conformal and quasiconformal mappings. Identities. Kernel functions.} 
Handbook of complex analysis: geometric function theory. Vol. 2, 131--163, 
Elsevier Sci. B. V., Amsterdam, 2005.}


\bibitem{kuhbib}{\small \textsc{R. K\"uhnau,  } Bibliography of geometric function theory. 
Handbook of complex analysis: geometric function theory. Vol. 2, 809--828,
Elsevier, 2005. }

\bibitem{kuz}{\small \textsc{G.V.
 Kuz\'{}mina, } Moduli of families of curves and quadratic differentials. A translation of Trudy Mat. Inst. Steklov. 139 (1980). Proc. Steklov Inst. Math. 1982, no. 1, vii+231 pp.}


\bibitem{ky}{\small \textsc{P. K. Kythe,} \emph{Handbook of conformal 
mappings and applications.} CRC Press, Boca Raton, FL, 2019. xxxv+906 pp.}


\bibitem{lan}{\small \textsc{N.S. Landkof,} \emph{Foundations of modern potential theory.} Translated from the Russian by A. P. Doohovskoy. Die Grundlehren der mathematischen Wissenschaften, Band 180. Springer--Verlag, New York--Heidelberg, 1972. x+424 pp.}

\bibitem{LaSa}{\small \textsc{V.I. Lavrik and V.N. Savenkov,} 
\emph{A handbook on conformal mappings, } (Russian), "Naukova Dumka'', Kiev, 1970. 252 pp.}
\bibitem{law}{\small \textsc{D.F.
Lawden,  } \emph{Elliptic functions and applications.  } Applied Mathematical Sciences, 80. Springer-Verlag, New York, 1989. xiv+334 pp.}

\bibitem{lmny}{\small \textsc{
K.W. Lee, A.H.M. Murid, M.M.S. Nasser, and S.H. Yeak,} Fast implementation of generalized Koebe's iterative method, Mathematics, 13 (2025) 1920.}

\bibitem{lv}{\small \textsc{O. Lehto and K.I. Virtanen,} \emph{ Quasiconformal mappings in the plane.} Second edition. Translated from the German by K. W. Lucas. Die Grundlehren der mathematischen Wissenschaften, Band 126. Springer--Verlag, New York--Heidelberg, 1973. viii+258 pp.}



\bibitem{LSN17} {\small  \textsc{J. Liesen, O. S\`ete, and M.M.S. Nasser,} 
Fast and accurate computation of the logarithmic capacity of compact sets. 
Comput. Methods Funct. Theory  17 (2017), 689--713.}




\bibitem{LSN23} {\small  \textsc{
J. Liesen, M.M.S. Nasser, and O.S\`ete,}  Computing the logarithmic capacity of compact sets having (infinitely) many components with the charge simulation method. Numer. Algorithms 93 (2023), no. 2, 581--614.}


\bibitem{AMth} {\small  \textsc{A. Mahenthiram,} 
\emph{Harmonic-Measure Distribution Functions, and Related Functions, for Simply Connected and Multiply Connected Two-Dimensional Regions (Ph.D. thesis).} 
University of South Australia, 2024.}



\bibitem{AMta} {\small  \textsc{A. Mahenthiram,} 
Computing -functions of Some Planar Simply Connected Two-dimensional Regions.
Taiwanese J. Math. 27 (2003) 931--952.}

\bibitem{mar} {\small  \textsc{D.E. Marshall,} 
Conformal welding for finitely connected regions.
Comput. Methods Funct. Theory 11 (2011) 655--669.}

\bibitem{mr07}{\small \textsc{D.E.
Marshall and S. Rohde, } Convergence of a variant of the zipper algorithm for conformal mapping. SIAM J. Numer. Anal. 45 (2007), no. 6, 2577--2609.}


\bibitem{Mik64}{\small \textsc{S.G. Mikhlin,} \emph{Integral equations and 
their applications to certain problems in mechanics, mathematical physics and technolog.} 
2nd ed., Pergamon Press, Oxford, 1964.}



\bibitem{Mit}
{\small \textsc{I.P. Mityuk,} A generalized reduced module and some of its applications. Izv. Vyssh. Uchebn. Zaved., Mat. 2 (1964) 110--119.} [in Russian].


\bibitem{mc} {\small \textsc{H. Miyoshi and D.G. Crowdy,} Estimating conformal capacity using asymptotic matching. IMA J. Appl. Math. 88 (2023), no. 3, 472--497.}


\bibitem{MuNR}{\small \textsc{A.H.M. Murid, M.Z. Nashed and M.R.M. Razali,} 
Numerical conformal mapping for exterior regions via the Kerzman-Stein kernel, J. Integral Equations Appl. 10 (1998) 517--532.}


\bibitem{mn03}{\small \textsc{A.H.M. Murid and M.M.S. Nasser,} 
Eigenproblem of the generalized Neumann kernel. 
Bull. Malaysian Math. Sc. Soc. 26 (2003) 13--33.}


\bibitem{Nas-Siam1}{\small \textsc{M.M.S. Nasser, } Numerical conformal 
mapping via a boundary integral equation with the generalized Neumann kernel. 
SIAM J. Sci. Comput. 31 (2009), 1695--1715.}

\bibitem{Nas-jmaa11}{\small \textsc{M.M.S. Nasser, } Numerical conformal mapping of multiply connected regions onto the second, third and fourth categories of Koebe's canonical slit domains. 
J. Math. Anal. Appl. 382 (2011) 47--56.}

\bibitem{Nas-jmaa13}{\small \textsc{M.M.S. Nasser, } Numerical conformal mapping of multiply connected regions onto the fifth category of Koebe's canonical slit regions. 
J. Math. Anal. Appl. 398 (2013) 729--743.}


\bibitem{Nas-ETNA}{\small  \textsc{M.M.S. Nasser, } Fast solution of boundary integral 
equations with the generalized Neumann kernel. Electron. Trans.  Numer. Anal. 44 (2015), 189--229.}




\bibitem{Nas-cmft15}{\small  \textsc{M.M.S. Nasser, } Fast computation of the circular map. 
Comput. Methods Funct. Theory 15 (2015) 187--223.}


\bibitem{Nas-PlgCir}{\small  \textsc{M.M.S. Nasser, } PlgCirMap: 
A MATLAB toolbox for computing the conformal mapping from polygonal 
multiply connected domains onto circular domains. SoftwareX 11 (2020), 100464, arXiv 2019, arXiv:1911.01787.}



\bibitem{Nas-Siam2}{\small \textsc{M.M.S. Nasser and F.A.A. Al-Shihri, } A fast boundary integral equation method for conformal mapping of multiply connected regions. 
SIAM J. Sci. Comput. 35(3) (2013) A1736--A1760.}

 
\bibitem{ngv}  {\small \textsc{M.M.S. Nasser,  G. Green,  and  M. Vuorinen,}   {Fast computation of analytic capacity}.
{Comput. Methods Funct. Theory} 
 \href{https://doi.org/10.1007/s40315-024-00547-2}{https://doi.org/10.1007/s40315-024-00547-2},
 \href{http://arxiv.org/abs/2307.01808}{arXiv:2307.01808}.}


\bibitem{Nas-cr} {\small  \textsc{M.M.S. Nasser, A.H.M. Murid, and Z. Zamzamir,} 
A boundary integral method for the Riemann-Hilbert problem in domains with 
corners. Complex Var. Elliptic Equ. 53 (2008) 989--1008.}


\bibitem{nnv}{\small  \textsc{M. M.S. Nasser, S. Nasyrov,  and  M. Vuorinen,}
{ Level sets of potential functions bisecting unbounded quadrilaterals}.
Anal. Math. Phys. 12 (2022) 149.}


\bibitem{nrrvwyz} {\small
  \textsc{M.M.S. Nasser, O. Rainio, A. Rasila, M. Vuorinen, T. Wallace, H. Yu, and X. Zhang},
{Polycircular domains, numerical conformal mappings, and moduli of quadrilaterals.}  Adv. Comput. Math.  (2022) 48:58.}


\bibitem{nrv} {\small
  \textsc{M.M.S. Nasser, O. Rainio, and M. Vuorinen},
{Condenser capacity and hyperbolic perimeter.}  Comput. Math. Appl. 105 (2022) 54--74.}




\bibitem{nv1} {\small  \textsc{M.M.S. Nasser  and M.~Vuorinen,} Numerical 
computation of the capacity of generalized condensers. 
J. Comput. Appl. Math. 377 (2020) 112865.}



\bibitem{nv2} {\small  \textsc{ M. M.S. Nasser and  M. Vuorinen,} 
 {Isoperimetric properties of condenser capacity}.  
 \href{https://doi.org/10.1016/j.jmaa.2021.125050}{J. Math. Anal. Appl. 499 (2021),  125050}, arXiv:2010.09704.}


\bibitem{nv3} {\small  \textsc{M.M.S. Nasser  and M.~Vuorinen,} Conformal 
invariants in simply connected domains. 
Comput. Methods Funct. Theory 20 (2020) 747--775.}




\bibitem{nv4} {\small  \textsc{M.M.S. Nasser  and M.~Vuorinen,} Computation 
of conformal invariants. Appl. Math. Comput. 389 (2021), 125617.}




\bibitem{nsv} {\small  \textsc{S. Nasyrov, T. Sugawa, and  M. Vuorinen},
Moduli of quadrilaterals and quasiconformal reflection. J. Math. Anal. Appl. 
524 (2023), 127092.}

\bibitem{neh}{\small \textsc{Z. Nehari,}  \emph{Conformal Mapping.} Dover Publications, 
New York, 1952.}



\bibitem{o}{\small \textsc{M. Ohtsuka,} \emph{ Extremal length and precise functions.} With a preface by Fumi-Yuki Maeda. GAKUTO International Series. Mathematical Sciences and Applications, 19. Gakk\~{o}tosho Co., Ltd., Tokyo, 2003. vi+343 pp. }


\bibitem{p1}{\small \textsc{N. Papamichael,} 
 Dieter Gaier's contributions to numerical conformal mapping. Comput. Methods Funct. Theory 3 (2003), no. 1--2,  1--53.}


\bibitem{p2}{\small \textsc{N. Papamichael,} Lectures on Numerical Conformal Mapping,
Department of Mathematics and Statistics, University of Cyprus,
March 28, 2008.}


\bibitem{p3}{\small \textsc{N. Papamichael,}
Orthonormalization methods for numerical conformal mapping. Topics in modern function theory, 166--266, Ramanujan Math. Soc. Lect. Notes Ser., 19, Ramanujan Math. Soc., Mysore, 2013.}

\bibitem{pasa}{\small \textsc{N. Papamichael and E. B. Saff, eds.,}
\emph{
Computational complex analysis.} J. Comput. Appl. Math. 46 (1993), no. 1-2. Elsevier Science B.V., Amsterdam, 1993. pp. i--iv and 1--326.}

\bibitem{ps10} {\small  \textsc{N. Papamichael and N. Stylianopoulos,}
 \emph{Numerical conformal mapping: Domain decomposition and the mapping of quadrilaterals.}
  World Scientific Publishing Co. Pte. Ltd., Hackensack, NJ, 2010. xii+229 pp.}



 
\bibitem{ps}{\small \textsc{G. P\'{o}lya  and G. Szeg\"{o}, }  \emph{Isoperimetric Inequalities 
in Mathematical Physics.} Annals of Mathematics Studies, no. 27, Princeton University Press, 
Princeton, N. J., 1951. xvi+279 pp.}

%
\bibitem{por07}{\small \textsc{R.M.
Porter,} History and recent developments in techniques for numerical conformal mapping. Quasiconformal mappings and their applications, 207--238, Narosa, New Delhi, 2007.}

\bibitem{por}{\small \textsc{R.M. Porter, }  On the art of calculating accessory 
parameters of conformal mappings of circular arc polygons--general considerations 
and special situations. Teichm\"uller theory and moduli problem, 549--576, 
Ramanujan Math. Soc. Lect. Notes Ser., 10, Ramanujan Math. Soc., Mysore, 2010.}


\bibitem{PP} 
{\small  \textsc{P. Pyrih,} 
Logarithmic capacity is not subadditive -- a fine topology approach. 
Comment. Math. Univ. Carolin. 33  (1992) 67--72.}


\bibitem{rnvk} 
{\small  \textsc{ O. Rainio, M.M.S. Nasser, M. Vuorinen, 
and R. Kl\'en}, {Image augmentation with conformal mappings 
for a convolutional neural network}. 
 \href{https://doi.org/10.1007/s40314-023-02501-9}{Comp. Appl. Math. 42, 361 (2023)},
DOI: 10.1007/s40314-023-02501-9,
 \href{http://arxiv.org/abs/2212.05258}{arXiv:2212.05258}}.

\bibitem{r}{\small \textsc{T. Ransford,} \emph{Potential theory in the complex plane.} 
London Mathematical Society Student Texts, 28. Cambridge University Press, Cambridge, 1995. x+232 pp.}


%
\bibitem{rr}{\small \textsc{Th. Ransford and J. Rostand,} 
Computation of capacity. (English summary)
Math. Comp. 76 (2007), no. 259, 1499--1520.}


\bibitem{RNMu}{\small \textsc{M.R.M. Razali, M.Z. Nashed and A.H.M. Murid,} 
Numerical conformal mapping via the Bergman kernel, J. Comput. Appl. Math. 82 (1997) 333--350.}


\bibitem{st}{\small \textsc{E.B. Saff and V. Totik,} \emph{Logarithmic potentials with external fields.} Appendix B by Thomas Bloom. Grundlehren der mathematischen Wissenschaften [Fundamental Principles of Mathematical Sciences], 316. Second ed. Springer-Verlag, Berlin, 2024. xvi+594 pp.}

\bibitem {sam}{\small \textsc{K. Samuelsson}, \emph{ Adaptive algorithms for finite element methods approximating flow problems.}  Royal Institute of Technology, Stockholm, Sweden, TRITA-NA-96-04, Department of Numerical      Analysis and Computing Science,(1996).}

\bibitem {sv}{\small \textsc{K. Samuelsson and M. Vuorinen}, { Computation of capacity in 3D by means of a posteriori estimates for an adaptive FEM.}  Royal Institute of Technology, Stockholm, Sweden, TRITA-NA-9508, Preprint, Department of Numerical Analysis and Computing Science, 34 pp.}


\bibitem{sg}{\small \textsc{G. Sansone and J. Gerretsen,} 
 \emph{Lectures on the theory of functions of a complex variable. II: Geometric theory.} 
 Wolters-Noordhoff Publishing, Groningen, 1969. x+700 pp.}

\bibitem{Sch} {\small \textsc{M. Schiffer,} 
Some recent developments in the theory of conformal mapping. 
Appendix to: R. Courant, Dirichlet’s Principle, Conformal Mapping and Minimal Surfaces. Interscience, New York, 1950.}

\bibitem{ScLa}{\small \textsc{R.
Schinzinger and P.A.A. Laura,} \emph{
Conformal mapping: methods and applications.} Elsevier Science Publishers, B.V., Amsterdam, 1991. xx+581 pp.}
%


%
\bibitem{shaw} {\small \textsc{W.T.
Shaw, } \emph{Complex analysis with Mathematica®}. With 1 CD-ROM (Windows, Macintosh and UNIX). Cambridge University Press, Cambridge, 2006. xxvi+571 pp.}


\bibitem{SnWa} {\small  \textsc{M.A. Snipes and L.A. Ward,} Harmonic measure distribution functions of planar domains: A survey. 
J. Anal. 24 (2016) 293--330.}


\bibitem{solvu} {\small  \textsc{A.Yu.
Solynin and M.  Vuorinen, }
Extremal problems and symmetrization for plane ring domains. (English summary)
Trans. Amer. Math. Soc. 348 (1996), no. 10, 4095--4112.}
\bibitem{step} {\small  \textsc{K. Stephenson,} Introduction to circle packing. The theory of discrete analytic functions. Cambridge University Press, Cambridge, 2005. xii+356 pp.}




\bibitem{to} {\small \textsc{J. Todd, ed.,} \emph{
Experiments in the computation of conformal maps.}
National Bureau of Standards Applied Mathematics Series, No. 42. U.S. Government Printing Office, Washington, DC, 1955.}


\bibitem{tr} {\small  \textsc{L.N. Trefethen, ed.,}  \emph{
Numerical conformal mapping.} Reprint of J. Comput. Appl. Math. 14 (1986), no. 1--2.  North-Holland Publishing Co., Amsterdam, 1986.}






\bibitem{Tr} {\small  \textsc{L.N. Trefethen,} Numerical Conformal Mapping with Rational 
Functions. Comput. Methods Funct. Theory 20 (2020), 369--387.}


\bibitem{Tr2} {\small  \textsc{L.N. Trefethen,} Numerical analytic continuation. Jpn. J. Ind. Appl. Math. 40 (2023), no. 3, 1587--1636.}

\bibitem{Tr3} {\small  \textsc{L.N. Trefethen,} Rational approximation. Notices of AMS, 72 (2025), 78--81.
DOI:doi.org/10.1090/noti3066}.
%
\bibitem{tre5}{\small \textsc{L. N. Trefethen,}
Numerical conformal mapping,
arXiv:2507.14872 .}

\bibitem{TrD} {\small  \textsc{L.N. Trefethen and T. A.  Driscoll},
Schwarz-Christoffel mapping in the computer era. Proceedings of the International Congress of Mathematicians, Vol. III (Berlin, 1998). Doc. Math. 1998, Extra Vol. III, 533--542.}

\bibitem{Tre} {\small  \textsc{L.N. Trefethen and J.A.C. Weideman}, 
The exponentially convergent trapezoidal rule. 
SIAM Rev., 56 (2014), 385--458.}


\bibitem{tsu}{\small \textsc{M. Tsuji, }  \emph{Potential Theory in Modern 
Function Theory.} Chelsea Publishing Co., New York, 1975.}




\bibitem{vais}
{\small \textsc{J.   V\"ais\"al\"a},   \emph{
Lectures on $n$-dimensional quasiconformal mappings.}
Lecture Notes in Mathematics, Vol. 229. Springer--Verlag, Berlin-New York, 1971.   }

\bibitem{Vas02}
{\small \textsc{A. Vasil'ev},   \emph{Moduli of Families of Curves for Conformal and Quasiconformal Mappings.}
Springer-Verlag, Berlin, 2002.  }





\bibitem{vla}{\small \textsc{L.I.
Volkovyski\v{i}, G.L. Lunts, and I.G. Aramanovich,} \emph{A collection of problems on complex analysis.}
 Translated from the 1960 Russian original by J. Berry. Reprint of the 1965 English translation. 
 Dover Publications, Inc., New York, 1991. x+426 pp.}


\bibitem{vs}{\small \textsc{W. von Koppenfels and F.  Stallmann,} \emph{
Praxis der konformen Abbildung.} (German)
Die Grundlehren der mathematischen Wissenschaften, Band 100. Springer--Verlag, Berlin--G\"ottingen--Heidelberg, 1959. xiii+375 pp.}


\bibitem{WaWa} {\small \textsc{B.L. Walden and L.A. Ward,} 
Distributions of harmonic measure for planar domains. In: I. Laine, O. Martio (Eds.), Proceedings of 16th Nevanlinna Colloquium, Joensuu, Walter de Gruyter, Berlin,  289--299, 1996.}

\bibitem{walsh} {\small  \textsc{J.L. Walsh,} 
On the conformal mapping of multiply connected regions. 
Trans. Am. Math. Soc. 82 (1956) 128--146.}

\bibitem{weg01} {\small \textsc{R. Wegmann,} 
Fast conformal mapping of multiply connected regions.
J. Comput. Appl. Math. 130 (2001) 119--138.}


\bibitem{weg05} {\small \textsc{R. Wegmann,} 
Methods for numerical conformal mapping. In: R. K{\"u}hnau (ed.), Handbook of Complex Analysis: 
Geometric Function Theory, Vol. 2, Elsevier B. V.,  351--477, 2005.}


\bibitem{wmn} {\small  \textsc{R. Wegmann, A.H. Murid, and  M.M.S. Nasser,} The Riemann-Hilbert problem and the generalized Neumann kernel. J. Comput. Appl. Math. 182 (2005), no. 2, 388--415.}

\bibitem{wn} {\small  \textsc{R. Wegmann and  M.M.S.  Nasser,} The 
{R}iemann-{H}ilbert 
problem and the generalized {N}eumann kernel on multiply connected regions. 
J. Comput. Appl. Math. 214 (2008), 36--57.}


\bibitem{wen} {\small  \textsc{C.C. Wen,} \emph{Conformal mappings and boundary value problems.} Translated from the Chinese by Kuniko Weltin. Translations of Mathematical Monographs, 106. American Mathematical Society, Providence, RI, 1992.}


\bibitem{yr} {\small  \textsc{M.
Younsi and Th. Ransford,} Computation of analytic capacity and applications to the subadditivity problem. Comput. Methods Funct. Theory 13 (2013), no. 3, 337--382.}

\bibitem{zs} {\small  \textsc{D.G. Zill and P.D. Shanahan,} \emph{A first course in complex analysis with applications.} Jones and Bartlett Publishers, Sudbury, MA, 2003.}


\end{thebibliography}
\end{document}